\documentclass[11pt]{article}
\usepackage{color}
\usepackage[colorlinks=true,citecolor=red]{hyperref}
\usepackage{amsmath,amsfonts,amssymb,amsthm}
\usepackage{mathtools}
\usepackage{graphicx,subfig}
\graphicspath{{figs/}}
\usepackage[english]{babel}
\usepackage[utf8]{inputenc}
\usepackage[margin=0.90in]{geometry}
\usepackage[margin=0.2in]{caption}
\usepackage{times,fourier}
\usepackage{tikz-cd}

\makeatletter
\renewcommand\section{\@startsection {section}{1}{\z@}%
	{-2.2ex \@plus -1ex \@minus -.2ex}%
	{1ex \@plus.1ex}%
	{\normalfont\bf\sffamily}}
\renewcommand\subsection{\@startsection{subsection}{2}{\z@}%
	{-2ex\@plus -0.4ex \@minus -.2ex}%
	{0.6ex \@plus .1ex}%
	{\normalfont\small\bf\sffamily}}
\renewcommand\subsubsection{\@startsection{subsubsection}{3}{\z@}%
	{-0.6ex\@plus -0.2ex \@minus -.2ex}%
	{0.4ex \@plus .1ex}%
	{\normalfont\normalsize\it}}
\renewcommand\paragraph{\@startsection{paragraph}{4}{\z@}%
	{0.2ex \@plus0.2ex \@minus0.1ex}{-0.5em}%
	{\normalfont\normalsize\bfseries}}
\def\ps@headings{%
	\let\@oddfoot\@empty
	\let\@evenfoot\@empty
	\def\@evenhead{\small\sffamily\thepage\hfil\slshape\leftmark}%
	\def\@oddhead{\small\sffamily{\slshape\rightmark}\hfil\thepage}%
	\let\@mkboth\markboth
	\def\chaptermark##1{\markboth{{\ifnum \c@secnumdepth >\m@ne
				\if@mainmatter \@chapapp\ \thechapter. \ \fi \fi ##1}}{}}%
	\def\sectionmark##1{\markright {{\ifnum \c@secnumdepth >\z@
				\thesection. \ \fi ##1}}}}
\makeatother

\advance\textwidth -8mm
\advance\textheight -5mm
\advance\hoffset 4mm
\advance\voffset -2mm
\advance\parskip 2pt
\newdimen\figwidth
\newtheorem{theorem}{Theorem}[section]
\newtheorem{corollary}[theorem]{Corollary}
\newtheorem{lemma}[theorem]{Lemma}

\newtheorem{remark}[theorem]{Remark}
\newtheorem{definition}[theorem]{Definition}
\numberwithin{equation}{section}

\def\maketitle{\par\noindent{\LARGE\bf\sffamily\thetitle}\\[1.6ex]{\large\theauthor}\\[1.0ex]
\textit{\thetextinfo}\\[0.6ex]{\small\today}\par\vglue1.4\bigskipamount}
\def\title#1{\def\thetitle{#1}}
\def\author#1{\def\theauthor{#1}}
\def\textinfo#1{\def\thetextinfo{#1}}

\def\be{\begin{equation}}
\def\ee{\end{equation}}
\def\bse{\begin{subequations}}
\def\ese{\end{subequations}}

\definecolor{deeppurple}{rgb}{0.5, 0, 0.7}

\def\half{{\textstyle\frac12}}
\def\halfi{{\textstyle\frac{\rm i}2}}
\def\txtfrac#1#2{{\textstyle\frac{#1}{#2}}}

\def\halft{{\textstyle\frac t2}}
\def\halfA{{\textstyle\frac A2}}

\def\halfk{{\textstyle\frac m2}}

\def\sech{\mathop{\rm sech}\nolimits}
\def\diag{\mathop{\rm diag}\nolimits}
\def\dn{\mathop{\rm dn}\nolimits}
\def\sn{\mathop{\rm sn}\nolimits}
\def\cn{\mathop{\rm cn}\nolimits}
\def\am{\mathop{\rm am}\nolimits}

\def\Natural{\mathbb{N}}
\def\N{\mathbb{N}}
\def\Real{\mathbb{R}}
\def\R{\mathbb{R}}
\def\Complex{\mathbb{C}}
\def\C{\mathbb{C}}

\def\Integer{\mathbb{Z}}
\def\Z{\mathbb{Z}}
\def\i{\mathrm{i}}
\def\ep{\epsilon}
\def\x{\xi}
\def\z{\zeta}

\def\Re{\mathop{\rm Re}\nolimits}
\def\Im{\mathop{\rm Im}\nolimits}
\def\tr{\mathop{\rm tr}\nolimits}

\def\d{\mathrm{d}}
\def\e{\mathop{\rm e}\nolimits}
\def\@#1{{\mathbf{#1}}}
\def\_#1{{\mathsf{#1}}}
\let\~=\tilde
\let\==\overline
\let\^=\hat
\let\<=\langle
\let\>=\rangle

\def\max{\mathop{\rm max}\nolimits}

\def\T{{\mathbb T}}

\def\D{\Delta}
\def\l{\lambda}
\def\s{\sigma}
\def\g{\gamma}
\def\G{\Gamma}
\def\1{{\bf 1}}
\def\am{{\rm am}}

\def\T{{\rm T}}
\def\BC{\mathrm{BC}}
\def\ra{\rightarrow}
\makeatother

\let\trueparagraph=\paragraph
\def\paragraph#1{\par\smallskip\trueparagraph{\rm\textbf{#1}}}

\allowdisplaybreaks

\def\bb{\begingroup\color{blue}}

\def\eb{\endgroup}
\def\note[#1]{\marginpar{\color{red}[#1]}}

\begin{document}
\pagestyle{plain}
	
\title{Elliptic finite-band potentials of a non-self-adjoint\\[0.2ex]Dirac operator}
\author{Gino Biondini$^1$, Xu-Dan Luo$^2$, Jeffrey Oregero$^3$, Alexander Tovbis$^3$}
\textinfo{{\small1:} Department of Mathematics, State University of New York at Buffalo, Buffalo, New York 14260\\
{\small2:} Key Laboratory of Mathematics Mechanization, Academy of Mathematics and Systems Science, Chinese Academy of Sciences, Beijing 100190, China\\
{\small3:} Department of Mathematics, University of Central Florida, Orlando, Florida 32816}
\maketitle
	
\begin{abstract}
We present an explicit two-parameter family of finite-band Jacobi elliptic potentials given by
$q\equiv A\dn(x;m)$, where $m\in(0,1)$ and $A$ can be taken to be positive without loss of generality,  
for a non-self-adjoint Dirac operator $L$, 
which connects two well-known limiting cases of the plane wave ($m=0$) and of the $\sech$ potential ($m=1$).
We show that, if $A\in\N$, then the spectrum consists of $\R$ plus $2A$ Schwarz symmetric segments (bands) on $\i\R$.
This characterization of the spectrum is obtained by relating the periodic and antiperiodic eigenvalue problems for the Dirac operator 
to corresponding eigenvalue problems for tridiagonal operators acting on Fourier coefficients in a weighted Hilbert space,
and to appropriate connection problems for Heun's equation. 
Conversely, if $A\not\in\N$, then the spectrum of $L$ consists of infinitely many bands in $\C$.
When $A\in\N$, the corresponding potentials generate finite-genus solutions for all the positive and negative 
flows associated with the focusing nonlinear Schr\"odinger hierarchy,
including the modified Korteweg-deVries equation and the sine-Gordon equation.
\end{abstract}
\vspace{5mm}
\tableofcontents

\bigskip
\section{Introduction and main results}
\label{s:introduction}

\subsection{Background}
\label{s:background}
In this work we study a non-self-adjoint Dirac operator with a Jacobi elliptic potential,
namely,
\be
\label{e:Diraceigenvalueproblem}
L\phi = z\phi\,, \quad z\in\C\,,
\ee
where $\phi(x;z)=(\phi_1,\phi_2)^\T$, 
the superscript ``$\T$'' denoting matrix transpose,
\be
L := \i\sigma_{3}(\partial_{x} - Q(x))\,, 
\qquad Q(x) = \begin{pmatrix} 0 & q(x) \\ -\overline{q(x)} & 0 \end{pmatrix}\,, \quad x\in\R\,, 
\label{e:Diracoperator}
\ee
the potential $Q(x)$ is $l$-periodic,
$\sigma_3:=\diag(1,-1)$ (cf. Appendix~\ref{a:notation})
and overline denotes the complex conjugate.
In particular, let
\be
q(x;A,m) = A\dn(x;m)\,,
\label{e:ellipticpotential}
\ee
where $\dn(x;m)$ is one of the three basic Jacobi elliptic functions (cf.~\cite{Gradshteyn,NIST}),
and $m\in(0,1)$ is the elliptic parameter.
Finally, $A$ is an arbitrary constant, which one can take to be real and positive without loss of generality.
(It is easy to see that $\arg A\ne0$ leaves the spectrum invariant.)
We will do so throughout this work.
Recall that $\dn(x;m)$ has minimal period $l = 2K$ along the real $x$-axis, 
where $K:=K(m)$ is the complete elliptic integral of the first kind~\cite{Gradshteyn,NIST}.
Also recall that
$\dn(x;0) \equiv 1$ and 
$\dn(x;1) \equiv \sech x$.
Both of the limiting cases $m=0$ and $m=1$ are exactly solvable 
(i.e., the spectrum is known in closed form),
and therefore provide convenient ``bookends'' for the results of this work.

There are several factors that motivate the present study.
A first one is that Dirac operators arise naturally in quantum field theory \cite{IZ1980,Weinberg1995},
and therefore the identification of exactly solvable potentials is relevant in that context.
A second one is the obvious similarity between the study of~\eqref{e:Diraceigenvalueproblem} and that of eigenvalue problems for the time-independent Schr\"odinger equation, namely
\vspace*{-1ex}
\be
(-\Delta + V(x))\,\phi = \lambda\,\phi\,,
\label{e:schrodinger_ndim}
\ee
where $\Delta$ denotes the $n$-dimensional Laplacian operator and $\phi:\Real^n\to\Complex$,
which has been an integral component of mathematical physics since its first appearance in the 1920's
(e.g., see \cite{GDGPS2006,Messiah,ReedSimon}),
and which received renewed interest in the late 1960's and 1970's
(e.g., see \cite{AS1981,DeiftTrubowitz,Krichever1975,NMPZ,Trubowitz})
thanks to the connection with infinite-dimensional integrable systems.
Namely, the fact that the one-dimensional time-independent Schr\"odinger equation 
[i.e., \eqref{e:schrodinger_ndim} with $n=1$] 
is the first half of the Lax pair for the \textit{Korteweg-deVries} (KdV) equation~\cite{GGKM,Lax}.
As a result,
the study of direct and inverse spectral problems for the Schr\"odinger operator played a key role in the development of the so-called \textit{inverse scattering transform} (IST) to solve the initial value problem for the KdV equation \cite{GGKM,Lax}.
The direct and inverse scattering theory was later made more rigorous, 
and generalizations of the theory were also studied
\cite{BealsCoifman,BealsDeiftTomei,date,DeiftTrubowitz,dubrovin,ItsMatveev1974,
ItsMatveev1975,lax2,trubowitz,mclaughlinnabelek,novikov,trogdon}.
In particular,
the so-called finite-gap (or finite-band) solution became a primary object of study. 

Similar problems have been considered for~\eqref{e:Diraceigenvalueproblem},  
since it comprises the first half of the Lax pair associated to the
\textit{nonlinear Schr\"odinger} (NLS) equation, namely,
the partial differential equation (PDE)
\begin{equation}
\label{e:nls}
\i q_t + q_{xx} + 2s |q|^2 q = 0\,.
\end{equation}
Here $q:\Real\times\Real\to\Complex$, subscripts $x$ and $t$ denote partial differentiation and,
as usual, the sign $s = \pm1$ denotes the focusing and defocusing cases, respectively.
Similarly to the KdV equation, the NLS equation is an infinite-dimensional Hamiltonian system.
Also, similarly to the KdV equation, the NLS equation is a ubiquitous physical model.
In particular,
\eqref{e:nls} is a universal model describing the slow modulations of a weakly monochromatic dispersive wave envelope,
and therefore appears in many physical contexts,
such as deep water waves, nonlinear optics, plasmas, ferromagnetics and Bose-Einstein condensates 
(e.g., see~\cite{AS1981,Agrawal,PitaevskiiStringari}).
Therefore, the study of the NLS equation is of both theoretical and applicative interest.

In 1972 \cite{ZS1972},
Zakharov and Shabat showed 
that \eqref{e:nls} is the compatibility condition of the
matrix Lax pair
\vspace*{-1ex}
\bse
\label{e:LP}
\begin{gather}
\phi_x = (-\i z\,\sigma_3 + Q(x,t))\,\phi\,,
\label{e:ZS}
\\
\phi_t = (-2\i z^2\sigma_3 + H(x,t;z))\,\phi\,,
\end{gather}
\ese
with $\sigma_3$ as above, and
\begin{gather}
Q(x,t) =  \begin{pmatrix} 0 & q(x,t) \\  -s \={q(x,t)} & 0 \end{pmatrix}\,,\qquad
H(x,t;z) = 2zQ - \i\sigma_3(Q^2 - Q_x)\,.
\end{gather}
Following \cite{ZS1972}, \eqref{e:ZS} [i.e., the first half of the Lax pair]
is referred to as the \textit{Zakharov-Shabat} (ZS) scattering problem.
It is easy to see that \eqref{e:ZS}[with $s=1$]
is equivalent to \eqref{e:Diraceigenvalueproblem}.
Thus, the solution to~\eqref{e:nls} [with $s=1$] comprises the scattering potential $q$ in~\eqref{e:Diraceigenvalueproblem}.
Moreover, one can also show that time evolution of $q$ according to the focusing NLS equation~\eqref{e:nls} [with $s=1$]
amounts to an isospectral deformation of the potential for the Dirac operator~\eqref{e:Diracoperator}.

Scattering theory for the Zakharov-Shabat system have been studied extensively over the years.
In \cite{ZS1972} the IST for~\eqref{e:nls} in the focusing case with localized data,
i.e., with $q(x,t=0)\in L^1(\Real)$, was formulated.
Corresponding results for the defocusing case with constant boundary conditions (BCs), i.e., $|q(x,t)|\to q_o\ne0$ as $x\to\pm\infty$,
were obtained in \cite{ZS1973}.
The theory was then revisited and elucidated in \cite{AS1981,FaddeevTakhtajan,NMPZ}.
When $q\in L^1(\Real)$, 
the isospectral data is composed of two pieces:
an absolutely continuous spectrum,
and a set of discrete eigenvalues.
When $q$ is periodic, however, the isospectral data is purely absolutely continuous and has a band and gap structure.

Of particular interest is the effort to find classes of potentials for which the scattering problem can be solved exactly.
Satsuma and Yajima \cite{SatsumaYajima} considered the case of $q(x) = A\,\sech x$,
with $A$ an arbitrary positive constant,
and obtained a complete representation of eigenfunctions and scattering data.
Their work was  later generalized by Tovbis and Venakides \cite{TovbisVenakides}
to potentials of the type $q(x) = A\,\sech x\,\e^{-\i a\,\log(\cosh x)}$, with $A$ as above and $a$ an arbitrary real constant.
These results were then used in~\cite{KMM2003,TVZ2004} to study the behavior of solutions of the focusing NLS equation
in the semiclassical limit.
More recently, Trillo et alii \cite{FratalocchiTrillo2008} obtained similar results 
for potentials of the type $q(x) = A\,\tanh x$ 
in the defocusing case.
In all of these cases, 
the ZS scattering problem is reduced to connection problems for the hypergeometric equation.
Finally, Klaus and Shaw \cite{klausshaw1,klausshaw2}
identified classes of ``single-lobe'' potentials for which the point spectrum is purely imaginary.

The above-mentioned works considered potentials that are either localized
or tend to constant boundary conditions as $|x|\to\infty$.
Spectral problems for the Schr\"odinger operator with a periodic potential similar to the one considered here 
are also a classical subject, 
and their study goes back to Lam\'e \cite{Lame},
and Ince \cite{Ince,Ince2,Ince3}
where the spectrum for a two-parameter family of potentials was studied, and necessary and sufficient conditions
in order for such potentials to give rise to a spectrum with a finite number of gaps were derived (see also \cite{Arscott,eastham2,Eastham,GW1996,MW1966}).
More recently, these results were generalized in~\cite{TreibichVerdier} and~\cite{Takemura}, and in seminal work
a characterization of all elliptic algebro-geometric solutions of the KdV and AKNS hierarchies was given by Gesztesy and Weikard in \cite{GW1996,GW1998,GW1998a}.  

Finite-band potentials for the focusing and defocusing ZS scattering problems have also been studied 
\cite{BBEIM1994,GW1998,ItsKotlyarov,Krichever1976,Smirnov}.
In particular, the special case of genus-one potentials was explicitly considered in \cite{CarterSegur,Kamchatnov},
and the stability of those solutions was recently studied in \cite{DeconinckSegal}.
On the other hand, 
the identification of exactly solvable cases for periodic potentials is generally challenging,
and few families of finite-band potentials for~\eqref{e:Diracoperator} have been studied in detail (see \cite{GW1998,GW1998a}).

Here we present an explicit, two-parameter family of finite-band potentials of the focusing ZS system and we characterize the resulting spectrum.
We also show that~\eqref{e:Diraceigenvalueproblem} with potential \eqref{e:ellipticpotential} can be reduced to certain connection problems for Heun's equation.
Unlike the case of the hypergeometric equation,
the connection problem for Heun's equation has not been solved in general \cite{Ronveaux}.
Still, 
special cases can be solved exactly. For example, for
certain classes of periodic potentials it turns out that Hill's equation 
[i.e., \eqref{e:schrodinger_ndim} with $n=1$ and periodic potential] 
can be mapped to a Heun equation.
Classical works~\cite{Ince,Ince2,MW1966}
where the spectrum of Hill's equation for a multi-parameter family of potentials was studied,
resulted in the derivation of necessary and sufficient conditions for such potentials to give rise to a spectrum with a finite number of bands and gaps.
Importantly,
the absence of a gap in the spectrum of the Hill operator corresponds uniquely to the coexistence of solutions,
namely,
the existence of two linearly independent periodic,
or antiperiodic,
solutions to the given ordinary differential equation (ODE)~\cite{MW1966}. 
More recently, those results were strengthened in~\cite{volkmer,volkmer2} and \cite{veselov}.
The results of this work provide a direct analogue of all these results for the Dirac operator~\eqref{e:Diracoperator}
as well as for the Hill operator with PT-symmetric potential.

\subsection{Main results}

We first introduce some definitions in order to state the main results of this work.
Various notations and standard definitions are given in Appendix~\ref{a:notation}.
We view the matrix-valued differential expression $L$, introduced in~\eqref{e:Diracoperator}, as a 
densely and maximally defined closed linear operator 
acting on $L^2(\R;\C^2)$,
i.e.,
the space of square integrable complex-valued vector functions~\cite{GW1998}. 
\begin{definition}
(Lax spectrum)
The Lax spectrum of the Dirac operator $L$ in~\eqref{e:Diracoperator} is the set
\be
\Sigma(L) := \{z\in\Complex : L\phi=z\phi,\,\, 0\neq \| \phi \|_{\infty}<\infty \}\,.
\label{e:laxspec}
\ee
\end{definition}

It is well-known that if $q$ is $l$-periodic, then $\Sigma(L)$
is purely absolutely continuous and comprised of an at most countable collection of regular analytic arcs,
referred to as \textit{bands},
in the spectral plane~\cite{GW1998,rofebek}.
Further properties of the Lax spectrum are discussed in Section~\ref{s:prelim}.
If the potential $q$ is such that there are at most finitely many bands we say that 
$q$ is a \textit{finite-band potential}
(see Definition~\ref{d:bandsandgaps}).
The class of finite-band potentials plays a key role in the IST for the NLS equation on the torus~\cite{BBEIM1994,Gesztesy,MaAblowitz}.
In particular, it was shown in \cite{ItsKotlyarov} that the potential 
can be reconstructed from the knowledge of two key spectral data:
(i) the periodic and antiperiodic eigenvalues of $L$
(i.e., the set of values~$z$ associated with periodic or antiperiodic eigenfunctions, respectively), 
which correspond to endpoints of spectral bands, 
and 
(ii) the Dirichlet (or auxiliary) eigenvalues of $L$, 
defined as the set of zeros of the\ 1,2 entry of the monodromy matrix
(see Section~\ref{s:prelim} for precise definitions of all these quantities).
To specify the dependence of the spectrum associated with~\eqref{e:ellipticpotential} on the parameters $A,m$, we will also occasionally 
use the notation $\Sigma(L;A,m)$ for the Lax spectrum.

\begin{theorem}
\label{t:mainresult}
Consider \eqref{e:Diraceigenvalueproblem} with
$q\equiv A\dn(x;m)$,
$m\in(0,1)$, and $A>0$.
Then the potential $q$ is finite-band if and only if $A\in\Natural$.
Moreover, if $A\in\Natural$, then:
\be
\Sigma(L;A,m) \subset \R \cup (-\i A, \i A)\,,
\label{e:subsetspec}
\ee
and $q$ is a $2\!A$-band (i.e., a genus $2A-1$) potential of the Dirac operator \eqref{e:Diracoperator}. 
\end{theorem}

\noindent 
(Of course it is well known that $\Sigma(L)$ is Schwarz symmetric and $\R\subset\Sigma(L;A,m)$ \cite{MaAblowitz,mclaughlinoverman}.)
Theorem~\ref{t:mainresult} is a consequence of the following more detailed description of the spectrum:

\begin{theorem}
\label{t:Ainteger}
Assume the conditions of Theorem~\ref{t:mainresult}. If $A\in\Natural$ then: 
\begin{enumerate}
\advance\itemsep-4pt
\item 
For any $m\in(0,1)$, 
the non-real part of the Lax spectrum, $\Sigma(L;A,m)\setminus \R$, is a proper subset of~~$(-\rm{i}\textit{A},\rm{i}\textit{A})$.
(For $m=0$, the Lax spectrum is $\Sigma(L;A,0)=\R\cup [-\rm{i}\textit{A},\rm{i}\textit{A}]$.)
\item 
For any $m\in(0,1)$,
there are exactly $2A$ symmetric bands of $\Sigma(L;A,m)$ along $(-\rm{i}\textit{A},\rm{i}\textit{A})$, 
separated by $2A-1$ open gaps. 
The central gap (i.e., the gap surrounding the origin) contains an eigenvalue at $z=0$, which is  
periodic when $A$ is even and antiperiodic when $A$ is odd.
\item 
For any $m\in[0,1)$,
$\Real\subset\Sigma(L;A,m)$ contains infinitely many interlaced  
periodic and antiperiodic eigenvalues, symmetrically located with respect to $z=0$.
\item 
Each periodic/antiperiodic eigenvalue $z\in\Real$ has geometric multiplicity two 
and each periodic/\break antiperiodic eigenvalue $z\in (-\rm{i}\textit{A},\rm{i}\textit{A})\setminus\{0\}$ has geometric multiplicity one.
\item
Each periodic/antiperiodic eigenvalue $z\in\Real$ is simultaneously a Dirichlet eigenvalue. 
All these Dirichlet eigenvalues are immovable.
\item
Each of the open $2A-1$ gaps on  $(-\rm{i}\textit{A},\rm{i}\textit{A})$ contains exactly one movable Dirichlet eigenvalue.
Thus, all of the $2A-1$ movable Dirichlet eigenvalues of the finite-band solution with genus $2A-1$ are located
in the gaps of the interval $(-\rm{i}\textit{A},\rm{i}\textit{A})$.
\end{enumerate}
\end{theorem}

Recall that a movable Dirichlet eigenvalue is a Dirichlet eigenvalue whose location changes 
when changing the normalization of 
the monodromy matrix,
whereas the location of immovable Dirichlet eigenvalues is independent of the normalization of the monodromy matrix.
For an $N$-band potential 
there are a total of $N-1$ movable Dirichlet eigenvalues
(cf.\ Definition~\ref{d:Dirichlet} and~\cite{ForestLee,Gesztesy}).

\begin{theorem}
\label{t:Anoninteger}
Assume the conditions of Theorem \ref{t:mainresult}.
If $A\not\in\Natural$, then:
\begin{enumerate}
\advance\itemsep-4pt
\item
For any $m\in(0,1)$, each periodic or antiperiodic eigenvalue has geometric multiplicity one. 
\item 
There are no periodic or antiperiodic eigenvalues on $\R$.
\item 
There are infinitely many spines (spectral bands emanating transversally from the real axis) 
at the real critical points of the Floquet discriminant (i.e., the trace of the monodromy matrix). 
\end{enumerate}
\label{t:mainresult-det}
\end{theorem}

Time evolution according to the NLS equation is an isospectral deformation of a potential of~\eqref{e:Diracoperator}.
Thus, by Theorem~\ref{t:mainresult}, 
if $A\in\Natural$, the initial condition $q(x,0) = A\dn(x;m)$ 
generates a genus $2A-1$ solution of the focusing NLS equation;
conversely, if $A\not \in \Natural$, the corresponding solution is not finite-genus.

Preceeded by Preliminaries in Section \ref{s:prelim},
the proof of Theorems~\ref{t:mainresult} and \ref{t:Ainteger} involves several steps:
\vspace*{-1ex}
\begin{itemize}
\advance\itemsep-6pt
\item 
In Section~\ref{s:Hill} we map~\eqref{e:Diraceigenvalueproblem} into 
Hill's equation with a complex potential, 
and in Section~\ref{s:trig} we map Hill's equation into a second-order trigonometric ODE. 
\item
In Section~\ref{s:recurrence} we map the trigonometric ODE into a three-term recurrence relation for the Fourier coefficients.  
\item
In Section~\ref{s:Fourier_ascendingdescending}
we demonstrate that, when $A\in\Natural$, each periodic or antiperiodic eigenvalue of $L$ is associated to a corresponding ascending or descending semi-infinite Fourier series.
\item
In Section~\ref{s:heun} 
we map the trigonometric ODE into Heun's equation and relate the periodic and antiperiodic eigenvalue problems for~\eqref{e:Diraceigenvalueproblem} with potential~\eqref{e:ellipticpotential} to a connection problem for Heun's equation.
\item
Moreover, in Section~\ref{s:heun} we show that the periodic and antiperiodic eigenvalues of~\eqref{e:Diracoperator} with potential~\eqref{e:ellipticpotential} correspond to the eigenvalues of certain tridiagonal operators that encode the recurrence relations for the coefficients of the Frobenius series solution of Heun's equation at the origin and at infinity.
\item
In Section~\ref{s:realevalues1}
we establish that all eigenvalues of the above-mentioned tridiagonal operators are real.
\end{itemize}
The determination of the precise number of spectral bands for any $m\in(0,1)$ is proved in Section~\ref{s:genus}.
Finally,
Theorem~\ref{t:Anoninteger} is proved in Section~\ref{s:noninteger}.
Notation, standard definitions, several technical statements and additional results and observations are relegated to the appendices.

\section{Preliminaries}
\label{s:prelim}

We begin by briefly reviewing basic properties of the Lax spectrum.
Unless stated otherwise, 
all statements in this section hold for operators $L$ with arbitrary continuous $l$-periodic potentials. 

\subsection{Bloch-Floquet theory}
\label{s:floquet}

While it is natural to pose \eqref{e:Diraceigenvalueproblem} on the whole real $x$-axis, 
all of the requisite information for the spectral theory is contained in the period interval of the potential,
namely,
$I_{x_o}:=[x_o,x_o+l]$,
where $x=x_o$ is an arbitrary base point. 
Consider the \textit{Floquet boundary conditions} (BCs):
\be
\BC_\nu(L):=\{\phi~:~\phi(x_o+l;z) = \e^{\i\nu l}\phi(x_o;z)\,,\,\, \nu\in\R\}\,.
\label{e:floqBC}
\ee
\begin{definition}
\label{d:floquet}
(Floquet eigenvalues of the Dirac operator)\,
Let the operator $L:H^1(I_{x_o};\C^2) \to L^2(I_{x_o};\C^2)$
be defined by~\eqref{e:Diracoperator}.
Let~${\rm dom}(L) := \{\phi\in H^1(I_{x_o};\C^2) : \phi\in \BC_\nu(L)\}$.
The set of Floquet eigenvalues of $L$ is given by  
\be
\Sigma_\nu(L) := \{z\in\C: \exists\phi\not\equiv 0 \in {\rm dom}(L)\,\,\, {\rm s.t.}\,\,\, L\phi=z\phi\}\,.
\label{e:floqspec}
\ee
In particular,
$\nu=2n\pi/l$, $n\in\Z$, identifies periodic eigenfunctions,
while $\nu=(2n-1)\pi/l$, $n\in\Z$, identifies antiperiodic eigenfunctions. 
We will call the corresponding eigenvalues periodic and antiperiodic, respectively,
and we will denote the set of periodic and antiperiodic eigenvalues by $\Sigma_\pm(L)$, respectively.
\end{definition}
\noindent 
($H^1$ denotes the space of square-integrable functions with square-integrable first derivative.)
It is well-known that $\Sigma_\nu(L)$ is discrete and countably infinite \cite{eastham2,Mityagin}. 

Next we review the 
theory of linear homogeneous ODEs with periodic coefficients and important connections to the Lax spectrum.
We set the base point $x_o=0$ without loss of generality.
Recall, the \textit{Floquet solutions (or Floquet eigenfunctions)} of~\eqref{e:Diraceigenvalueproblem} are solutions such that 
\be
\phi(x+l;z) = \mu\,\phi(x;z)\,,
\label{e:normalsoln}
\ee 
where $\mu:=\mu(z)$ is the~Floquet multiplier. 
Then by Floquet's Theorem (see~\cite{eastham2,Floquet}) 
all bounded (in $x$) Floquet solutions of~\eqref{e:Diraceigenvalueproblem} have the form
$\phi(x;z) = \e^{\i\nu x}\psi(x;z)$,
where $\psi(x+l;z)=\psi(x;z)$ and $\nu:=\nu(z) \in \Real$.
Thus, 
a solution of~\eqref{e:Diraceigenvalueproblem} is bounded for all $x\in\R$ if and only if $|\mu|=1$, 
in which case one has the relation
\be
\mu = \e^{\i\nu l}\,,
\label{e:Floquetexponent}
\ee
with $\nu\in\Real$. 
The quantity $\i\nu$ is the Floquet exponent. 
(With a slight abuse of terminology, we will often simply refer to $\nu$ as the Floquet exponent for brevity.)
The Floquet multipliers are the eigenvalues of the \textit{monodromy matrix} $M:=M(z)$, 
defined by
$Y(x+l;z)=Y(x;z)M(z)$,
where $Y(x;z)$ is any fundamental matrix solution of~\eqref{e:Diraceigenvalueproblem}.
It is well-known that the monodromy matrix is entire as a function of $z$~\cite{MaAblowitz,mclaughlinoverman}.
Note that $\det M(z)\equiv 1$ $\forall z\in\Complex$ by
Abel's formula,
since~\eqref{e:Diracoperator} is traceless. 
Thus, the eigenvalues of $M$ are given by the roots of the quadratic equation
$\mu^2 - 2\D\,\mu + 1 = 0$,
where $\D:=\D(z)$ is the \textit{Floquet discriminant},
i.e.,
\be
\D(z) = \half\tr M(z)\,.
\label{e:DeltaM}
\ee
Further,
$\mu_{\pm} = \D \pm \sqrt{\D^2 - 1}$.
Thus~\eqref{e:Diraceigenvalueproblem} admits bounded solutions if and only if
$-1\le \D \le 1$.
\begin{remark}
\label{r:asymp}
For $q\in C(\Real)$ one has
$\D(z) = \cos(zl) + o(1)$ as $z\to\infty$ along the real $z$-axis (see~\cite{MaAblowitz,mclaughlinoverman}).
\end{remark}

The above considerations yield
an equivalent representation of the Lax spectrum (see \cite{eastham2,Eastham,rofebek}):

\begin{theorem}
The Lax spectrum $\Sigma(L)$ is given by 
\be
\Sigma(L)=\{z\in\C : \D(z)\in [-1,1]\}\,.
\label{e:laxspec2}
\ee
Additionally, 
for any
fixed $\nu\in\R$ the Floquet eigenvalues are given by
\be
\Sigma_\nu(L)=\{z\in\Complex : \D(z)=\cos(\nu l)\}\,.
\label{e:bloch}
\ee
\unskip
For each $\nu\in\Real$ the set $\Sigma_\nu(L)$ is discrete and the only accumulation point occurs at infinity.
Moreover,
\be
\Sigma(L) = \bigcup_{\nu \in [0,2\pi/l)} \Sigma_\nu(L)\,.
\label{e:floquet}
\ee
\label{t:rofebeketov}
\end{theorem}

\begin{remark}
By \eqref{e:normalsoln}, \eqref{e:Floquetexponent} and \eqref{e:bloch},
the values $z\in\C$ for which 
$\D(z)=\pm 1$ are the periodic and antiperiodic eigenvalues $z\in\Sigma_\pm(L)$ (see Definition~\ref{d:floquet}), respectively.
The periodic and antiperiodic eigenvalues correspond to band edges of the Lax spectrum.
Further,
$\Sigma_\nu(L) \cap \Sigma_{\nu'}(L) = \emptyset$ for all $\nu\ne\nu'\mod~2\pi/l$. 
\end{remark}

\subsection{General properties of the Lax spectrum}
\label{s:generalproperties}

Owing to~\eqref{e:laxspec2},
the Lax spectrum~\eqref{e:laxspec} is located along the zero level curves 
of $\Im\D(z)$, i.e.,
$\Gamma := \{z\in\C : \Im\D(z) = 0\}$.
Moreover, 
$\Gamma$ is the union of an at most countable set of regular analytic curves $\Gamma_n$
\cite{GW1998}, 
each starting from infinity and ending at infinity: 
\vspace*{-1ex}
\be
\Gamma = \cup_{n\in\Natural}\Gamma_n\,.
\ee
%

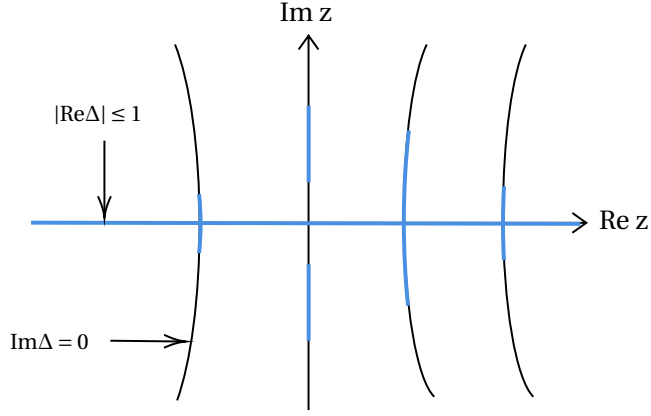
\begin{figure}[t!]
\begin{center}
\tikzset{every picture/.style={line width=0.75pt}}        
\begin{tikzpicture}[x=0.75pt,y=0.75pt,yscale=-1,xscale=1]
\draw  (201,140.1) -- (480,140.1)(339.94,45.6) -- (339.94,234.6) (473,135.1) -- (480,140.1) -- (473,145.1) (334.94,52.6) -- (339.94,45.6) -- (344.94,52.6)  ;
\draw    (100,52.6) ;
\draw  [draw opacity=0] (453.35,227.87) .. controls (444.6,219.49) and (438,181.59) .. (438,136.1) .. controls (438,97.48) and (442.76,64.32) .. (449.55,50.06) -- (457,136.1) -- cycle ; \draw   (453.35,227.87) .. controls (444.6,219.49) and (438,181.59) .. (438,136.1) .. controls (438,97.48) and (442.76,64.32) .. (449.55,50.06) ;  
\draw  [draw opacity=0] (403.35,227.87) .. controls (394.6,219.49) and (388,181.59) .. (388,136.1) .. controls (388,97.48) and (392.76,64.32) .. (399.55,50.06) -- (407,136.1) -- cycle ; \draw   (403.35,227.87) .. controls (394.6,219.49) and (388,181.59) .. (388,136.1) .. controls (388,97.48) and (392.76,64.32) .. (399.55,50.06) ;  
\draw  [draw opacity=0] (273.66,229.71) .. controls (280.5,210.65) and (285,178.16) .. (285,141.3) .. controls (285,102.53) and (280.02,68.61) .. (272.59,50.04) -- (259.5,141.3) -- cycle ; \draw   (273.66,229.71) .. controls (280.5,210.65) and (285,178.16) .. (285,141.3) .. controls (285,102.53) and (280.02,68.61) .. (272.59,50.04) ;  
\draw [color={rgb, 255:red, 74; green, 144; blue, 226 }  ,draw opacity=1 ][line width=1.5]    (200,140) -- (477,140.6) ;
\draw [color={rgb, 255:red, 74; green, 144; blue, 226 }  ,draw opacity=1 ][line width=1.5]    (340,81) -- (340,119.6) ;
\draw [color={rgb, 255:red, 74; green, 144; blue, 226 }  ,draw opacity=1 ][line width=1.5]    (340,161) -- (340,199.6) ;
\draw  [draw opacity=0][line width=1.5]  (438.5,159.01) .. controls (438.17,153.04) and (438,146.89) .. (438,140.6) .. controls (438,134.14) and (438.18,127.81) .. (438.53,121.68) -- (468,140.6) -- cycle ; \draw  [color={rgb, 255:red, 74; green, 144; blue, 226 }  ,draw opacity=1 ][line width=1.5]  (438.5,159.01) .. controls (438.17,153.04) and (438,146.89) .. (438,140.6) .. controls (438,134.14) and (438.18,127.81) .. (438.53,121.68) ;  
\draw  [draw opacity=0][line width=1.5]  (390.14,181.92) .. controls (388.76,168.86) and (388,154.58) .. (388,139.6) .. controls (388,123.16) and (388.92,107.56) .. (390.56,93.49) -- (417,139.6) -- cycle ; \draw  [color={rgb, 255:red, 74; green, 144; blue, 226 }  ,draw opacity=1 ][line width=1.5]  (390.14,181.92) .. controls (388.76,168.86) and (388,154.58) .. (388,139.6) .. controls (388,123.16) and (388.92,107.56) .. (390.56,93.49) ;  
\draw  [draw opacity=0][line width=1.5]  (284.73,155.4) .. controls (285.23,150.69) and (285.5,145.73) .. (285.5,140.6) .. controls (285.5,135.39) and (285.22,130.34) .. (284.71,125.56) -- (264.5,140.6) -- cycle ; \draw  [color={rgb, 255:red, 74; green, 144; blue, 226 }  ,draw opacity=1 ][line width=1.5]  (284.73,155.4) .. controls (285.23,150.69) and (285.5,145.73) .. (285.5,140.6) .. controls (285.5,135.39) and (285.22,130.34) .. (284.71,125.56) ;  
\draw    (237,98.6) -- (237,135.6) ;
\draw [shift={(237,137.6)}, rotate = 270] [color={rgb, 255:red, 0; green, 0; blue, 0 }  ][line width=0.75]    (10.93,-3.29) .. controls (6.95,-1.4) and (3.31,-0.3) .. (0,0) .. controls (3.31,0.3) and (6.95,1.4) .. (10.93,3.29)   ;
\draw    (240,199.6) -- (276,199.98) ;
\draw [shift={(278,200)}, rotate = 180.6] [color={rgb, 255:red, 0; green, 0; blue, 0 }  ][line width=0.75]    (10.93,-3.29) .. controls (6.95,-1.4) and (3.31,-0.3) .. (0,0) .. controls (3.31,0.3) and (6.95,1.4) .. (10.93,3.29)   ;

\draw (485,132) node [anchor=north west][inner sep=0.75pt]   [align=left] {Re z};
\draw (324,27) node [anchor=north west][inner sep=0.75pt]   [align=left] {Im z};
\draw (210,79.4) node [anchor=north west][inner sep=0.75pt]  [font=\footnotesize]  {$|\text{Re} \Delta |\leq 1$};
\draw (188,194.4) node [anchor=north west][inner sep=0.75pt]  [font=\footnotesize]  {$\text{Im} \Delta =0$};
\end{tikzpicture}
\end{center}
\caption{Schematic diagram of the Lax spectrum for a generic potential.}
\label{f:f1}
\end{figure}

\noindent (The precise details of the map $n\mapsto\Gamma_n$ are not important for the present purposes.) 
Different curves $\Gamma_i \ne \Gamma_j$ (and therefore different spectral bands) can 
intersect at saddle points of $\D(z)$.
However, two distinct $\Gamma_n$ can intersect at most once,
as a result of the fact that 
each $\G_n$ is a level curve of the harmonic function $\Im\D(z)$.
Thus the Lax spectrum $\Sigma(L)$ cannot contain any closed curves in the finite $z$-plane.

\begin{definition}
\label{d:bandsandgaps} 
(Spectral band)
A spectral band is a maximally connected regular analytic arc along~$\Gamma_n$ where
$\D(z)\in[-1,1]$ holds.
Each finite portion of ~$\Gamma_n$ where $|\Re \D(z)| > 1$, delimited by a band on either side, is called a spectral gap. 
\end{definition}
\begin{lemma}
(\cite{BOT2020})
\label{l:connected}
The real $z$-axis is the only band extending to infinity;
$\Sigma(L)$ contains no closed curves in the finite $z$-plane;
and the resolvent set $\varrho:=\C\setminus\Sigma(L)$ is comprised of two connected components.
\end{lemma}
With the above definition, 
the Lax spectrum can be decomposed into bands and gaps along each $\Gamma_n$ as in a self-adjoint problem,
with the crucial difference that here
bands and gaps are not restricted to lie along the real $z$-axis
(as they would be in a self-adjoint problem),
but lie instead along arcs of~$\Gamma_n$.
Figure~\ref{f:f1} provides a schematic illustration of these concepts.

We call a spectral band intersecting the real or imaginary $z$-axis transversally a \textit{spine} \cite{mclaughlinoverman}.
Generically, 
the Lax spectrum of the operator~\eqref{e:Diracoperator} includes infinitely many spines emanating from (infinitely many) critical points that extend to ($\pm$) infinity along the real $z$-axis \cite{Mityagin,mclaughlinoverman},
in which case we call $q$ an \textit{infinite-band potential}.
Otherwise,
we call $q$ a \textit{finite-band potential}.
Specifically,
if there are $N$ bands (not including the real $z$-axis) we say that $q$ is an $N$-\textit{band potential}.
The corresponding solutions of the focusing NLS equation are described in terms of Riemann $\Theta$-functions determined by hyperelliptic Riemann surfaces of genus $G=N-1$ (see \cite{BBEIM1994,ForestLee,Gesztesy,ItsKotlyarov,TracyChen}).
For example, $q\equiv A$ is a genus-0 (i.e., a 1-band) potential of the Dirac operator~\eqref{e:Diracoperator},
and $q\equiv \dn(x;m)$ is a genus-1 (i.e., a 2-band) potential.

\begin{remark}
The following sets play a key role in the analysis: 
\vspace*{-1ex}
\begin{itemize}
\advance\itemsep-4pt
\item 
Periodic/antiperiodic points: points $z_{\pm}\in\C$ such that
$\D(z_{\pm})= \pm1$ (note $z_{\pm}\in\Sigma_\pm(L)$);
\item 
Critical points: 
points $z_c\in\C$ such that $\partial_{z}\D(z_c)=0$.
\end{itemize}
\end{remark}

We denote by $\Phi(x;z)$
the fundamental matrix solution of~\eqref{e:Diraceigenvalueproblem} normalized so that $\Phi(0;z)\equiv\1$,
where $\1$ is the $2\times 2$ identity matrix.
The trace and the eigenvalues of the monodromy matrix $M(z)$ are independent of the particular fundamental matrix solution
chosen, and therefore so is the Floquet discriminant $\D(z)$ and the Lax spectrum $\Sigma(L)$.
Nonetheless,
it will be convenient to use $\Phi(x;z)$, so that $M(z)$ is simply given by
\be
M(z)=\Phi(l;z)\,.
\label{e:Msymnormalizeddef}
\ee

\begin{remark}
\label{r:Msymmetries}
It is straightforward to see that, for all $z\in\Complex$, the monodromy matrix satisfies the same symmetries as the 
scattering matrix for the IST on the line (e.g., see~\cite{AS1981,APT2004,MaAblowitz})
\bse
\begin{gather}
M^{-1}(z) = \sigma_2 M^\T(z) \,\sigma_2\,,
\\
\overline{M(\overline{z})} = \sigma_2 M(z) \,\sigma_2\,.
\label{e:ZSschwarzsymmetry}
\end{gather}
\ese
Moreover, it is also straightforward to verify the following additional symmetries 
(e.g., see~\cite{BOT2020}).
If $q$ is real, then
\bse
\label{e:Msymsymmetries}
\be
M(-\overline{z}) = \overline{M(z)}\,,\qquad z\in\Complex\,.
\label{e:Msymreal}
\ee
Moreover, if $q$ is even, then
\be
M(-\overline{z}) = \sigma_1\overline{M^{-1}(z)}\sigma_1\,,\qquad z\in\Complex\,,
\label{e:Msymeven}
\ee
while if $q$ is odd, then
\be
M(-\overline{z}) = \sigma_2\overline{M^{-1}(z)}\sigma_2\,,\qquad z\in\Complex\,,
\label{e:Msymodd}
\ee
\ese
where $\sigma_1$ and $\sigma_2$ are the first and second Pauli spin matrices, respectively (see Appendix~\ref{a:notation}).
\end{remark}

The symmetry~\eqref{e:ZSschwarzsymmetry} for the monodromy matrix
implies that the Floquet discriminant satisfies the Schwarz symmetry
\be
\overline{\D(\overline{z})} = \D(z)\,,\qquad z\in\Complex\,.
\ee
Moreover, if $q$ is real or even or odd, \eqref{e:Msymsymmetries} implies additionally that $\D(z)$ is an even function:
\be
\D(-z) = \D(z)\,,\qquad z\in\Complex\,.
\ee
As a result, one has:
\begin{lemma}
\label{l:quartets}
If $q$ is real or even or odd, 
$\Sigma(L)$ is symmetric about the real and imaginary $z$-axes.
Thus, the Floquet eigenvalues come in symmetric quartets $\{z,\overline{z},-z,-\overline{z}\}$.
\end{lemma}

For $q$ real and even, it follows from~\eqref{e:Msymreal} and~\eqref{e:Msymeven} that
\begin{equation}
M(z) = \D(z)\,\1 + c(z)\,\s_3 - \i s(z)\,\s_2\,, \qquad z\in\Complex\,.
\label{e:Msym}
\end{equation}
Obviously,~\eqref{e:Msym} together with the fact that $\det M(z)\equiv1$, imply the relation 
\be
\label{det-eq}
\D^2(z) = 1 + c^2(z) - s^2(z)\,,\qquad z\in\Complex\,.
\ee
Equation~\eqref{e:Msymreal} also implies that, if $q$ is real, $M(z)$ is real when $z\in \i\Real$. 
Moreover, for $q$ real and even, one has:
\begin{equation}
\label{Dsc}
\D(z)=\D(\overline{z})=\overline{\D(z)}, \quad  
s(z)=s(\overline{z})=\overline{s(z)},  \quad  
c(z)=-c(\overline{z})=\overline{c(z)}\,,\qquad
z\in \i\R\,.
\end{equation}
That is, $\D(z)$, $s(z)$ and $c(z)$ are all real for $z\in \i\R$.
For $z\in\R$, $\D(z)$ and $s(z)$ are real, whereas $c(z)$ is purely imaginary. 
Finally, since $M(z)$ is entire, \eqref{Dsc} also implies
\be\label{even-sc}
s(-z) = s(z)\,,\qquad c(-z) = -c(z)\,,\qquad
z\in\Complex\,.
\ee

Next we show that the Lax spectrum of~\eqref{e:Diracoperator} with a non-constant potential is confined to an \textit{open} strip in the spectral plane.
The following Lemma is proved in Appendix~\ref{a:newboundproof},
and is instrumental for this work:
\begin{lemma}
\label{l:combinedlemma}
Suppose $q\in C(\R)$ is $l$-periodic.
(i) If $q$ is not constant and 
$z\in\Sigma(L)$, then $|\Im z| < \|q\|_{\infty}$.
(ii) If $q$ is real or even or odd, and $\Sigma_{\pm}(L)\subset\R\cup\i\R$,
then $\Sigma(L)\subset\R\cup\i\R$ and $q$ is finite-band.
\end{lemma}

To solve the inverse problem in the IST (namely, reconstructing the potential from the scattering data),
an auxiliary set of spectral data is also needed--the Dirichlet eigenvalues \cite{ForestLee,ItsKotlyarov}:
\begin{definition}
(Dirichlet eigenvalues)
Let $M(z)$ be defined by \eqref{e:Msymnormalizeddef}. 
The set of Dirichlet eigenvalues (see~\cite{ForestLee}) with base point $x_0=0$ is defined as
\be
\Sigma_{\rm Dir}(L;x_o=0) := \{z\in\C : s(z)=0\}\,.
\label{e:dirspec}
\ee
\end{definition}
In contrast to the Lax spectrum,
the Dirichlet eigenvalues are not invariant with respect to changes in the base point $x=x_o$,
or to time evolution of $q$ according to the focusing NLS equation.
Indeed, in the context of the integrability of NLS on the torus, 
the Dirichlet eigenvalues correspond to angle variables and are used to coordinatize the isospectral level sets.
As we discuss next,
the set of Dirichlet eigenvalues is discrete, 
consists of movable and immovable points, 
and the number of movable Dirichlet eigenvalues is tied to the genus of the corresponding Riemann surface (see~\cite{ForestLee,GW1998a}).

The monodromy matrix $M(z)$ in~\eqref{e:Msymnormalizeddef} was defined in terms of 
the fundamental matrix solution $\Phi(x;z)$ normalized as $\Phi(0;z)\equiv\1$.  
The monodromy matrix $M(z;x_o)$ associated with a ``shifted" solution $\~\Phi(x;x_o,z)$ normalized 
as $\~\Phi(x_o;x_o,z)\equiv\1$, with $x_0\in\R$, is given by
\be
\label{e:M(x_o)}
M(z;x_o)=\Phi(x_o;z)M(z)\Phi^{-1}(x_o;z).
\ee
Let $\Sigma_{\rm Dir}(L;x_o)$ be the corresponding set of Dirichlet eigenvalues. 
\begin{definition}
(Movable and immovable Dirichlet eigenvalues)
Let $z\in\Complex$ be a Dirichlet eigenvalue associated to the monodromy matrix $M(z;x_o)$ with a given base point $x=x_o$,
i.e., $z\in\Sigma_{\rm Dir}(L;x_o)$.
Following~\cite{ForestLee}, 
we say that $z$ is an \textit{immovable} Dirichlet eigenvalue if $z\in\Sigma_{\rm Dir}(L;x)$ for all $x\in\Real$.
Otherwise, we say $z\in\Complex$ is a \textit{movable} Dirichlet eigenvalue.
\label{d:Dirichlet}
\end{definition}
\begin{remark}
If $q$ is an $N$-band potential of the non-self-adjoint Dirac operator~\eqref{e:Diracoperator}, then the number of movable Dirichlet eigenvalues is $N-1$ (see~\cite{ForestLee,GW1998a}).
\end{remark}

An immediate consequence of~\eqref{det-eq} and the symmetries of $M(z)$, $\D(z)$, $c(z)$ and $s(z)$ is the following lemma, which will be useful later (see also \cite{ForestLee,mclaughlinoverman}):
\begin{lemma}
\label{lem-Dcs}
If $z\in\R$ and $|\D(z)|=1$, then $c(z)=s(z)=0$, so that $z$ is an  immovable Dirichlet eigenvalue.  
Conversely, if $s(z)=0$ with $z\in {\rm i}\R$, then $|\D(z)|\geq 1$.
\end{lemma}

\begin{lemma}
Let $z_\pm\in\Sigma_\pm(L)$.
If $\partial_{z}\D(z_\pm)\neq 0$, then the corresponding eigenspace has dimension one.
\label{l:opengap}
\end{lemma}
\textit{Proof.}
Suppose that there exists two linearly independent periodic (or antiperiodic) eigenfunctions.
Consider the normalized fundamental matrix solution $\Phi(x;z)$ of~\eqref{e:Diraceigenvalueproblem},
namely,
$ L\Phi(x;z) = z\Phi(x;z)$ with $\Phi(0;z)\equiv \1$.
Differentiating with respect to $z$ and using variation of parameters one gets
\vspace*{-1ex}
\be
\D_z = \half\tr\Big( -\i\Phi(l;z)\int_{0}^{l}\Phi^{-1}(x;z)\sigma_{3}\Phi(x;z)\,\d x\Big)\,.
\label{e:traceprime}
\ee
By Floquet's theorem $\Phi(l;z_{\pm}) = \pm\1$,
respectively.  
Then~\eqref{e:traceprime} yields $\partial_{z}\D(z_{\pm}) = 0$.
\qed

The following lemma is a direct consequence of Lemmas \ref{lem-Dcs}, and \ref{l:opengap}:
\begin{lemma}
If $z_\pm\in\Sigma_\pm(L)\cap\R$, then the 
geometric multiplicity is two and $\partial_{z}\D(z^{\pm})=0$, respectively.
\label{l:dim2}
\end{lemma}

\subsection{Limits $m\to0$ and $m\to1$; $z=0$}
\label{s:m=0&m=1}

The two distinguished limits $m\to0^+$ and $m\to1^-$  
of the two-parameter family of elliptic potentials~\eqref{e:ellipticpotential}
provide convenient limits of the results of this work. 
Interestingly, both of these limits yield 
exactly solvable models.
Here it will be convenient to keep track of the dependence on $m$ explicitly. 

Since $\dn(x;0) \equiv 1$, 
when $m=0$ the potential~\eqref{e:ellipticpotential} reduces to a constant background,
i.e.,
$q\equiv A$ with period $l=2K(0)=\pi$.
Thus, \eqref{e:Diraceigenvalueproblem} becomes a linear system of ODEs 
with constant coefficients, for which one easily obtains a fundamental matrix solution
\be
\Phi(x;z,m=0) = \e^{-\i(z\sigma_3-A\sigma_2)x}\,.
\label{e:constsoln}
\ee
Hence the monodromy matrix is 
\vspace*{-1ex}
\be
M(z,m=0) = \cos\big(\sqrt{z^2+A^2}\pi\big)\,\1 - \frac{\i\sin\big(\sqrt{z^2+A^2}\pi\big)}{\sqrt{z^2+A^2}}\big(z\sigma_3-A\sigma_2\big)\,,
\label{e:M_m=0}
\ee
implying $\Sigma(L;A,0) = \Real \cup [-\i A,\i A]$.
Further,
$z=\pm\i A$ are the only simple periodic eigenvalues;
all other periodic (resp. antiperiodic) eigenvalues are double points.
Hence, for any $A\ne0$,
$q\equiv A$ is a 1-band (i.e., genus-0) potential of~\eqref{e:Diraceigenvalueproblem}.
Moreover,
the associated solution of the focusing NLS equation 
[i.e., \eqref{e:nls} with $s=1$] is simply $q(x,t) = A\,\e^{2\i A^2t}$.

On the other hand, 
the limit $m\to1^-$ is singular, since $K(m)$, and therefore the period $l=2K(m)$ of the potential~\eqref{e:ellipticpotential}, 
diverges in this limit. 
Indeed, $\dn(x;1) \equiv \sech x$, so 
letting $m=1$ results in the eigenvalue problem~\eqref{e:Diraceigenvalueproblem} with potential 
$q\equiv A\sech x$. 
This case is also exactly solvable,
and was first studied by Satsuma and Yajima~\cite{SatsumaYajima}.
The point spectrum is comprised of a set of discrete eigenvalues located along the imaginary $z$-axis.
Moreover,
for $A\in\Natural$ the potential is reflectionless, and the point spectrum is given by 
$z_n = \i(n-1/2)$ for $n=1,\dots,A\,$.
That is, when $A\in\Natural$,
$q\equiv A\sech x$ corresponds to a pure bound-state $A$-soliton solution of the focusing NLS equation
\cite{SatsumaYajima}.
When $A=1$, the solution of the NLS equation~\eqref{e:nls} is simply $q(x,t) = \e^{\i t}\sech x$.
When $A>1$, the solutions are much more complicated \cite{BiondiniLiSchiebold,SatsumaYajima}.
Indeed, 
the potential $A\,\sech x$ was used to study the semiclassical limit of the focusing NLS equation 
in the pure soliton regime~\cite{KMM2003}.

Lastly, we discuss the origin $z=0$ of the spectral plane. 
When $z=0$, the ZS system~\eqref{e:ZS} admits closed-form solutions
(see Appendix~\ref{a:exactsolution}).
These solutions then allow one to obtain the following lemma, which is 
proved in Appendix~\ref{a:exactsolution}:
\begin{lemma}
\label{l:zequalzero}
Consider~\eqref{e:Diracoperator} with potential~\eqref{e:ellipticpotential} and $m\in [0,1)$.
If $A\in\N$ is even or odd, 
then $z=0$ is a periodic or antiperiodic eigenvalue, respectively, with geometric multiplicity two in each case.
\end{lemma}

\section{Transformation to Hill's equation}
\label{s:Hill}

In this section we introduce a transformation of \eqref{e:Diraceigenvalueproblem}
that will be instrumental in proving Theorem~\ref{t:mainresult}, and 
we consider the effect of this transformation on the Lax spectrum.


First we transform~\eqref{e:Diraceigenvalueproblem} to Hill's equation with a complex-valued potential
via the unitary linear transformation
\be
\phi \mapsto v = \Lambda \phi\,, \qquad \Lambda := \frac{1}{\sqrt{2}}\begin{pmatrix*}[r] 1 & \i \\ 1 & -\i \end{pmatrix*}\,,
\label{e:changevar_5}
\ee
where $v:=v(x;z^2)=(v^+, v^-)^\T$.
Differentiation of~\eqref{e:Diraceigenvalueproblem} and use of~\eqref{e:changevar_5}
show that, if $q$ in~\eqref{e:Diracoperator} is a real-valued differentiable potential, then~\eqref{e:changevar_5} maps~\eqref{e:Diraceigenvalueproblem} into the diagonal system
\vspace*{-1ex}
\be
Hv:= (-\partial_x^2 + Q^2 - {\rm i}Q_x\sigma_1)v = z^2v\,.
\label{e:uncoupled}
\ee
Or, in component form,
\be
\label{e:tindschrod5}
v_{xx}^{\pm} + (\pm{\rm i} q_{x} + z^2 + q^2)v^{\pm} = 0\,.
\ee
Equation~\eqref{e:tindschrod5} is Hill's equation with the complex (Riccati) potential
$V^{\pm}:=\mp{\rm i} q_x - q^2$.
Thus, \eqref{e:tindschrod5} amounts to the pair of eigenvalue problems
\vspace*{-1ex}
\be
H^{\pm}v^{\pm} = \lambda v^{\pm}\,, \quad \lambda := z^2\,,
\label{e:hillevalue}
\ee
where
\vspace*{-1ex}
\be
H^{\pm}:=-\partial_{x}^2 + V^{\pm}(x)\,.
\label{e:hilloperators}
\ee

\begin{remark}
If the potential $q$ in~\eqref{e:Diracoperator} is real and even, then 
$V^{\pm}(-x)=\overline{V^{\pm}(x)}$,
i.e.,
$V^\pm$ is PT-symmetric.
\end{remark}

Next, similarly to~\eqref{e:floqBC}, we introduce the corresponding~\textit{Floquet} BCs for $H^\pm$: 
\be 
\text{BC}_\nu(H^{\pm}) := \{v^{\pm} : 
v^{\pm}(l;\lambda)=\e^{\i\nu l}v^{\pm}(0;\lambda)\,,~~
v^{\pm}_x(l;\lambda)=\e^{\i\nu l}v^{\pm}_x(0;\lambda)\,,~~
\nu\in\R\}\,.
\ee
\begin{definition}
\label{d:floquethill}
(Floquet eigenvalues of Hill's operator)
Let the operators $H^{\pm}: H^2([0,l])\to L^2([0,l])$
be defined by~\eqref{e:hilloperators}.
Let ${\rm dom}(H^{\pm}):=\{v^\pm \in H^2([0,l]): v^\pm\in{\rm BC}_\nu(H^{\pm})\}$.
The set of Floquet eigenvalues of $H^\pm$ is given by
\be
\Sigma_\nu(H^{\pm}) := \{\lambda \in \C : \exists v^{\pm}\not\equiv 0 \in {\rm dom}(H^\pm) \,\,\, {\rm s.t.}\,\,\, H^{\pm}v^{\pm}=\lambda v^{\pm} \}\,.
\label{e:hillspec}
\ee
In particular,
$\nu=2n\pi/l$, $n\in\Z$, identifies periodic eigenfunctions,
while $\nu=(2n-1)\pi/l$, $n\in\Z$, identifies antiperiodic eigenfunctions. 
We will call the corresponding eigenvalues periodic and antiperiodic, respectively,
and we will denote the set of periodic and antiperiodic eigenvalues by $\Sigma_\pm(H^{\pm})$, respectively.
\end{definition}

\noindent ($H^2$ denotes the space of square-integrable functions with square-integrable first and second derivatives.) 
It is well-known that $\Sigma_{\nu}(H^{\pm})$ is discrete and countably infinite~\cite{eastham2,Mityagin,Eastham,MW1966}. 

\begin{lemma}
\label{l:equalspec}
If the potential $q$ in~\eqref{e:Diracoperator} is real and even, 
then $\Sigma_\nu(H^{+}) = \Sigma_{-\nu}(H^{-})$,
the dimension of the corresponding eigenspaces are equal,
and each of $\Sigma_{\nu}(H^{\pm})$ is symmetric about the real $\lambda$-axis.
\end{lemma}
\textit{Proof.}
Let $\lambda\in\Sigma_\nu(H^+)$ with eigenfunction $v^{+}(x;\lambda)$.
Since $q$ is even, it is easy to check $\~v:=v^{+}(-x;\lambda)$ satisfies $H^{-}\~v=\lambda \~v$.
Moreover,
\bse
\vspace*{-1ex}
\begin{gather}
\~v(l;\lambda) = v^{+}(-l;\lambda) = \e^{-\i\nu l}\~v(0;\lambda)\,,
\\
\~v_x(l;\lambda) = -v^{+}_x(-l;\lambda) = \e^{-\i\nu l}\~v_x(0;\lambda)\,.
\end{gather}
\ese
Hence,
$\lambda\in\Sigma_{-\nu}(H^-)$.
Conversely, if $\lambda\in\Sigma_{-\nu}(H^-)$ with eigenfunction $v^-(x;\lambda)$,
a completely symmetric argument shows that $\lambda\in\Sigma_{\nu}(H^+)$.
Finally,
since the map $v(x;\lambda) \mapsto v(-x;\lambda)$ is a (unitary) isomorphism, the dimension of the corresponding eigenspaces are the same.

Next, we prove the symmetry. 
Assume that $\lambda\in\Sigma_\nu(H^{\pm})$ with corresponding eigenfunction $v^{\pm}(x;\lambda)$, respectively.
Then it is easy to check that $\~v^{\pm}:=\overline{v(-x;\lambda)^{\pm}}$ satisfies $H^{\pm}\~v^{\pm}=\overline{\lambda}\~v^{\pm}$. Moreover,
\bse
\begin{gather}
\~v^{\pm}(l;\lambda) = \overline{v^{\pm}(-l;\lambda)} = \e^{\i\nu l}\~v^{\pm}(0;\lambda)\,,
\\
\~v^{\pm}_x(l;\lambda) = -\overline{v^{\pm}_x(-l;\lambda)} = \e^{\i\nu l}\~v^{\pm}_x(0;\lambda)\,.
\end{gather}
\ese
Thus,
$\overline{\lambda}\in \Sigma_\nu(H^{\pm})$ with eigenfunction $\~v^{\pm}(x;\lambda)$,
respectively.
\qed

\begin{remark}
It is easy to see that Lemma~\ref{l:equalspec} implies $\Sigma_\pm(H^+)=\Sigma_\pm(H^-)$, respectively.
\end{remark}

Next,
since the Lax spectrum $\Sigma(H^{\pm}) = \cup_{\nu\in[0,2\pi/l)}\Sigma_{\nu}(H^{\pm})$,
we have the following key equivalence:
\begin{lemma}
\label{l:ZSHill_spectrum_equal}
If the potential $q$ in~\eqref{e:Diracoperator} is real and even, then the unitary map~\eqref{e:changevar_5} implies:
\be
\Sigma(H^{+}) = \Sigma(H^{-}) = \{\lambda = z^2 : z\in\Sigma(L)\}\,.
\label{e:LequivHspectra}
\ee
That is,
the Lax spectrum of these three operators is related through the relation $\lambda = z^2$.
In particular, 
\be
z\in \Sigma_+(L) ~\Leftrightarrow~ \lambda = z^2 \in \Sigma_+ (H^\pm) \,,\quad
z\in \Sigma_-(L) ~\Leftrightarrow~ \lambda = z^2 \in \Sigma_- (H^\pm) \,.
\label{e:LequivHperiodicspectra}
\ee
Finally, for $z\ne 0$, the geometric multiplicity of an eigenvalue $z\in\Sigma_+(L)$ equals that of 
$\lambda=z^2\in\Sigma_+(H^\pm)$, and similarly for 
$z\in\Sigma_-(L)$ and $\lambda=z^2\in\Sigma_-(H^\pm)$.
\end{lemma}

\textit{Proof.}
If $z\in\Sigma(L)$,
the transformation~\eqref{e:changevar_5} implies that $v^{\pm}(x;\lambda)$ are both bounded solutions of Hill's ODE~\eqref{e:hillevalue},
respectively, implying $\lambda\in\Sigma(H^\pm)$.
Conversely,
if $v^{+}(x;\lambda)$ is a bounded solution of~\eqref{e:hillevalue} with the plus sign,
it follows that $\~v:=v^{+}(-x,\lambda)$ is a bounded solution of~\eqref{e:hillevalue} 
with the minus sign.
Further,
$\phi_1 = (v^{+}+v^{-})/\sqrt{2}$,
and $\phi_2 = \i(v^{-}-v^{+})/\sqrt{2}$ are both bounded, 
and the map~\eqref{e:changevar_5} then implies that
$\phi(x;z) = (\phi_1,\phi_2)^\T$ solves~\eqref{e:Diraceigenvalueproblem}, implying $z\in\Sigma(L)$. 
A similar argument follows if one starts with $v^-(x;\lambda)$ bounded.
Thus, \eqref{e:LequivHspectra} follows.
Equation \eqref{e:LequivHperiodicspectra} follows directly from Lemma~\ref{l:equalspec}.

It remains to show that, for $z\ne0$, the dimension of the corresponding eigenspaces are equal.
The argument follows~\cite{DjakovMityagin} where the self-adjoint case was studied.
To this end, let $E_{\pm}(L,z)$ denote the eigenspace associated with an eigenvalue $z\in\Sigma_{\pm}(L)$, 
and similarly for $L^2:=L\circ L$ and $H^\pm$.
First,
note that $\phi \mapsto \i\sigma_2\phi$ is a (unitary) isomorphism between the eigenspaces $E_{\pm}(L,z)$ and $E_{\pm}(L,-z)$.
Thus, applying the operator twice, for $z\ne0$ one easily gets
\be
\label{e:espace1}
{\rm dim}E_{\pm}(L^2,\lambda=z^2) = 2{\rm dim}E_{\pm}(L,z)\,.
\ee
Next,
note that $L^2$ is (unitary) equivalent to the diagonal system~\eqref{e:uncoupled},
i.e., 
$H = \half\Lambda L^2\Lambda^{-1}$ 
Moreover,
\be
\label{e:espace3}
H = \begin{pmatrix} H^{+} & 0 \\ 0 & 0 \end{pmatrix} + \begin{pmatrix} 0 & 0 \\ 0 & H^{-} \end{pmatrix}\,,
\ee
and so
\be
\label{e:espace4}
E_{\pm}(H,\lambda) = (E_{\pm}(H^+,\lambda)\oplus 0)\oplus (0\oplus E_{\pm}(H^-,\lambda))\,.
\ee
Hence,
by~\eqref{e:espace1}--\eqref{e:espace4} and Lemma~\ref{l:equalspec} it follows ${\rm dim}E_{\pm}(L,z) = {\rm dim}E_{\pm}(H^\pm,\lambda)$, respectively. 
\qed

\begin{remark}
By Lemma~\ref{l:ZSHill_spectrum_equal}, 
the spectrum of the Dirac operator $L$ in~\eqref{e:Diracoperator} with real and even potential is associated to that of the spectrum of the Hill operators $H^\pm$ in~\eqref{e:hilloperators}. 
Importantly,
note that the final statement of Lemma~\ref{l:ZSHill_spectrum_equal} does not hold at $z=0$;
that is,
the geometric multiplicity of the periodic (or antiperiodic) eigenvalue $z=0$ 
of the Dirac and Hill operators need not be equal (see Appendix~\ref{a:exactsolution}).
\end{remark}

All of the above results hold for generic real and even potentials.
Moving forward, we restrict our attention to the Jacobi elliptic potential~\eqref{e:ellipticpotential}.
By Lemma~\ref{l:ZSHill_spectrum_equal} we fix 
$y:=v^{-}(x;\lambda)$
without loss of generality (dependence on $A$ and $m$ is omitted for brevity).
Then Hill's equation $H^-v^-=\l v^-$ is given by 
\be
y_{xx} + (\i Am\sn(x;m)\cn(x;m) + \lambda + A^2\dn^2(x;m))y = 0\,.
\label{e:reduction1}
\ee
\begin{remark}
Since $\dn^2(x;m)\equiv 1-m\sn^2(x;m)$, \eqref{e:reduction1}
can be viewed as an imaginary deformation of the celebrated Lam\'e equation~\cite{Arscott,erdelyi,Ince2,MW1966}, $y_{xx} + (\lambda + V(x))\,y = 0$
up to a shift of the eigenvalue $\lambda$.
The Lam\'{e} equation has the remarkable property that solutions can coexist 
if and only if $A^2 = n(n+1)$ 
where $n$ is an integer~\cite{Arscott,erdelyi,Ince2,MW1966}.
Recall that solutions coexist if two linearly independent periodic (or respectively antiperiodic) solutions exist for a given $\lambda$.
In the case of Hill's equation with a real potential this amounts to a ``closed gap'' in the spectrum
(corresponding to finite gap potentials).
\end{remark}

\section{Transformation to a trigonometric ODE}
\label{s:trig}

In this section we introduce a second transformation of \eqref{e:Diraceigenvalueproblem}.
By part (ii) of Lemma~\ref{l:combinedlemma} moving forward we only need to consider the periodic and antiperiodic eigenfunctions. 

\subsection{Second-order ODE with trigonometric coefficients}
\label{s:trigcoef}
Consider the following change of independent variable:
\be
x \mapsto t:= 2\am(x;m)\,,
\label{e:amp}
\ee
where $\am(x;m)$ is the Jacobi amplitude~\cite{BirdFriedman,Gradshteyn}.
Equation~\eqref{e:amp} establishes a conformal map between the strip $|\Im x\,|<K(1-m)$ and the complex $t$-plane cut along the rays $(2j+1)\pi \pm 2\i\tau r$, 
$\tau\ge 1$,
$j\in\Z$,
where $r=\ln[(2-m)/m]/2$ \cite{Gradshteyn,NIST}.
We then arrive at our second reformulation of the Dirac eigenvalue problem:
\be
4(1-m\sin^2 \halft)y_{tt} - (m\sin t)y_{t} + 
(\lambda + A^2(1-m\sin^2 \halft) + \halfi Am\sin t)y = 0\,.
\label{e:trigonometricODE}
\ee
(The independent variable $t$ introduced above should not be confused with the time variable 
of the NLS equation~\eqref{e:nls}.)

\begin{remark}
Equation~\eqref{e:trigonometricODE} can be written as the eigenvalue problem
\vspace*{-0.6ex}
\be
By=\lambda y\,,
\label{e:trigvalue}
\ee
where the operator $B:H^2([0,2\pi])\to L^2([0,2\pi])$ is defined by
\vspace*{-0.6ex}
\be
B := -4(1-m\sin^2\halft)\partial_{t}^2 + (m\sin t)\partial_{t} - (A^2(1-m\sin^2\halft) + \halfi Am\sin t))\,.
\label{e:inceoperator}
\ee
The coefficients are now $2\pi$-periodic and as before $\Sigma_{\pm}(B)$ will denote the periodic and antiperiodic eigenvalues of the operator $B$, respectively (see Definition~\ref{d:floquethill}).
\end{remark} 
 
This leads to the following result which connects the periodic/antiperiodic eigenvalues of Hill's equation~\eqref{e:reduction1} to the periodic/antiperiodic eigenvalues of the trigonometric equation~\eqref{e:trigonometricODE}. 
\begin{lemma} 
Let $B$ be the trigonometric operator~\eqref{e:inceoperator}.
Then $\lambda\in \Sigma_{\pm}(B)$ 
if and only if 
$\lambda\in \Sigma_{\pm}(H^{-})$. 
\label{l:twopi}
\end{lemma}  
\textit{Proof.} 
By \ref{e:trigonometricODE} one gets $B y = \lambda y$ if and only if $H^{-}\~y = \lambda \~y$,
with $\~y(x;\lambda) = y(t;\lambda)$
and $t=2\am(x;m)$ as per \eqref{e:amp}.
Next,
note that $\am(x;m)$ is monotonic increasing for $x\in(0,2K)$,
$\am(x+2K;m)=\am(x;m)+\pi$, and $\am(0;m)=0$. 
Hence, the map \eqref{e:amp}
is a bijection between 
$x\in[0,2K]$ and $t\in[0,2\pi]$.
Moreover,
$\~y(0;\l) = \pm \~y(2K;\l)$ if and only if $y(0;\l) = \pm y(2\pi;\l)$.
Similarly, 
$\~y_{x}(0;\l) = \pm \~y_{x}(2K;\l)$ if and only if $y_{t}(0;\l) = \pm y_{t}(2\pi;\l)$.
Thus,
$2K$-periodic (resp. antiperiodic) solutions of~\eqref{e:reduction1} map to $2\pi$-periodic (resp. antiperiodic) solutions of~\eqref{e:trigonometricODE},
and vice versa.
\qed 

\begin{remark}
\label{r:inceconnection}
The trigonometric ODE~\eqref{e:trigonometricODE} can be viewed as a complex deformation of Ince's equation (see Chapter~7 of~\cite{MW1966} for more details).
Namely,
one can write~\eqref{e:trigonometricODE} as
\be
(1 + a\cos t)y_{tt} + (b\sin t)y_{t} + (h + d\cos t + {\rm i} e\sin t)y = 0\,,
\label{e:ince1}
\ee
where
$a = m/(2-m)$,
$b = -a/2$,
$h = \lambda/(4-2m) + A^2/4$,
$d = A^2a/4$,
$e = Aa/4$.
To the best of our knowledge this is the first example of a non-self-adjoint version of Ince's equation arising from applications.
\end{remark}

\subsection{Fourier series expansion and three-term recurrence relation}
\label{s:recurrence}
 
Recall that any Floquet solution $y(t;\l)$ of~\eqref{e:trigonometricODE} bounded for all $t\in\R$ has the form $y(t;\l) = \e^{\i\nu t}f(t;\l)$ where $f(t+2\pi;\l)=f(t;\l)$ and $\nu\in\R$ (cf. Section~\ref{s:floquet}).
Moreover, since $f(t;\l)$ is $2\pi$-periodic, 
we can express it in terms of a Fourier series on $L^2(\mathbb{S}^1)$,
where $\mathbb{S}^1:=\R/\Z$ is the unit circle. 
By direct calculation,
let $y(t;\l)$ be a Floquet solution of~\eqref{e:trigonometricODE} given by
\be
y(t;\l) = \e^{{\rm i}\nu t}\sum_{n\in \Z}c_{n}\e^{{\rm i} nt}\,.
\label{e:fourier}
\ee
Then the coefficients $\{c_n\}_{n\in\Integer}$ are given by the following three-term recurrence relation:
\be
\alpha_{n}c_{n-1} + (\beta_{n}-\lambda)c_{n} + \gamma_{n}c_{n+1} = 0\,,
\qquad
n\in\Z\,,
\label{e:Fourierrecurrence}
\ee
where
\vspace*{-1ex}
\bse
\label{e:alphabetagammadef}
\begin{align}
\alpha_{n} &= - \txtfrac 14 m\,[A-(2n+2\nu-2)][A+(2n+2\nu-1)]\,, \\
\beta_{n} &= (1-\half m)[(2n+2\nu)^2 - A^2]\,, \\
\gamma_{n} &= - \txtfrac 14 m\,[A-(2n+2\nu+2)][A+(2n+2\nu+1)]\,.
\end{align}
\ese

\begin{remark}
In turn, the recurrence relation~\eqref{e:Fourierrecurrence} can be written as the eigenvalue problem
\be
B_\nu c=\lambda c\,,
\label{e:matrixevalue}
\ee
where $c=\{c_n\}_{n\in\Z}$
\be
B_\nu := 
\begin{pmatrix*} 
\ddots & \ddots & \ddots & & \\ 
& \alpha_{n} & \beta_{n} & \gamma_{n} & \\ 
& & \ddots & \ddots & \ddots   
\end{pmatrix*}\,.
\label{e:trigoperatorBnu}
\ee
Note: $\nu\in\Z$ corresponds to periodic, and $\nu\in\Z+\half$ to antiperiodic eigenfunctions of~\ref{e:trigvalue}.
\end{remark}
Next,
define the space
$\ell^{2,p}(\Integer) := \{c\in\ell^2(\Integer):\sum_{n\in\Z}|n|^{p}|c_{n}|^2 < \infty\}$.
The requirement that $c\in \ell^{2,4}(\Integer)$ ensures $By\in L^2([0,2\pi])$.
The reason why this is the case is that $B$ is a second-order differential operator, 
which implies that the Fourier coefficients of $By$ will grow $n^2$ faster as $|n|\to\infty$ 
than those of $y$.
\begin{definition}
\label{d:tridiagonalspec}
(Eigenvalues of the tridiagonal operator)
Let the operator $B_{\nu}: \ell^2(\Z)\to\ell^2(\Z)$ be defined by~\eqref{e:trigoperatorBnu}.
The set of eigenvalues is given by
\be
\Sigma(B_\nu) := \{\lambda\in\C : \exists c\not\equiv 0 \in \ell^{2,4}(\Integer)\,\,\,{\rm s.t.}\,\,\, B_\nu c = \lambda c\}\,.
\label{e:trispec}
\ee
\end{definition} 

We have the following important result:
\begin{lemma} 
If $\nu\in\Z$ or $\Z+\half$, then $\Sigma_{\pm}(B) = \Sigma(B_{\nu})$,
and the dimension of the corresponding eigenspaces are equal,
respectively. 
\label{l:connectspec}
\end{lemma}
\textit{Proof.}
By standard results in Fourier analysis~\cite{ReedSimon} one defines the bijective linear map 
\be
U:\ell^2(\Z)\to L^2(\mathbb{S}^1)\,, \qquad
(Uc)(t) = \sum_{n\in\Z}c_{n}\e^{\i nt}\,,
\label{e:isometry}
\ee
and the multiplication operator
\be
M_\nu: L^2([0,2\pi]) \to L^2([0,2\pi])\,, \qquad 
(M_{\nu}w)(t) = \e^{\i\nu t}w(t)\,.
\label{e:multoperator}
\ee
By construction 
$B_\nu = (M_{\nu}U)^{-1}B M_{\nu}U$ in the standard basis 
and $U$, $M_{\nu}$ are unitary.
Also, $y = M_{\nu}Uc$ (see~\eqref{e:fourier}). 
Hence,
it follows $\Sigma_{\pm}(B)=\Sigma(B_{\nu})$ and the dimensions of the corresponding eigenspaces are equal.
\qed

\begin{remark}
The Floquet exponent~$\nu$ can be shifted by any integer amount without loss of generality, 
since doing so simply corresponds to a shift in the numbering of the Fourier coefficients in~\eqref{e:fourier}.
So, for example, 
$\nu\mapsto\nu+s$ simply corresponds to
$(\alpha_n,\beta_n,\gamma_n)\mapsto(\alpha_{n+s},\beta_{n+s},\gamma_{n+s})$ for all $n\in\Integer$.
\label{r:Floquetshift}
\end{remark}

\subsection{Reducible tridiagonal operators and ascending and descending Fourier series}
\label{s:Fourier_ascendingdescending}
We show that the tridiagonal operator $B_\nu$ is reducible.
Recall that a tridiagonal operator is reducible if there exists a zero element along the subdiagonal, or superdiagonal \cite{Horne}.
\begin{lemma}
\label{l:reducible}
If $A\in\Natural$ and $\nu\in\Z$ or $\Z+\half$, then $B_\nu$ is reducible.
\end{lemma}

\textit{Proof.}
There are two cases to consider: 
(i) $\nu\in\Z$, corresponding to periodic eigenvalues, 
and (ii) $\nu\in\Z + \half$, corresponding to antiperiodic eigenvalues. 
In either case, however, when $A\in\Natural$ one has
\vspace*{-1ex}
\bse
\label{e:fourierzeros1}
\begin{align}
\alpha_{n} &= 0 ~\iff~ n = \halfA + 1 - \nu ~~~\vee~~~ n = -\halfA + \half - \nu \,,
\\
\gamma_{n} &= 0 ~\iff~ n = \halfA - 1 - \nu ~~~\vee~~~ n = -\halfA - \half - \nu\,.
\end{align}
\ese
In both cases,
one can find two values of $n$ that make $\alpha_{n}$ and $\gamma_{n}$ zero, respectively, 
but only one of them is an integer, depending on whether $A$ is even or odd.
Note also that $\beta_n=0$ for $n= - \nu \pm A/2$, 
but the corresponding value of $n$ is integer only if $A$ is even and $\nu\in\Z$ or $A$ is odd and $\nu\in\Z+\half$.
(The equalities in~\eqref{e:fourierzeros1} hold for all $\nu\in\Real$, 
but only when $\nu\in\Z$ or $\nu\in\Z + \half$ do they yield integer values of~$n$.)
\qed

We emphasize that, when $A\notin\Natural$, a similar statement (namely, that $B_\nu$ is reducible)
can be made for different values of $\nu$.  
The precise values of $\nu$ can be immediately obtained from the definition of the coefficients 
$\alpha_n$, $\beta_n$ and $\gamma_n$ in~\eqref{e:alphabetagammadef}.
On the other hand, the particular significance of integer and half-integer values of $\nu$ is that they are associated with
periodic and antiperiodic eigenvalues, which are the endpoints of the spectral bands.
In Section~\ref{s:Heunconnection} we will also see how the periodic and antiperiodic eigenvalues are related to the solution of a connection problem for a particular Heun ODE.

Consider the tridiagonal operator $B_\nu: \ell^2(\Z) \to \ell^2(\Z)$ in~\eqref{e:trigoperatorBnu}.
Let $\ell^2_+ = \ell^2(\Natural_o)$ 
($\Natural_o:=\Natural\cup\{0\}$)
and
$\ell^2_- = \ell^2(\Z\setminus\Natural_o)$,
so that $\ell^2(\Z) = \ell^2_- \oplus \ell^2_+$,
and denote by $P_\pm$ orthogonal projectors from $\ell^2(\Integer)$ onto $\ell^2_\pm$ respectively.
Finally, introduce the block decomposition
\be
B_\nu = \begin{pmatrix} B_- & A_- \\ A_+ & B_+ \end{pmatrix}, 
\label{e:Bnu_pmdecomposition}
\ee
where the semi-infinite tridiagonal operators $B_\pm$ are defined as
\be
\label{e:Bpm1}
B_{-} :=
\begin{pmatrix*}
\ddots & \ddots & \ddots &   \\
  & \alpha_{-2} & \beta_{-2} & \gamma_{-2} \\
  &   & \alpha_{-1}  & \beta_{-1}
\end{pmatrix*}\,,
\hspace{9mm}
B_{+} :=
\begin{pmatrix*}
\beta_{0} & \gamma_{0} &  &  \\
\alpha_1 & \beta_1 & \gamma_1 &  \\
& \ddots & \ddots & \ddots  
\end{pmatrix*}\,,
\ee
and $A_\pm$ only have one nontrivial entry each, equal to $\alpha_{0}$ and $\gamma_{-1}$ respectively,  
in their upper right corner and lower left corner, respectively.
If $A\in\Natural$ and $\nu=(1-A)/2$ 
(corresponding to the case of periodic eigenvalues when $A$ is odd and antiperiodic eigenvalues when $A$ is even),
it is easy to see that $\alpha_0 = \gamma_{-1} = 0$ and therefore $A_\pm\equiv 0$, 
which implies that $B_\nu = B_- \oplus B_+$ and
$\ell_\pm^2$ are invariant subspaces of $B_\nu$.
The above considerations imply the following:

\begin{lemma}
\label{l:Bdecomposition_nu=-half-Ahalf}
If $A\in\Natural$ and $\nu=(1-A)/2$,
then $\Sigma(B_\nu)=\Sigma(B_-)\cup\Sigma(B_+)$,
where $B_\pm$ are given by~\eqref{e:Bpm1}.
\end{lemma}

The case $A\in\Natural$ and $\nu = A/2$ is similar, but more complicated. 
In this case, it is necessary to also introduce 
a second block decomposition of $B_\nu$ in addition to~\eqref{e:Bnu_pmdecomposition}, namely:
\be
B_\nu = \begin{pmatrix} \~B_- & \~A_- \\ \~A_+ & \~B_+ \end{pmatrix}, 
\ee
where
\be
\label{e:Bpm2}
\~B_{-} :=
\begin{pmatrix*}
\ddots & \ddots & \ddots &   \\
& \alpha_{-1} & \beta_{-1} & \gamma_{-1} \\
&   & \alpha_{0} & \beta_{0}
\end{pmatrix*}\,,
\hspace{9mm}
\~B_{+} :=
\begin{pmatrix*}
\beta_1 & \gamma_1 &  &   \\
\alpha_2 & \beta_2 & \gamma_2 & \\
 & \ddots & \ddots & \ddots  
\end{pmatrix*}\,,
\ee
and $\~A_\pm$ only have one nontrivial entry each, equal to $\alpha_{1}$ and $\gamma_0$ respectively,  
in their upper right corner and lower left corner, respectively.
If $A\in\Natural$ and $\nu=A/2$ 
(corresponding to the case of periodic eigenvalues when $A$ is even and antiperiodic eigenvalues when $A$ is odd),
it is easy to see that $\gamma_{-1} = \beta_0 = \alpha_1 = 0$ and therefore $A_- = \~A_+ \equiv 0$.
On the other hand, $A_+$ and $\~A_-$ are not identically zero.
Thus, $B_\nu$ cannot be split into a direct sum of two semi-infinite tridiagonal operators. 
Nevertheless, an analog of Lemma~\ref{l:Bdecomposition_nu=-half-Ahalf} still holds.

\begin{lemma} 
\label{l:Bdecomposition_nu=Ahalf}
If $A\in\Natural$ and $\nu=A/2$, then 
$\Sigma(B_\nu)=\Sigma(B_-)\cup\Sigma(B_+) = \Sigma(\~B_-)\cup\Sigma(\~B_+)$,
where ${B}_\pm$ and $\~B_\pm$ 
are given by~\eqref{e:Bpm1} and~\eqref{e:Bpm2}, respectively. 
\end{lemma}

\textit{Proof.}
We first show that $\Sigma(B_\nu)\subset\Sigma(B_-)\cup\Sigma(B_+)$.
Recall that $A_-=0$ but $A_+\ne0$.
Let $\l$ and $c$ be an eigenpair of $B_\nu$, and let $c_\pm = P_\pm c$,
so that $c = (c_-,c_+)^\T$.
If $c_-\ne0$, we have $B_-c_-=\l c_-$, and therefore $\l\in\Sigma(B_-)$.
Otherwise, $c_-=0$ implies $c_+\ne0$ and $c= (0,c_+)^\T$, and $B_+c_+ = \l c_+$, i.e., $\l\in\Sigma(B_+)$.

We show that $\Sigma(B_-)\cup\Sigma(B_+)\subset \Sigma(B_\nu)$.
Suppose that $\l$ and $c_+\ne0$ are an eigenpair of $B_+$, and let $c = (0,c_+)^\T$.
Then $B_\nu c = \l c$, implying $\l\in\Sigma(B_\nu)$.
Finally, suppose that $\l\in\Sigma(B_-)\setminus\Sigma(B_+)$, with associated eigenvector $c_-\ne0$.
In this case, let $c = (c_-,p)^\T$.
We choose $p$ such that $p = - (B_+-\l)^{-1}A_+c_-$.
One can show 
(similarly to Lemma~\ref{l:volkmer1}) 
that it is always possible to do so
since $B_+$ is closed with compact resolvent.
Therefore, the operator $(B_+-\l)^{-1}$ exists and is bounded, and 
$\l\notin\Sigma(B_+)$ implies that $\l$ is in the resolvent set of $B_+$.
But then we have $B_\nu c = \l c$, which implies $\l\in\Sigma(B_\nu)$.

The proof that $\Sigma(B_\nu)=\Sigma(\~B_-)\cup\Sigma(\~B_+)$ is entirely analogous,
but we report it because it is useful later. 
If $\l$ and $c$ are an eigenpair of $B_\nu$, let $\~c_\pm = \~P_\pm c$, 
with $\~P_\pm$ defined similarly as $P_\pm$.
If $\~c_+\ne0$, we have $\~B_+\~c_+ = \l\~c_+$ and therefore $\l\in\Sigma(\~B_+)$,
since $\~A_+=0$.
Otherwise, similar arguments as before show that $\~B_-\~c_-=\l\~c_-$ and therefore $\l\in\Sigma(\~B_-)$. 
We therefore have $\Sigma(B_\nu)\subset\Sigma(\~B_-)\cup\Sigma(\~B_+)$.
Finally, to show that $\Sigma(\~B_-)\cup\Sigma(\~B_+)\subset\Sigma(B_\nu)$,
we first observe that if $\l$ and $\~c_-$ are an eigenpair of $\~B_-$,
and $c = (\~c_-,0)^\T$, one has $B_\nu c = \l c$ and therefore $\l\in\Sigma(B_\nu)$.
Conversely, if $\l\in\Sigma(\~B_+)\setminus\Sigma(\~B_-)$, with eigenvector $\~c_+$,
it is always possible to choose $p$ such that $p=-(\tilde{B}_--\l)^{-1}\tilde{A}_-\tilde{c}_+$ 
(again, cf.\ Lemma~\ref{l:volkmer1}), 
and therefore 
$c = (p,\~c_+)^\T$ satisfies $B_\nu c = \l c$, implying $\l\in\Sigma(B_\nu)$.
\qed

\begin{remark}
\label{r:decomposition}
If $A\in\Natural$ and $\nu=A/2$, then $B_+$ and $\~B_-$ can be decomposed as
\be
B_+ = \begin{pmatrix} 0 & \gamma_0 \\ 0 & \~B_+ \end{pmatrix}\,,\qquad
\~B_- = \begin{pmatrix} B_- & 0 \\ \alpha_0 & 0 \end{pmatrix}.
\ee
\end{remark} 

\begin{corollary}
\label{c:furtherdecomp}
If $A\in\Natural$ and $\nu=A/2$, then 
$\Sigma(B_\nu)= \Sigma(B_-)\cup \Sigma(\~B_+)\cup \{0\}$.
\end{corollary}

\noindent
Importantly, the proofs of Lemmas~\ref{l:Bdecomposition_nu=-half-Ahalf} and~\ref{l:Bdecomposition_nu=Ahalf} also imply the following:
\begin{theorem}
\label{t:ascendingseries}
If $A\in\Natural$ and $\lambda\in\Complex\setminus\{0\}$ is a periodic or antiperiodic eigenvalue of the trigonometric operator~\eqref{e:inceoperator},
then there exists an associated eigenfunction generated by either an ascending or descending Fourier series.
\end{theorem}
\textit{Proof.}
The proof is trivial when $\nu = (1-A)/2$, since in this case $B_\nu = B_-\oplus B_+$.
On the other hand, the case $\nu = A/2$ requires more care. 
The proof of Lemma~\ref{l:Bdecomposition_nu=Ahalf} shows that, if 
$\l$ and $c_+$ are an eigenvpair of $B_+$, 
then $c=(0, c_+)^\T$ is a corresponding eigenvector of $B_\nu$.
Next,
if $\l$ and $\~c_-$ are an eigenvpair of $\~B_-$,
then $c=(\~c_-,0)^\T$ is a corresponding eigenvector of $B_\nu$.
Finally, note $\Sigma(\~B_-) = \Sigma(B_-)\cup\{0\}$.
\qed

\begin{corollary}
If $A\in\Natural$ and $\lambda\in\Complex$ is a periodic or antiperiodic eigenvalue with geometric multiplicity two, 
then a first eigenfunction can be written in terms of an ascending Fourier series, 
while a second linearly independent eigenfunction is given by a descending Fourier series.
\end{corollary}

\section{Transformation to a Heun ODE}
\label{s:heun}

We now introduce a final change of independent variable that maps the trigonometric ODE~\eqref{e:trigonometricODE} into Heun's equation.
All the results of Sections~\ref{s:Fourier2Heun} and~\ref{s:frobenius} below
will hold for integer as well as non-integer values of $A$ except where expressely indicated.
This further reformulation 
allows us to interpret the Dirac problem \eqref{e:Diraceigenvalueproblem} as a connection problem for Heun's ODE. 

\subsection{Transformation from the trigonometric ODE to Heun's equation}
\label{s:Fourier2Heun}

Recall that Heun's equation is a second-order linear ODE with four regular singular points \cite{erdelyi,Ince1956,Ronveaux}.
We first rewrite~\eqref{e:trigonometricODE} using Euler's formula.  
Then we perform the following change of independent variable:
\be
t \mapsto \z := \e^{\i t}\,.
\label{e:exp}
\ee
We then obtain a third reformulation of our spectral problem,
since the transformation~\eqref{e:exp} maps the trigonometric ODE~\eqref{e:trigonometricODE}
(and therefore \eqref{e:Diraceigenvalueproblem} with elliptic potential~\eqref{e:ellipticpotential})
into the following Heun ODE:
\be
\z^2F(\z;m)y_{\z\z} + \z G(\z;m)y_{\z} + H(\z;\lambda,A,m)y = 0\,,
\label{e:HeunODE}
\ee
where
\vspace*{-1ex}
\bse
\begin{gather}
F(\z;m) := -m\z^2 + (2m-4)\z - m\,,
\label{e:p1}
\\
G(\z;m) := -\txtfrac32 m \z^2 + (2m-4)\z - \half m\,,
\label{e:p2}
\\
H(\z;\lambda,A,m) := \txtfrac14 A(A+1)m\,\z^2 + \big(\lambda + A^2(1 - \halfk)\big)\,\z + \txtfrac14 A(A-1)m\,.
\label{e:p3}
\end{gather}
\ese
Note that the trigonometric ODE~\eqref{e:trigonometricODE} does not explicitly contain the Floquet exponent~$\nu$.
The role of $\nu$ for Heun's ODE will be played by the Frobenius exponents discussed below.

Equation~\eqref{e:HeunODE} has three regular singular points in the finite complex plane 
plus a regular singular point at infinity.  
Specifically,
in the finite complex plane one has a regular singular point at $\z=0$ and two additional regular singular points 
where $F(\z;m)=0$, i.e., when
\be
\z^2 - 2(1- \txtfrac{2}{m})\,\z + 1 = 0\,,
\label{e:Heunquadratic}
\ee
which is satisfied for 
\be
\label{e:singpoints}
\z_{1,2} = \frac{m-2 \pm 2\sqrt{1-m}}{m}\,.
\ee
Note that $\z_{1,2}<0$ for all $m\in(0,1)$, and
$\z_{2} = 1/\z_{1}$.
Without loss of generality, we take $|\z_{1}| < 1 < |\z_{2}|$.
Summarizing, the four real regular singular points are at $0,\z_1,\z_2,\infty$, with
$\z_2\in(-\infty,-1)$ and $\z_1\in(-1,0)$.

\begin{remark}
One can equivalently map the first-order ZS system \eqref{e:ZS} into a first-order
Heun system with the same four singular points 
using the same change of independent variable~\eqref{e:exp} [cf.\ Appendix~\ref{a:heun}]:
\be
\label{zs4}
\z w_\z= - \left[\frac A2 \s_3 + 
\begin{pmatrix*}
0  &\frac{1}{2}  \\
\frac{2\l\z}{4\z+m(\z-1)^2}  & \frac{(\z^2-1)m}{2(4\z+m(\z-1)^2)}  
\end{pmatrix*}\right]w\,,
\ee
where $w(\z;\lambda)=(w_1,w_2)^{\T}$.
\end{remark}

\subsection{Frobenius analysis of Heun's ODE}
\label{s:frobenius}

Next we apply the method of Frobenius to \eqref{e:HeunODE} at the regular singular points $\z=0$ and $\z=\infty$.
Then we construct half-infinite tridiagonal operators whose eigenvalues coincide with those of the tridiagonal operators
discussed in Section~\ref{s:Fourier_ascendingdescending}.
By direct calculation,
one can easily check that the Frobenius exponents of~\eqref{e:HeunODE} are as in Table~\ref{t:Frobenius}.
\begin{table}[b!]
\begin{center}
\begin{tabular}{|l|l|l|l|l|}
\hline
& $\z=0$ & $\z=\z_1$ & $\z=\z_2$ & $\z=\infty$ \\
\hline
$\rho_1$ & $\rho_1^o = A/2$ & $\rho_1^1=0$ & $\rho_1^2=0$ & $\rho_1^{\infty}=A/2$ \\
\hline
$\rho_2$ & $\rho_2^o=-(A-1)/2$ & $\rho_2^1=1/2$ & $\rho_2^2=-1/2$ & $\rho_2^{\infty}=-(A+1)/2$ \\
\hline
\end{tabular}
\end{center}
\caption{Frobenius exponents corresponding to the Heun ODE~\eqref{e:HeunODE}.}
\label{t:Frobenius}
\kern-2\bigskipamount
\end{table} 
The Frobenius exponents $\rho_{1,2}$ at $\z=0$ and $\z=\infty$ are obtained by looking for solutions of~\eqref{e:HeunODE} in the form
\vspace*{-1ex}
\bse
\label{e:frobeniusseries}
\be
y_o(\z;\lambda) = \z^{\rho}\sum_{n=0}^{\infty}c_{n}\,\z^{\,n}\,,
\label{e:frobenius1}
\ee
and
\vspace*{-1ex}
\be
y_\infty(\z;\lambda) = \z^{\rho}\sum_{n=0}^{\infty} c_{n}\,\z^{-n}\,,
\label{e:frobenius2}
\ee
\ese
respectively, with $c_0 \neq 0$ in each case.
Note, 
when $A$ is even, $\rho_1^o$ and $\rho_1^\infty$ are integer while $\rho_2^o$ and $\rho_2^\infty$ are half-integer,
and vice versa when $A$ is odd.
Note also that $\rho^{o}_2-\rho^{o}_1=1/2-A$ and $\rho_1^{\infty}-\rho_2^{\infty} = 1/2+A$,
so when $A\in\Natural$, these differences are never integer, and no exceptional cases (i.e., resonances) arise.

Next we study the three-term recurrence relations at $\z = 0$ and $\z = \infty$, since
they are key to proving the reality of the $\l$ eigenvalues.
We begin by plugging~\eqref{e:frobenius1} and~\eqref{e:frobenius2} into \eqref{e:HeunODE}.
The coefficients of the Frobenius series~\eqref{e:frobenius1} at $\z=0$ solve the following three-term recurrence relations.
For $\rho = \rho_1^o = A/2$:
\bse
\label{e:threetermrecurrence@z=0}
\begin{gather}
-\lambda c_0 + \txtfrac{m}{2}(2A+1)\,c_1 = 0\,, \qquad n=0\,,
\\
P_nc_{n-1} + (R_n - \lambda)c_n + S_nc_{n+1} = 0\,, \qquad n\ge 1\,,
\end{gather}
where
\vspace*{-1ex}
\begin{gather}
P_n = \txtfrac{m}{2}(n-1)(2A+2n-1)\,,~~
R_n = (1- \txtfrac{m}{2})((A+2n)^2-A^2)\,,~~
S_n = \txtfrac{m}{2}(n+1)(2A+2n+1)\,.
\end{gather}
\ese
For $\rho = \rho_2^o = -(A-1)/2$:
\bse
\label{e:heunrecur2}
\begin{gather}
\big[(\txtfrac{m}{2}-1)(2A-1)-\lambda\big]\,c_0 - \txtfrac{m}{2}(2A-3)\,c_1 = 0\,, \qquad n=0\,,
\\
\widetilde{P}_nc_{n-1} + (\widetilde{R}_n - \lambda)c_n + \widetilde{S}_nc_{n+1} = 0\,, \qquad n\ge 1\,,
\end{gather}
where
\vspace*{-1ex}
\begin{gather}
\label{e:heunrecurelements}
\widetilde{P}_n = -\txtfrac{m}{2}n(2A-2n+1)\,,~~
\widetilde{R}_n = (1-\txtfrac{m}{2})((2n+1-A)^2-A^2)\,,~~
\widetilde{S}_n = -\txtfrac{m}{2}(n+1)(2A-2n-3)\,.
\end{gather}
\ese
Similarly, 
the coefficients of the Frobenius series~\eqref{e:frobenius2} at $\zeta=\infty$ are given by the following three-term recurrence relations. For $\rho=\rho_{1}^{\infty} = A/2$:
\bse
\label{e:heunrecur3}
\begin{gather}
-\lambda c_0 - \txtfrac{m}{2}(2A-1)c_1 = 0\,, \qquad n=0\,,
\\
X_nc_{n-1} + (Y_n - \lambda)c_n + Z_nc_{n+1} = 0\,, \qquad n\ge 1\,,
\end{gather}
where
\vspace*{-1ex}
\begin{gather}
X_n = -\txtfrac{m}{2}(n-1)(2A-2n+1)\,,~~
Y_n = (1-\txtfrac{m}{2})((2n-A)^2-A^2)\,,~~
Z_n = -\txtfrac{m}{2}(n+1)(2A-2n-1)\,.
\end{gather}
\ese
For $\rho=\rho_{2}^{\infty} = -(A+1)/2$:
\bse
\label{e:heunrecur4}
\begin{gather}
\big[(1-\txtfrac{m}{2})(2A+1)-\lambda\big]\,c_0 + \txtfrac{m}{2}(2A+3)\,c_1 = 0\,, \qquad n=0\,,
\\
\widetilde{X}_nc_{n-1} + (\widetilde{Y}_n - \lambda)c_n + \widetilde{Z}_nc_{n+1} = 0\,, \qquad n\ge 1\,,
\end{gather}
where
\vspace*{-1ex}
\begin{align}
\widetilde{X}_n = \txtfrac{m}{2}n(2A+2n-1)\,,~~
\widetilde{Y}_n = (1-\txtfrac{m}{2})((2n+1+A)^2-A^2)\,,~~
\widetilde{Z}_n = \txtfrac{m}{2}(n+1)(2A+2n+3)\,.
\end{align}
\ese

\begin{remark}
\label{rem-eig-heun}
The three-term recurrence relations at $\zeta=0$ can be written as the eigenvalue problems
\be
T^{\pm}_{o}c = \lambda c\,,
\label{e:evalue01}
\ee
where $T_o^{\pm}:\ell^{2,4}(\N_o)\subset\ell^2(\N_o) \to \ell^2(\N_o)$,
and
\be
T^{-}_{o} := \begin{pmatrix*} R_0 & S_0 &  & \\ P_1 & R_1 & S_1 & \\  & \ddots & \ddots & \ddots \end{pmatrix*}\,,
\qquad
T^{+}_{o} := \begin{pmatrix*} \widetilde{R}_0 & \widetilde{S}_0 &  & \\ \widetilde{P}_1 & \widetilde{R}_1 & \widetilde{S}_1 & \\ & \ddots & \ddots & \ddots \end{pmatrix*}\,. 
\label{e:matrix01}
\ee
Similarly, the 
three-term recurrence relations at $\zeta=\infty$ can be written as the eigenvalue problems
\be
T^{\pm}_{\infty}c = \lambda c\,,
\label{e:evalue11}
\ee
where $T_{\infty}^{\pm}: \ell^{2,4}(\N_o)\subset\ell^2(\N_o)\to \ell^2(N_o)$,
and
\be
T^{-}_{\infty} := \begin{pmatrix*} Y_0 & Z_0 & & \\ X_1 & Y_1 & Z_1 & \\  & \ddots & \ddots & \ddots \end{pmatrix*}\,,
\qquad
T^{+}_{\infty} := \begin{pmatrix*} \widetilde{Y}_0 & \widetilde{Z}_0 & & \\ \widetilde{X}_1 & \widetilde{Y}_1 & \widetilde{Z}_1 & \\ & \ddots & \ddots & \ddots \end{pmatrix*}\,. 
\label{e:matrix11}
\ee
\end{remark}

\subsection{Relation between Fourier series and the connection problem for Heun's ODE}
\label{s:Heunconnection}

Recall that: 
(i)~If $\lambda\in\Complex$ is a periodic or antiperiodic eigenvalue of~\eqref{e:trigvalue},
one has $\nu\in\Z$ or $\nu \in \Z+\half$, respectively.
(ii)~The Floquet exponents can be shifted by an arbitrary integer amount
by shifting the indices of the Fourier coefficients
(cf.\ Remark~\ref{r:Floquetshift}).
(iii)~By Theorem~\ref{t:ascendingseries}, each periodic or antiperiodic eigenvalue has an associated
ascending or descending Fourier series when $A\in\Natural$.
(iv) The transformation $\z=\e^{\i t}$ maps the Frobenius series~\eqref{e:frobeniusseries} 
to ascending or descending Fourier series~\eqref{e:fourier},
and vice versa.
(v)~Finally, when $A\in\Natural$, 
the values of the Frobenius exponents for the expansions at $\z=0$ and at $\z=\infty$
are either integer or half-integer.

Moreover, the Floquet exponents $\nu = (1-A)/2$ and $\nu = A/2$ in Lemmas~\ref{l:Bdecomposition_nu=-half-Ahalf}
and \ref{l:Bdecomposition_nu=Ahalf} coincide exactly with the Frobenius exponents $\rho_2^o$ and $\rho_1^o$ at $\z=0$, respectively.
The Frobenius exponents $\rho_2^\infty$ and $\rho_1^\infty$ at $\z=\infty$ 
are also equivalent to the above Floquet exponents upon a shift of indices.
As a result, the recurrence relations 
\eqref{e:threetermrecurrence@z=0}, \eqref{e:heunrecur2}, \eqref{e:heunrecur3} and \eqref{e:heunrecur4} 
associated to the Frobenius series \eqref{e:frobeniusseries} of Heun's ODE~\eqref{e:HeunODE}
are equivalent to those associated to the Fourier series solutions of the trigonometric ODE~\eqref{e:trigonometricODE}.
More precisely:

\begin{lemma}
\label{l:fslemma}
If $A\in\Natural$ and $\lambda\in\Sigma(B_\nu)$ is either a periodic or antiperiodic eigenvalue,
(i.e., $\nu$ integer or half-integer, respectively) 
then the following identities map the recurrence relations generated by the Frobenius series solution of~\eqref{e:HeunODE} at $\z=0$ and $\z=\infty$ to the ascending and descending recurrence relations generated by the Fourier series solution of~\eqref{e:trigonometricODE}, respectively.
Namely:
\begin{itemize}
\item[(i)]
For $\nu=\rho^o_{1,2}$ one has, respectively:
\vspace*{-1ex}
\bse
\begin{gather}
(\alpha_n,\,\beta_n,\,\gamma_n) = (P_n,\,R_n,\,S_n)\,, \quad n\ge 0\,,
\\
(\alpha_n,\,\beta_n,\,\gamma_n) = (\widetilde{P}_n,\,\widetilde{R}_n,\,\widetilde{S}_n)\,, \quad n\ge 0\,.
\end{gather}
\ese
\item[(ii)]
For $\nu = \rho_{1,2}^\infty$ one has, respectively:
\vspace*{-1ex}
\bse
\begin{gather}
(\alpha_{-n},\,\beta_{-n},\,\gamma_{-n}) = (Z_{n},\,Y_{n},\,X_{n})\,, \quad n\ge 0\,,
\\
(\alpha_{-n-1},\,\beta_{-n-1},\,\gamma_{-n-1}) = (\widetilde{Z}_n,\,\widetilde{Y}_n,\,\widetilde{X}_n),, \qquad n\ge 0\,,
\end{gather}
\ese
\end{itemize}
\end{lemma}
\textit{Proof.}
When $\nu=\rho_{1,2}^o$, the result follows immediately by direct comparison.
Likewise when $\nu = \rho_1^\infty$.
Finally, when $\nu = \rho_2^{\infty}$ we can simply shift $\nu \mapsto \nu+ 1$, 
which sends $n \mapsto -n-1$.
\qed

\begin{corollary}
\label{r:Aevenodd_periodicantiperiodic}
If $A$ is odd, 
then the eigenvalues of $T_o^+$ and $T_\infty^+$ correspond to the periodic eigenvalues of the Dirac operator~\eqref{e:Diracoperator} and $T_o^-$ and $T_\infty^-$ to the antiperiodic ones, via the map $\lambda = z^2$.
Conversely, if $A$ is even, then the eigenvalues of $T_o^-$ and $T_\infty^-$ correspond to the periodic eigenvalues of the Dirac operator 
and those of $T_o^+$ and $T_\infty^+$ to the antiperiodic ones.
\end{corollary}

\begin{remark}
We emphasize that, when $A\in\Natural$, Lemma \ref{l:fslemma} only holds for periodic or antiperiodic 
solutions of~\eqref{e:trigonometricODE} (i.e., $\nu$ integer or half-integer). 
On the other hand, even when $A\notin\Natural$, a similar conclusion holds for certain Floquet solutions of~\eqref{e:trigonometricODE}. 
Namely, even for generic values of $A$, one can establish a one-to-one correspondence between certain 
Floquet exponents and ascending or descending Floquet eigenfunctions of~\eqref{e:trigonometricODE}, 
and in turn with Frobenius series solutions of~\eqref{e:HeunODE}.
\end{remark}

So far we have analyzed the properties of solutions corresponding to
periodic and antiperiodic eigenvalues of the problem.
We now turn to the question of identifying these eigenvalues.
Doing so yields the desired characterization of the Lax spectrum of \eqref{e:Diracoperator}. 

\begin{figure}[t!]
\centerline{\includegraphics[trim= 120 380 20 140,clip,width=0.6\textwidth]{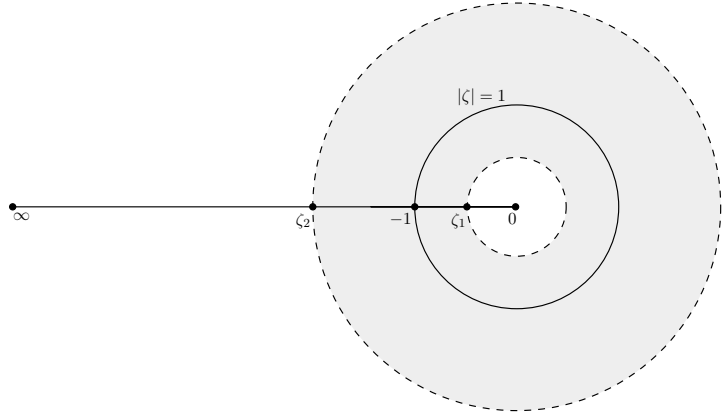}}
\caption{The singular points $\zeta = 0,\zeta_1,\zeta_2$, and $\infty$ and the region $|\zeta|=1$ in the complex $\zeta$ plane.}
\label{f:zetaplane}
\end{figure}

\begin{remark}
A periodic/antiperiodic eigenfunction of \eqref{e:Diraceigenvalueproblem} with potential \eqref{e:ellipticpotential} corresponds to a Fourier series solution \eqref{e:fourier} of the trigonometric ODE~\eqref{e:trigonometricODE} that is convergent for $t\in\Real$.
The transformation \eqref{e:exp} given by $\z = \e^{{\rm i} t}$,
which maps the real $t$-axis onto the unit circle $|\z|=1$
(cf.\ Fig.~\ref{f:zetaplane}),
maps these solutions
into a Laurent series representation for the solutions of Heun's ODE~\eqref{e:HeunODE}.
The question of identifying which solutions of Heun's ODE define periodic/antiperiodic eigenfunctions of \eqref{e:Diraceigenvalueproblem} is discussed next.
\end{remark}

\begin{lemma}
\label{lem-Heuns-e-vs}
Let $T_o$ be either one of the operators $T_o^\pm$ defined in Remark~\ref{rem-eig-heun} and 
let $y_o(\z) = \z^\rho w_o(\z)$ be a corresponding 
Frobenius series solution of Heun's equation
at $\z=0$.
Then:
\vspace*{-1ex}
\begin{enumerate}
\advance\itemsep-4pt
\item[(i)]
$\l$ is an eigenvalue of $T_o$ if and only if $w_o(\z)$ is analytic in the disk $|\z|<|\z_2|$;
i.e., if $y_o(\z)$ is analytic up to a branch cut when the Frobenius exponent $\rho_o$ is not integer.
\end{enumerate} 
Similarly, 
let $T_\infty$ be either one of the operators $T_{\infty}^\pm$ defined in Remark~\ref{rem-eig-heun} and 
let $y_\infty(\z) = \z^\rho w_\infty(\z)$ be a corresponding 
Frobenius series solution of Heun's equation at $\z=\infty$.  
Then:
\vspace*{-1ex}
\begin{enumerate} 
\advance\itemsep-4pt
\item[(ii)]
$\l$ is an eigenvalue of $T_\infty$ if and only if $w_\infty(\z)$ is analytic in the exterior disk $|\z|>|\z_1|$;
i.e., if $y_\infty(\z)$ is analytic up to a branch cut when the Frobenius exponent $\rho_\infty$ is not integer.
\end{enumerate}
\end{lemma}

\textit{Proof.}  
We consider~$T_o$ first.
The radius of convergence of the Frobenius series representing $y(\z)$ in a neighborhood of $\z=0$ is at least $|\z_1|$.  
Moreover, $\l\in\C$ is an eigenvalue of $T_o$ if and only if the corresponding eigenvector $c\in \ell^{2,4}(\Natural_o)$ (see Remark~\ref{rem-eig-heun}). 
Since the entries of $c$ coincide with the coefficients of the Frobenius power series representing $y(\z)$, we conclude that $\l\in\C$ is an eigenvalue of $T_o$ if and only if the radius of convergence of this series is at least one.  
In this case $y(\z)$ is analytic in the disk $|\z|<|\z_2|$ (up to a possible branch cut), 
since there are no singular points of the Heun's equation in the annulus $|\z_1|<|\z|<|\z_2|$.
The proof for $T_\infty$ follows along the same lines.
\qed

\begin{remark}
\label{rem-Heun-conn}
The above results relate the existence of eigenvalues to the connection problem for Heun's equation~\eqref{e:HeunODE}. 
For simplicity, consider the case of periodic eigenvalues.
Assume $A\in \N$.
The Frobenius analysis of Section~\ref{s:frobenius} yields two linearly independent solutions of Heun's ODE near each of the four singular points.
Let $y_{1,2}^o(\z;\l)$ be the Frobenius series with base point $\z=0$ and $y_{1,2}^1(\z;\l)$ those with base point $\z=\z_1$.
Both  $y_{1,2}^o(\z;\l)$ and $y_{1,2}^1(\z;\l)$ form a basis for the solutions of Heun's ODE~\eqref{e:HeunODE} 
in their respective domains of convergence.
Since these domains overlap, in the intersection region one can express one set of solutions in terms of the other, i.e., 
$(y_1^1,y_2^1) = (y_1^o,y_2^o)\,C$, with a constant non-singular connection matrix~$C$.
The Frobenius exponents at $\z=\z_1$ are 0 and $\half$,
and, when $A\in\N$, one of the Frobenius exponents at $\z=0$ is integer and the other is half-integer.
Therefore, the values of $\lambda$ for which the analytic solution at $\z=0$ converges up to $|\z|=|\z_2|$
are precisely those values for which the Frobenius series with integer exponent at $\z=0$ 
is exactly proportional to that with integer exponent at $\z=\z_1$.
Similar arguments hold for the solutions near $\z=\z_2$ and $\z=\infty$.
In other words, when $\l$ is a periodic eigenvalue, 
the analytic solutions at $\z=0$ and $\z=\z_1$ or those at $\z=\z_2$ and $\z=\infty$ must 
be proportional. 
This is the manifestation of an eigenvalue in terms of the connection problem for the Heun's equation~\eqref{e:HeunODE}.  
If both pairs of analytic solutions are proportional to each other, $\l$ is a double eigenvalue, otherwise $\l$ is a simple eigenvalue.
(In Section~\ref{s:realevalues1} we will also see that all positive eigenvalues have multiplicity two and all negative eigenvalues have multiplicity one.) 
Similar results hold 
for the antiperiodic eigenvalues once the square root branch cut resulting from the half-integer Frobenius exponent is taken into account.
\end{remark}

We also mention that there is an alternative but in a sense equivalent way to look at the problem, which is to study the convergence of the 
Frobenius series solutions~\eqref{e:frobeniusseries} using Perron's rule~\cite{Perron}.
This connection is briefly discussed in Appendix~\ref{s:perron}.

\section{Real eigenvalues of the operators $T_{o}^{\pm}$ and $T_{\infty}^{\pm}$}
\label{s:realevalues1} 

Thus far we have shown that the periodic and antiperiodic eigenvalues of \eqref{e:Diracoperator} with Jacobi elliptic potential~\eqref{e:ellipticpotential} 
and amplitude $A\in\Natural$ can be obtained from the eigenvalues of certain unbounded tridiagonal operators,
namely, 
$T_o^\pm$ and $T_\infty^\pm$ defined in Section~\ref{s:frobenius}.
We now prove that all eigenvalues of these operators are real.
We do so in two steps:
First, in Section~\ref{s:realHeuntruncation}, we show that finite truncations of these operators have purely real eigenvalues.
Then, in Section~\ref{s:conteigenvalues}, we use semicontinuity to show these operators have purely real eigenvalues. 

\subsection{Real eigenvalues of the truncated operators $T_{o,N}^{\pm}$ and $T_{\infty,N}^{\pm}$}
\label{s:realHeuntruncation}

Here we show that finite truncations of the operators $T_{o}^{\pm}$, $T_{\infty}^{\pm}$ have purely real eigenvalues.
We form the truncations by considering only the first $N-1$ terms of the corresponding three-term recurrence relations. 
To this end let $T^{\pm}_{o,N}$ and $T^{\pm}_{\infty,N}$ be the $N\times N$ truncations of $T^{\pm}_{o}$ and $T^{\pm}_{\infty}$,
respectively.

\begin{lemma}
If $A\in\Natural$ and $m\in(0,1)$, then for any $N>0$
the matrices $T^{-}_{o,N}$ and $T^{\pm}_{\infty,N}$ have purely real eigenvalues.
\label{l:truncationslemma}
\end{lemma}
 
The result is a consequence of the fact that 
$P_{n+1}S_{n}\ge 0$,
$X_{n+1}Z_{n}\ge 0$,
and $\widetilde{X}_{n+1}\widetilde{Z}_{n}>0$, 
$n\ge 0$,
which makes it possible to symmetrize $T^{-}_{o,N}$ and $T^{\pm}_{\infty,N}$
via a similarity transformation 
(see~\cite{Golub,Horne}).
The result does not apply to $T_{o,N}^{+}$, since 
there exists an $n>0$ such that $\widetilde{P}_{n+1}\widetilde{S}_{n} < 0$,
and, as a result,
some of the entries of the resulting symmetrized matrix would be complex.
Thus,
another approach is needed to show the eigenvalues of $T^{+}_{o,N}$ are all real.
To this end we introduce the following definition~\cite{Horne}:

\begin{definition}
(Irreducibly diagonally dominant)
An $N\times N$ tridiagonal matrix is irreducibly diagonally dominant if
(i) it is irreducible;
(ii) it is diagonally dominant, i.e.,
$|a_{ii}| \geq \sum_{j\neq i}|a_{ij}|$\,,
for all $i\in \{0,\dots,N-1\}$;
and
(iii)
there exists an $i\in\{0,\dots,N-1\}$ such that $|a_{ii}|>\sum_{j\ne i} |a_{ij}|$.
Here $a_{ij}$ denotes the entry in the $i$-th row and $j$-th column of the matrix.
\end{definition}
\begin{theorem}
\label{t:veselic}
(Veselic,~\cite{veselic}~p.~171)
Let $T_{o,N}^{+}$ be an $N\times N$ tridiagonal matrix which is irreducibly diagonally dominant and such that ${\rm sign}(\widetilde{P}_{n}\widetilde{S}_{n-1}) = {\rm sign}(\widetilde{R}_{n}\widetilde{R}_{n-1})$ for $n=1,\dots,N-1$.
Then $T_{o,N}^+$ has $N$ real simple eigenvalues.  
\end{theorem}

Next we show that $T_{o,N}^+$ satisfies the hypotheses of Theorem~\ref{t:veselic} and thus has only real eigenvalues.

\begin{lemma}
If $A\in\Natural$ and $m \in (0,1)$,
then for any $N>0$ all eigenvalues of $T_{o,N}^+$ are real and distinct.
\label{l:Veselic}
\end{lemma}

\textit{Proof.} 
First,
$A\in \Natural$ implies $\widetilde{P}_{n}\widetilde{S}_{n-1} \neq 0$ for $n \ge 1$.
Thus,
$T^{+}_{o,N}$ is irreducible. 
Next,
$\widetilde{R}_{n} < 0$ when $n \leq \lfloor A-\half \rfloor$. 
Similarly,
$\widetilde{R}_{n-1} < 0$ when $n \leq \lfloor A+\half \rfloor$.
(Here,
$\lfloor x \rfloor$ denotes the greatest integer less than or equal to $x$.)
Thus,
${\rm sign}(\widetilde{R}_{n} \widetilde{R}_{n-1}) < 0$ if and only if $A=n$.
Likewise,
$\widetilde{P}_{n} < 0$ when $n \leq \lfloor A +\half \rfloor$,
and $\widetilde{S}_{n-1} < 0$ when $n \leq \lfloor A-\half \rfloor$. 
Thus,
sign$(\widetilde{P}_{n}\widetilde{S}_{n-1}) < 0$ if and only if $A=n$.
Hence,
\be
{\rm sign}(\widetilde{P}_{n}\widetilde{S}_{n-1}) =
{\rm sign}(\widetilde{R}_{n}\widetilde{R}_{n-1})\,, \qquad n\ge 1\,.
\label{e:sign2}
\ee

Finally,
consider the transpose $(T^{+}_{o,N})^{\T}$. 
Note~\eqref{e:sign2} remains valid.
For $n = 0$ one easily gets $|\widetilde{R}_0| > |\widetilde{P}_{1}|$.
Moreover, for $n\ge1$, one has
$|\widetilde{R}_{n}| = (1-\halfk)(2n+1)|2n + 1 - 2A|$ and 
$|\widetilde{P}_{n+1}| + |\widetilde{S}_{n-1}| = \halfk(2n+1)|2n+1-2A|$.
Thus,
$|\widetilde{R}_{n}| > |\widetilde{P}_{n+1}| + |\widetilde{S}_{n-1}|$ for $n\ge 1$.
Hence
$(T^{+}_{o,N})^{\T}$ is an $N\times N$ irreducibly diagonally dominant tridiagonal matrix and satisfies~\eqref{e:sign2}.
The result follows from Theorem~\ref{t:veselic}. 
\qed

\begin{theorem}
If $A\in\Natural$ and $m \in (0,1)$,
then for any $N>0$ all eigenvalues of the tridiagonal matrices $T_{o,N}^\pm$ and $T_{\infty,N}^\pm$ are real and have geometric multiplicity one. 
Moreover, all eigenvalues of $T_{o,N}^+$ and $T_{\infty,N}^+$, and all nonzero eigenvalues of $T_{o,N}^{-}$ and $T_{\infty,N}^-$ are simple.
\label{t:realsimple_truncated_eigs}
\end{theorem}

\textit{Proof.} 
For $T_{o,N}^-$ and $T_{\infty,N}^\pm$,
the reality of all eigenvalues was proved in Lemma~\ref{l:truncationslemma},
and for $T_{o,N}^+$ it was proved in Lemma~\ref{l:Veselic}.
Moreover, 
Lemma~\ref{l:Veselic} also proved that the eigenvalues of $T_{o,N}^+$ are simple.

Let $\lambda$ be an eigenvalue of~$T_{o,N}^-$
and $c = (c_0,\dots,c_{N-1})^\T$ the corresponding eigenvector. 
Assume $c_0=0$. 
Then it follows from the three-term recurrence relation~\eqref{e:threetermrecurrence@z=0} that $c_n = 0$ for $n\ge 1$. 
(Note that $S_0$ is nonzero.)
Since $c$ is an eigenvector this is a contradiction.
Hence the first component of the eigenvector is necessarily nonzero.
Next, let $c$ and $\tilde{c}$ be two eigenvectors corresponding to the same eigenvalue of $T_{o,N}^{-}$. 
Consider the linear combination $b = \alpha c + \tilde{\alpha}\tilde{c}$. 
Then there exists $(\alpha,\tilde{\alpha})\ne 0$ such that $b_0=0$.
By the first part of the argument $b\equiv 0$.
Hence the eigenvectors $c$ and $\tilde{c}$ are linearly dependent.
The proofs for $T_{o,N}^+$ and $T_{\infty,N}^\pm$ are identical.
\qed 

\subsection{Generalized convergence and reality of periodic and antiperiodic eigenvalues}
\label{s:conteigenvalues}

In Section~\ref{s:realHeuntruncation} we showed that, for $A\in\Natural$, the $N\times N$ truncations of the tridiagonal operators have only real eigenvalues. 
It remains to show that the tridiagonal operators $T_o^\pm$ and $T_\infty^\pm$ also have only real eigenvalues. 
This result will follow from the fact that the eigenvalues of the tridiagonal operators possess certain continuity properties as the truncation parameter $N$ tends to infinity.
Some of the proofs in this section follow from Volkmer~\cite{volkmer3}.
For brevity we only present the details of the analysis for $T_{o}^{+}$.

\begin{lemma}
\label{l:assumption2}
Consider the operator $T_o^+$.
There exists $\theta \in (0,1)$ and $n_* \in \Natural$ such that
\be
2\max(\widetilde{P}_{n}^{\,2} + \widetilde{S}_{n}^{\,2},\,\,\, \widetilde{P}_{n+1}^{\,2} + \widetilde{S}_{n-1}^{\,2}) \leq \theta^2\widetilde{R}_{n}^{\,2}\,,\qquad n\ge n_*\,.
\label{e:assumption2}
\ee
The same estimate holds for the operators $T_o^-$ and $T_\infty^{\pm}$.
\end{lemma}
\textit{Proof.}
It follows from the definition of 
$\widetilde{P}_n$, 
$\widetilde{R}_n$, 
and $\widetilde{S}_n$ in~\eqref{e:heunrecurelements} that
\bse
\begin{gather}
\widetilde{P}_{n}^2 + \widetilde{S}_{n}^2 = 2m^2n^4(1+o(1))\,,\\
\widetilde{P}_{n+1}^2 + \widetilde{S}_{n-1}^2 = 2m^2n^4(1+o(1))\,,\\
\widetilde{R}_{n}^2 = (4-2m)^2n^4(1+o(1))\,,
\end{gather}
\ese
as $n\to\infty$.
Hence, 
let $\theta=m$. 
For $n$ sufficiently large one gets
$4m^2n^4 \leq \theta^2\widetilde{R}_{n}^{2} = m^2(4-2m)^2n^4$.
The result holds for $m\in(0,1)$.
It is easy to check that the same estimate holds also for the operators $T_o^-$, $T_\infty^{\pm}$.
\qed

Next we decompose $T_o^{\pm}$, $T_{\infty}^{\pm}$ into their diagonal and off-diagonal parts.
Namely, if $T$ is any one of the operators $T_o^{\pm}$, $T_\infty^{\pm}$, 
we write
\be
\label{e:decompose}
T := T_{D} + T_{O}\,,
\ee
where $T_D$ is the diagonal, and $T_O$ the off-diagonal. 
This decomposition is instrumental in proving the following:

\begin{lemma}
The operators $T_o^{\pm}$ and $T_\infty^{\pm}$ are closed with compact resolvent.
\label{l:volkmer1}
\end{lemma}

\textit{Proof.}
The proof follows closely that of the analogous result in~\cite{volkmer3}.
We provide details of the proof for the operator $T_o^+$. 
First,
by replacing $\widetilde{R}_{n}$ by $\widetilde{R}_{n} + \omega$ with sufficiently large $\omega$,
we may assume,
without loss of generality,
that $\widetilde{R}_{n} > 0$ and that~\eqref{e:assumption2} holds for all $n\ge 0$.
Then
\be
\|T_{O}c\| \leq \theta\|T_{D}c\| \quad \forall c\in \ell^{2,4}(\N_o)\,. 
\label{e:operatorbound}
\ee
Since $0 < \widetilde{R}_{n} \to \infty$ it follows
$T_D^{-1}$ exists and is a compact operator.
Moreover,
by~\eqref{e:operatorbound} it follows that
$\|T_O T_D^{-1}\|\leq\theta<1$.
Hence $T^{-1} = T_D^{-1}(I + T_OT_D^{-1})^{-1}$ is a compact operator (see~\cite{kato} p. 196).
Therefore, 
$T^{+}_{o}$ is a closed operator with compact resolvent. 
The proofs for $T_o^{-}$, $T_\infty^{\pm}$ are identical.
\qed

The proof of the next lemma is identical to that of Theorem~\ref{t:realsimple_truncated_eigs}.

\begin{lemma}
\label{l:geom1}
All eigenvalues of the operators $T_o^\pm$ and $T_\infty^\pm$ have geometric multiplicity one. 
\end{lemma}

Next we begin to address reality of the eigenvalues.
By Lemma~\ref{l:volkmer1} it follows $\Sigma(T_o^{\pm})$ and $\Sigma(T_\infty^{\pm})$ 
are comprised of a set of discrete eigenvalues with finite multiplicities.
Recall that,
for $z\in\Sigma(L)$, we have
$|\Im z| < \|q\|_{\infty}$.
Moreover, we also have
$|\Im z|\,|\Re z| \leq \half\|q_x\|_{\infty}$
for any $z\in\Sigma(L)$
(see \cite{BOT2020}).
Hence,
by the correspondence between the Dirac and Hill equations (see Section~\ref{s:Hill}) we have
\vspace*{-1ex}
\be
\Re\lambda \geq -\|q\|_{\infty}^2\,,
\qquad
|\Im\lambda| \leq \|q_x\|_{\infty}\,.
\ee
Thus there exists a curve $\mathcal{C}$ such that the region in the complex $\lambda$-plane bounded by $\mathcal{C}$ contains finitely many periodic (resp. antiperiodic) eigenvalues of Hill's equation with complex elliptic potential~\eqref{e:reduction1} counting multiplicity.
This suggests to apply the concept of generalized convergence of closed linear operators 
(see Appendix~\ref{a:genconv} for a discussion of generalized convergence).
In particular, we will use the following result:
\begin{theorem}
\label{t:kato1}
(Kato,~\cite{kato} p.~206)
Let $T$, $T_{n} \in \mathfrak{C}(\mathcal{X},\mathcal{Y})$, 
$n=1,2,\dots$ the space of closed operators between Banach spaces.
If $T^{-1}$ exists and belongs to $\mathfrak{B}(\mathcal{X},\mathcal{Y})$,
the space of bounded operators,
then $T_{n} \to T$ in the generalized sense if and only if $T_{n}^{-1}$ exists and is bounded for sufficiently large $n$ and $\|T_{n}^{-1} - T^{-1}\| \to 0$.
\end{theorem}

Theorem~\ref{t:kato1} implies the semicontinuity of a finite system of eigenvalues counted according to multiplicity (\cite{kato} p.~213).
To this end, 
we introduce a sequence of tridiagonal operators:
\be
T_n := T_D + P_nT_O\,,
\label{e:projection}
\ee
where $P_n$ is the orthogonal projection of $\ell^2(\Natural_o)$ onto ${\rm span}\{e_0,e_1,\dots,e_{n-1}\}$,
with 
$\{e_i\}_{n\in\Natural_o}$ being the canonical basis. 
Thus,
for example,
$T_n$ is determined by,
say, 
$T_o^+$ with the off-diagonal entries $\widetilde P_{j}$, $\widetilde S_{j}$
replaced by zeros for $j\ge n$.
Clearly,
$\Sigma(T_N) = \Sigma(T_{o,N}^+)\cup\{\widetilde{R}_n\}_{n\ge N}$ for any $N\in\N$.
The following result is obtained from Theorem 2 in \cite{volkmer3},
the difference being the additional zero column for $T_o^-$ and $T_{\infty}^-$.
Once the first column and row are deleted, the proof is identical.
We therefore omit the proof for brevity.

\begin{lemma}
Let $T$ be any one of the operators $T_o^{\pm}$, $T_\infty^{\pm}$.
If $T_n$ is defined by~\eqref{e:projection}, with $T_D$ and $T_O$ defined by~\eqref{e:decompose}, 
then $T_{n} \to T$ in the generalized sense (see Theorem \ref{t:kato1} above).
\label{l:genconv}
\end{lemma}

Using convergence in the generalized sense, 
we are now ready to
show that $T_o^{\pm}$ and $T_\infty^{\pm}$ have real eigenvalues only:
\begin{lemma}
\label{l:truncateconvergence}
If $\lambda_{n} \in \Sigma(T_o^{\pm})$
or $\lambda_{n} \in \Sigma(T_\infty^{\pm})$,
then $\lambda_{n} \in \R$.
\end{lemma}
\textit{Proof.}
Let $T$ be any one of the operators $T_{o}^{\pm}$ or $T_{\infty}^{\pm}$. 
Count eigenvalues according to their multiplicity.
Fix $n\in\Natural$, 
and let $\lambda_n\in\Sigma(T)$. 
Let $\ep>0$ and $C_\ep := \{\lambda\in\Complex : |\lambda-\lambda_n|=\ep\}$.
Since $\| T_{n}^{-1} - T^{-1} \| \to 0$ as $n\to\infty$, 
we know for each $\delta > 0$ there exists $N \in \Natural$ such that  
$\| T_{N}^{-1} - T^{-1}\|<\delta$.
By semicontinuity of a finite system of eigenvalues (see~\cite{kato} p.~212), 
we can choose $\delta > 0$ such that 
$C_\ep$ contains an eigenvalue of $T_{N}$.  
Call this eigenvalue $\lambda_N$.
Since $\epsilon$ is arbitrary,
and $\lambda_N$ is real for any $N$ it follows $\lambda_n\in\R$.
\qed

Summarizing,
we have shown that the periodic (resp. antiperiodic) eigenvalue problems for \eqref{e:Diracoperator} with Jacobi elliptic potential~\eqref{e:ellipticpotential} can be mapped to eigenvalue problems for four tridiagonal operators obtained from a Frobenius analysis of the Heun equation~\eqref{e:HeunODE}.
Moreover,
all eigenvalues of the tridiagonal operators are real with geometric multiplicity one.
Putting everything together, we are now ready to prove the first part of Theorem~\ref{t:mainresult}:

\begin{theorem}
\label{t:mainresult_p1}
Consider \eqref{e:Diracoperator} with potential
\eqref{e:ellipticpotential} and $m\in(0,1)$.
If $A\in\Natural$, then 
$\Sigma(L;A,m) \,\subset \,\Real \cup (-\rm{i} \textit{A}, \rm{i} \textit{A})$,
and $q$ is a finite-band potential.
\end{theorem}

\textit{Proof.}
Let $z\in\Sigma_{\pm}(L)$.
Recall that we have established a direct correspondence between the periodic (resp. antiperiodic) eigenvalues of the tridiagonal operators $T_{o}^{\pm}$, $T_{\infty}^{\pm}$ and the periodic (resp. antiperiodic) eigenvalues of the Dirac operator~\eqref{e:Diracoperator} with elliptic potential~\eqref{e:ellipticpotential}.
Also,
$\Sigma(T_o^{\pm})\cup\Sigma(T_{\infty}^{\pm}) \subset \R$.
Hence, by Lemma~\ref{l:ZSHill_spectrum_equal}, 
and since $\lambda = z^2$, it follows that $\Sigma(L;A,m)\subset \Real\cup (-\i A,\i A)$ (see also Lemma~\ref{l:combinedlemma}).
Thus,
the periodic (resp. antiperiodic) eigenvalues of the Dirac operator~\eqref{e:Diracoperator} with elliptic potential~\eqref{e:ellipticpotential} are real or purely imaginary.
Then, by symmetry (see Lemmas~\ref{l:quartets} and~\ref{l:combinedlemma}), the entire Lax spectrum is only real and purely imaginary.
Finally,
that $q$ is finite-band for all $A\in\Natural$ and $m\in(0,1)$ follows from 
Lemma~\ref{l:combinedlemma}.
\qed

\begin{lemma}
\label{l:posevalues}
If $\lambda \in \Sigma(T_o^-) \cup \Sigma(T_{\infty}^{+})$,
then $\lambda \geq 0$. 
\end{lemma}

\textit{Proof.}
All entries of the tridiagonal operator $T_{\infty}^{+}$ are positive.
Consider the truncation $T_{\infty,N}^{+}$.
Without loss of generality take the transpose.
A simple calculation shows that $(T_{\infty,N}^{+})^\T$ is strictly diagonally dominant.
Hence, 
by the Gershgorin circle theorem all eigenvalues of $(T_{\infty,N}^{+})^\T$ are strictly positive.
By semicontinuity, 
in the limit $N\to \infty$ it follows that $\Sigma(T_{\infty}^{+}) \subset [0,\infty)$.

Next,
note that the first column of $T_o^{-}$ is comprised of all zeros.
Thus,
$\Sigma(T_{o}^{-}) = \Sigma(\tilde{T}_{o}^{-}) \cup \{0\}$,
where $\tilde{T}_o^{-}$ is defined by $T_o^-$ with the first row and the first column removed.
Moreover, $\tilde{T}_{o}^{-}$ has strictly positive entries,
and the transpose is diagonally dominant.
Arguing as in the previous case gives the result.
\qed 

\section{Lax spectrum for non-integer values of $A$}
\label{s:noninteger}

All of the results in this work up to the Fourier series expansion and the three-term recurrence relation
in Section~\ref{s:recurrence} hold independently of whether or not $A$ is integer.
The same holds for the Frobenius analysis in Section~\ref{s:frobenius}.
On the other hand, the reducibility of the tridiagonal operator $B_\nu$ with integer and half-integer
Floquet exponents $\nu$ in Section~\ref{s:Fourier_ascendingdescending} only holds when 
$A\in\Natural$ (because it is only in that case that zeros appear in the upper and lower diagonal entries).
Similarly, 
the Frobenius exponents at $\z=0$ and $\z=\infty$ in Section~\ref{s:frobenius} are integer or half-integer
only when $A\in\Natural$.
We next show that these are not just technical difficulties, but instead reflect a fundamental difference
in the properties of the Lax spectrum of \eqref{e:Diracoperator} when $A\notin\Natural$.

\begin{lemma}
If $A\notin\Natural$, and $m\in(0,1)$, then all periodic and antiperiodic eigenvalues of~\eqref{e:Diracoperator} with Jacobi elliptic potential~\eqref{e:ellipticpotential} have geometric multiplicity one.
\label{l:Anoninteger_nodoublepoints}
\end{lemma}

\textit{Proof.}
The proof proceeds by contradiction.
For simplicity, 
we focus on the periodic eigenvalues.
Suppose that for $A\notin\N$ and $\nu\in\Integer$ there exist two linearly independent eigenfunctions.
Then the transformation $\z = \e^{\i t}$ yields two linearly independent solutions of 
Heun's ODE~\eqref{e:HeunODE} on $|\z|=1$.
Let us denote these solutions as $\^y_1(\z;\l)$ and $\^y_2(\z;\l)$.
Note all points on $|\z|=1$ are ordinary points for Heun's ODE and, 
therefore, 
both $\^y_1(\z;\l)$ and $\^y_2(\z;\l)$ are analytic and single-valued in the annulus 
$|\z_1|<|\z|<|\z_2|$
(cf.~Fig.~\ref{f:zetaplane}). 
Moreover, recall that the Frobenius exponents at $\z=\z_1$ are $\rho_1^1=0$ and $\rho_2^1 =1/2$.  
Let $y_1^1(\z;\l)$ and $y_2^1(\z;\l)$ denote the corresponding solutions.  
Since $\^y_1(\z;\l)$ and $\^y_2(\z;\l)$ are linearly independent solutions, we have
$y_1^1(\z;\l) = c_1\^y_1(\z;\l) + c_2\^y_2(\z;\l)$ for some constants $c_1$ and $c_2$.
Then $y_1^1(\z;\l)$ is analytic and single valued in the region $0<|\z|<|\z_2|$

On the other hand, $y_1^1(\z;\l)$ is a linear combination of the Frobenius solutions $\~y_1(\z;\l)$ and $\~y_2(\z;\l)$ defined at the singular poinr $\z=0$, with  Frobenius exponents $\rho_1^o = A/2$ and $\rho_2^o = (1-A)/2$ respectively,
neither of which is an integer.
Thus, no single-valued solution can exist around $\z=0$.
Therefore, there cannot be two linearly independent periodic eigenfunctions.
Similar considerations apply for the antiperiodic eigenvalues.
\qed

Note Lemma~\ref{l:Anoninteger_nodoublepoints} does not hold for $m=0$, as in the limit $m\to 0$
Heun's equation~\eqref{e:HeunODE} degenerates into a Cauchy-Euler equation (with two regular singular points at $\z=0$ and $\z=\infty$).
Still, together with Lemma~\ref{l:dim2}, Lemma~\ref{l:Anoninteger_nodoublepoints} implies:

\begin{corollary}
If $A\notin\Natural$, and $m\in(0,1)$, then $\Sigma_\pm(L) \cap \Real = \emptyset$.
\label{c:norealperiodic}
\end{corollary}

In turn, since both $\Sigma_\pm(L)$ are infinite (see~\cite{Mityagin}), 
and since the periodic and antiperiodic eigenvalues are the endpoints of the spectral bands,
Corollary~\ref{c:norealperiodic} directly implies:
\begin{corollary}
If $A\notin\N$, and $m\in(0,1)$, then $\Sigma(L)$ with Jacobi elliptic potential~\eqref{e:ellipticpotential} has an infinite number of spines along the real $z$-axis. 
\end{corollary}

We conclude that when $A\notin\Natural$, the potential $q$ in~\eqref{e:ellipticpotential}
is not finite-band according to Definition~\ref{d:bandsandgaps},
which proves the only if part of Theorem~\ref{t:mainresult},
namely that $A\in\Natural$ is not only sufficient, but also necessary in order for $q$ in \eqref{e:ellipticpotential} to be finite-band,
as well as Theorem~\ref{t:Anoninteger}.

\section{Further characterization of the spectrum and determination of the genus}
\label{s:genus}

It remains to prove the last part of Theorem~\ref{t:mainresult}, namely the determination of the genus.
To this end, we need a more precise characterization of the Lax spectrum for $A\in\Natural$,
which will also yield the proof of the remaining parts of Theorem~\ref{t:Ainteger}.
We turn to this task in this section.

\subsection{Multiplicity of imaginary eigenvalues}

\begin{theorem}
\label{t:z-mult-1} 
If $z\in (-{\rm i}A,{\rm i}A)\setminus \{0\}$ is a periodic or an antiperiodic eigenvalue of~\eqref{e:Diracoperator}
with potential \eqref{e:ellipticpotential} with $A\in\N$,
and $m\in(0,1)$, then it has geometric multiplicity one.
\end{theorem}
\textit{Proof.}
By Lemma~\ref{l:ZSHill_spectrum_equal}
it follows $z\in\Sigma_{\pm}(L)$ if and only if $\lambda=z^2\in\Sigma_{\pm}(H^-)$,
respectively.
Moreover,
for $z\ne 0$ the geometric multiplicity of the periodic (resp. antiperiodic) eigenvalues is the same.
Next, by the results of Section \ref{s:heun},
each periodic (resp. antiperiodic) eigenfunction of $H^-$ is associated with an eigenvector of $T_o^\pm$ or $T_\infty^\pm$,
and
\be
\Sigma(T^{\pm}_o) \cup \Sigma(T^{\pm}_{\infty}) = \{\lambda=z^2:z\in\Sigma_\pm(L)\}\,.
\ee
with $T_o^+$, $T_{\infty}^+$ yielding periodic eigenvalues and $T_o^-$, $T_{\infty}^-$ antiperiodic eigenvalues when $A$ is odd,
and vice versa when $A$ is even
(cf.\ Corollary~\ref{r:Aevenodd_periodicantiperiodic}).
By Lemma~\ref{l:geom1}, each eigenvalue of $T_o^\pm$, $T_\infty^\pm$
has geometric multiplicity one.
Therefore, a periodic (resp. antiperiodic) eigenvalue $z\in\Complex$ of $L$ can have geometric multiplicity two if and only if $\lambda = z^2$ 
is simultaneously an eigenvalue of both $T_o^+$ and $T_\infty^+$ or 
simultaneously an eigenvalue of both $T_o^-$ and $T_\infty^-$.
On the other hand, Lemma~\ref{l:posevalues} showed that the eigenvalues of $T_o^-$ and $T_\infty^+$ are non-negative.
Hence, by the relation $\lambda=z^2$
all periodic (resp. antiperiodic) eigenvalues $z\in(-\i A,\i A)\setminus\{0\}$ of \eqref{e:Diracoperator} with potential~\eqref{e:ellipticpotential} have geometric multiplicity one.
\qed

\begin{corollary}
\label{cor-FlDir}
For all $m\in(0,1)$, 
if $z\in (-{\rm i}A,{\rm i}A)\setminus \{0\}$ is a periodic or an antiperiodic eigenvalue of \eqref{e:Diracoperator} with potential \eqref{e:ellipticpotential}, then $s(z)\ne 0$. 
\end{corollary}
\textit{Proof.}
Recall that $s(z)$ is defined by  \eqref{e:Msym}.
If $s(z)=0$, the monodromy matrix $M(z)$ would be diagonal, but this would imply the existence of two periodic (resp. antiperiodic)
eigenfunctions, which would contradict Theorem \ref{t:z-mult-1}.
\qed

Note that the above results do not hold for $m=0$ (a constant background potential), 
since in that case all periodic and antiperiodic eigenvalues except $z=\pm\i A$ have geometric multiplicity two.

\subsection{Dirichlet eigenvalues and behavior of the Floquet discriminant near the origin}
\label{s:z=0behavior}

In this subsection we prove some technical but important results that will be used later in the  proof of Theorem~\ref{t:Ainteger}. 

As in Section~\ref{s:m=0&m=1},
here it will be convenient to explicitly keep track of the dependence on $m$ 
by writing the potential, fundamental matrix solution, and monodromy matrix respectively as
$q(x;m)$, $\Phi(x;z,m)$ and $M(z;m)$.
We begin by recalling some relevant information.
We will use the structure of the monodromy matrix $M(z;m)=\Phi(2K(m);z,m)$ introduced in \eqref{e:Msym}.
Also recall that, when $m=0$ (in which case $q(x,0) \equiv A$),
$M(z,0)$ is given by~\eqref{e:M_m=0}. 
(Recall that $l=2K(m)$ is the (real) period of $\dn(x;m)$, and $2K(0)=\pi$.)
Thus,~\eqref{e:Msym} implies
\be
\label{D,sm=0}
\D(z;0)=\cos\big(\sqrt{z^2+A^2}{\pi}\big)\,, \quad 
s(z;0)=-\frac{A}{\sqrt{z^2+A^2}}\sin\big(\sqrt{z^2+A^2}{\pi}\big)\,.
\ee
Recall from Section~\ref{s:generalproperties}
that $\D(z;m)$ and $s(z;m)$ are even functions of $z$ while $c(z;m)$ is an odd function of $z$.
Let $\D_{j}(m)$, $-\i c_{j}(m)$ and $s_{j}(m)$ denote, respectively, 
the coefficients of $z^{2j}$, $z^{2j+1}$ and $z^{2j}$ in the Taylor series of $\D(z;m)$, $c(z;m)$ and $s(z;m)$ around $z=0$.
Combining \eqref{e:M_m=0} and~\eqref{e:v_amx2}, we obtain the following expansions near $z=0$:
\bse
\label{exp-0}
\begin{align}
&\D(z;m)=(-1)^A +\D_1(m)z^2 + O(z^4)\,, \\
&c(z;m)=-\i c_0(m)z - \i c_1(m)z^3 + O(z^5)\,,\\ 
&s(z;m)=s_1(m)z^2+ O(z^4)\,.
\end{align}
\ese

We want to study in detail the  behavior of $\D(z;m)$ near $z=0$.  
We begin by looking at the dynamics of (closed) gaps as a function of $A$ at $m=0$,
to show how the number of bands grows as $A$ increases.
According to \eqref{D,sm=0}, the periodic and antiperiodic eigenvalues are, respectively,
\be\label{z_n-form}
z_n=\pm\sqrt{4n^2-A^2}\,,\qquad 
z_n=\pm\sqrt{(2n+1)^2-A^2}\,,\qquad n\in\Z\,.
\ee
It follows that $z=0$ is a  periodic or antiperiodic eigenvalue when $A\in\Z$. 
Direct calculations show that
\be\label{Del''m=0}
\D_{zz}(z;0)=-\frac{\pi^2z^2}{z^2+A^2}\cos(\sqrt{z^2+A^2}\pi)-\frac{\pi A^2}{(z^2+A^2)^{\frac 32}}\sin(\sqrt{z^2+A^2}\pi),
\ee
so that 
\be\label{Del''zm=0}
\D_{zz}(0;0)=-\frac{\pi\sin (A\pi)}{A}.
\ee

Observe that, as a function of~$A$, $\D_{zz}(0;0)$ changes sign as $A$ passes through an integer value. 
For example, if $A$ passes through an even value $n\in\Natural$, the sign of $\D_{zz}(0;0)$ changes from ``+'' to ``$-$'',
corresponding to the transition of a pair of critical points of $\D(z;0)$ from $\R$ to $[-\i A,\i A]$. 
Correspondingly, a pair of zero level curves of $\Im \D(z)=0$
intersecting $\R$ transversally will pass through $z=0$ and intersect $[-\i A,\i A]$ forming an extra closed gap on $[-\i A,\i A]$.
This is the mechanism of increase of the number of gaps on $[-\i A,\i A]$. 
Note that \eqref{Del''m=0} implies that $\D_{z}(z;0)$ has a third order zero at $z=0$ when $A\in\Z$.
Next we show that this mechanism works for any $m\in (0,1)$.  
This will be accomplished through several intermediate steps.

\begin{lemma}
\label{l:del_zz=0}
For fixed  $A\in\N$ we have $\D_{zz}(0;0)=0$, 
and $(-1)^A\D_{zz}(0;m)$ is a strictly monotonically decreasing function of~$m$ for $m\in [0,1)$.
\end{lemma}
\textit{Proof.}
The first statement follows from \eqref{Del''m=0}. The rest of the proof is devoted to show that 
\be
(-1)^A\D_{zz}(0;m)<0\,,
\ee
when $m\in (0,1)$.
Substitution of \eqref{exp-0} into \eqref{det-eq} yields
\be\label{Dzz_0}
\half\D_{zz}(0;m)=\D_1(m)= \half(-1)^{A+1}{c_0^2(m)}\,.
\ee
Note that $(-1)^A\D_{zz}(0;m)\leq 0$
since $M(z;m)$ is real on $z\in[-\i A,\i A]$. 
Thus, it remains to show that $c_0(m)\not = 0$ for $m\in (0,1)$.
In fact, we will show below that $c_0(m)$ is monotonically increasing on $m\in [0,1)$. 
That, combined with $c_0(0)=0$ (see \eqref{Dzz_0}), will complete the proof.

Recall that $\Phi=\Phi(x;z,m)$ is the solution of the ZS system \eqref{e:ZS} normalized as $\Phi(0;z,m) \equiv \1$. 
Differentiating  \eqref{e:ZS} with respect to $z$  we get the system
\be\label{Phi-zs1}
\Phi_{xz}=(-\i z\s_3 + \i q\s_2)\Phi_z - \i\s_3\Phi\,.
\ee
Considering system \eqref{Phi-zs1} as a  non-homogeneous ZS system
[i.e., treating the term $-\i\sigma_3\Phi$ as a ``forcing''] and integrating, we obtain the solution
\begin{gather}
\Phi_z(x;z,m)=-\i\Phi(x;z,m)\int_0^x \Phi^{-1}(\x;z,m)\s_3\Phi(\x;z,m)\,\d\x\,, 
\qquad z\in\Complex\,.		
\end{gather}
Also recall that the (real) period of $q$ in~\eqref{e:ellipticpotential} is $l = 2K(m)$
and that $M(z;m)=\Phi(l;z,m)$.
By Lemma \ref{l:zequalzero}, $A\in\N$ implies $\Phi(l;0,m)\equiv(-1)^A\1$.
At $z=0$, we therefore have
\be
M_z(0;m)=-\i(-1)^A\int_0^l \Phi^{-1}(\x)\s_3\Phi(\x)\,\d\x=-\i(-1)^A\int_0^l
\begin{pmatrix*}
u_1v_2 +u_2v_1 & 2v_1v_2 \\
-2u_1u_2  & -u_1v_2 -u_2v_1 
\end{pmatrix*}\d\xi\,.	
\ee
where we introduced the notation
\be
\Phi(x;z,m) = \begin{pmatrix*}
u_1 & v_1 \\
u_2  & v_2 
\end{pmatrix*}\,,
\label{e:Phi_u1u2v1v2def}
\ee
which we will use extensively below.
On the other hand, in light of \eqref{e:Msym} we have
\be\label{det-sym'}
M_z=\D_z \1+c_z\s_3-\i s_z\s_2\,,
\qquad z\in\Complex\,,
\ee
which implies 
\be
c_z(0;m)=-\i(-1)^A\int_0^l(u_1v_2 +u_2v_1)\,\d x\,, \quad s_z(0;m)=2\i(-1)^A\int_0^lu_1u_2\,\d x=2\i(-1)^A\int_0^lv_1v_2\,\d x\,.
\ee
Comparing this with the expansion \eqref{exp-0} we then have 
\be
c_0(m) = \i c_z(0;m) = (-1)^A\int_0^l(u_1v_2 +u_2v_1)\,\d x\,.
\ee
At $z=0$, according to Section~\ref{a:exactsolution}, we also have 
\be
\Phi(x;0,m)=\cos(A\am x)\,\1+\sin (A\am x)\,\i\s_2\,, 
\ee
where  $\am x=\am(x;m)$,  
so that
\be
\Phi^{-1}(x;0,m)\s_3\Phi(x;0,m)
= \cos(2A\am x)\s_3+ \sin (2A\am x)\s_1\,.
\ee
We obtain
\be\label{c_0m}
c_0(m)=(-1)^A\int_0^{2K(m)}\cos(2A\am x)\,\d x
=(-1)^A\int_0^\pi \frac{\cos 2Ay\,\d y}{\sqrt{1-m\sin^2y}}\,,
\ee
where we used $y=\am x$, $\d y=\sqrt{1-m\sin^2y}\,\d x$. 
Now, from \cite{BirdFriedman}, 806.01, for $m\in(0,1)$ we have
\be\label{BF-int}
(-1)^A\int_0^\pi \frac{\cos 2Ay\, \d y}{\sqrt{1-m\sin^2y}}=
\pi\sum_{j=A}^\infty\frac{[(2j-1)!]^2m^j}{4^{2j-1}(j-A)!(j+A)![(j-1)!]^2}>0\, 
\ee
since all the coefficients of the convergent Taylor series are positive.
\qed

\begin{corollary}
\label{cor-z=0}
One has $|\D(z;m)|>1$ in a deleted neighborhood of $z=0$ on $(-{\rm i}A,{\rm i}A)$ for any $A\in\Natural$ and $m\in(0,1)$.
Moreover, $z=0$ is a simple critical point of $\D(z;m)$ and $\D(0;m)=(-1)^A$. 
\end{corollary}

\begin{remark}
Corollary \ref{cor-z=0} shows that no critical points of $\D(z;m)$ can move from $\R$ to $\i\R$ when we vary $m\in(0,1)$ with a fixed $A\in\N$.  
Similarly to the case $m=0$, the change of genus in the case $m>0$,
happens when we vary $A$ (see also Figure \ref{f:spectrum_A}).
\end{remark}

Next, 
recall that the monodromy matrix $M(z;x_o)$ normalized at a base point $x_o$ is given by~\eqref{e:M(x_o)}.
Using \eqref{e:M(x_o)}, we prove the following lemma regarding Dirichlet eigenvalues.

\begin{lemma}\label{lem-Dir} Let $A\in\N$ and $m\in (0,1)$.
If an open gap $\g$ on $(-{\rm i}A,{\rm i}A)$ contains a zero of $s(z)$, then the associated Dirichlet eigenvalue is movable. 
\end{lemma}

\textit{Proof.} 
Equations~\eqref{e:Msym} and~\eqref{e:M(x_o)} and direct calculation show that
\be\label{Mx0a}
M(z;x_0) = 
\begin{pmatrix*}
\D+c(u_1v_2 +u_2v_1)+s(u_1u_2+v_1v_2) & -s(u_1^2+v_1^2)-2cu_1v_1 \\
s(u_2^2+v_2^2)+2cu_2v_2  & \D-c(u_1v_2 +u_2v_1)-s(u_1u_2+v_1v_2)
\end{pmatrix*}\,,
\ee
where $c=c(z)$ and $s=s(z)$ were defined in \eqref{e:Msym} and 
the functions $u_1$, $u_2$, $v_1$ and $v_2$, defined as in~\eqref{e:Phi_u1u2v1v2def}, are evaluated at $x = x_0$. 

Consider first an open gap $\g\subset\i\R$ that does not contain $z=0$.
From \eqref{det-eq} it follows that $c(z)\neq  0$ on~$\g$. 
Note $u_1(x_0)v_1(x_0)\neq 0$  for small $x_0>0$ 
[because $u_1(0)=1$ and $v_1(0)=0$ and $u_1$ and $v_1$ are analytic in $x$ as solutions of~\eqref{e:ZS} with 
the potential~\eqref{e:ellipticpotential}].
Therefore, it follows from \eqref{Mx0a} that 
$M_{12}(z;x_0)=0$ implies $s(z)\neq 0$.
But $M_{12}(z;0)=s(z)$.
Thus, each Dirichlet eigenvalue in such a gap is movable.

Consider now a gap $\g_0\subset \i\R$ containing $z=0$,
i.e., the \textit{central gap}. 
By Corollary \ref{cor-z=0}, such gap exists for any $m\in(0,1)$,
and by~\eqref{D,sm=0}, it does not exist when $m=0$.  
We consider $m\in(0,1)$.
Then by Lemma \ref{l:del_zz=0} (see~\eqref{Dzz_0}), 
$c(z)$ has a simple zero at $z=0$ and,
by~\eqref{exp-0}, $s(z)$ has at least a double zero at $z=0$. 
Thus, the condition $M_{12}(z;x_0)=0$ near $z=0$ becomes
\be
\label{M12=0}
-s_1z^2(u_1^2(x_0;z)+v_1^2(x_0;z)) +2\i c_0zu_1(x_0;z)v_1(x_0;z) = R(x_0;z)\,,
\ee
where $R(x_0;z)\in\R$ when $z\in \i\R$ and $R(x_0;z)=O(z^3)$ uniformly in small real $x_0$.
By Lemma~\ref{l:del_zz=0} (see~\eqref{Dzz_0}) we have $c_0\neq 0$. 
If $s_1\neq 0$,
\eqref{M12=0} shows that $M_{12}(z;x_0)$ 
has one fixed zero at $z=0$ whereas the location of the second zero depends on $x_0$ and is given by
\be
\label{z-Dir}
z = \frac{2\i c_0u_1(x_0;z)v_1(x_0;z)-\frac{R(x_0;z)}{z}}{s_1(u_1^2(x_0;z)+v_1^2(x_0;z))}=
 \frac{2\i c_0u_1(x_0;0)v_1(x_0;0)+O(z)}{s_1(u_1^2(x_0;0)+v_1^2(x_0;0)+O(z))}
\in \i\R\,, 
\ee
which is a point inside the central gap on $(-\i A,\i A)$. 
Indeed,
the requirement $\det M(z;x_0)\equiv 1$ and \eqref{Mx0a} imply that a Dirichlet eigenvalue can not be in the interior of any band located on $(-\i A,\i A)\setminus \{0\}$.

Equations~\eqref{M12=0} and~\eqref{z-Dir} 
show that a zero of $M_{12}(z;x_0)$ in the gap $\g_0\subset(-\i A,\i A)\setminus\{0\}$ is 
always fixed at $z=0$, and therefore corresponds to an immovable Dirichlet eigenvalue, whereas a second zero is located at a point changing with $x_0$, and is therefore a movable Dirichlet eigenvalue. 
Indeed, the point $z=z(x_0)$ defined by~\eqref{z-Dir} attains $z(0)=0$ and $z(x_0)\neq 0$ at least for small $x_0>0$ since $v_1(0)=0$  and $v_1(x_0)\neq 0$ in a deleted neighborhood of zero.

Finally, if $s_1=\cdots =s_{k-1} = 0$ and $s_k\ne0$, with $k>1$, the leading-order portion of each term in the 1,2 entry of~\eqref{M12=0} yields instead
\be\label{M12=0-H}
-s_kz^{2k}(u_1^2+v_1^2) + 2\i c_0zu_1v_1 = R\,,
\ee
where again $R = O(z^3)$ is real-valued for $z\in\i\Real$ and where for brevity we dropped the arguments. 
Repeating the same arguments as for \eqref{M12=0}, we see 
 that at least one of the roots in   \eqref{M12=0-H} is purely imaginary.
\qed  

\begin{remark}
\label{rem-Dir}
It follows from \eqref{M12=0} and~\eqref{M12=0-H} that $s_1(m)\neq 0$
if and only if there is excatly one movable Dirichlet eigenvalue  in a vicinity of $z=0$ for small $x_0\in\R$.
\end{remark}

\subsection{Proof of the remaining statements of Theorem~\ref{t:Ainteger}}
\label{S: Proofs of Thm 1.3 and 1.5}

Lemma \ref{lem-Dcs} and Theorem \ref{t:z-mult-1}
proves items 4 and~5 of Theorem~\ref{t:Ainteger}.
Theorem \ref{t:mainresult} together with Lemma \ref{r:asymp} and the symmetries \eqref{Dsc} implies item 3. 
Thus, it remains to prove items 2 and 6 only, namely:

\begin{theorem}
\label{t:Ainteger2}
Consider \eqref{e:Diracoperator} with Jacobi elliptic potential \eqref{e:ellipticpotential}.
If $A\in\Natural$ then: 
\begin{enumerate}
\advance\itemsep-4pt
\item 
For any $m\in(0,1)$,
there are exactly $2A$ symmetric bands of $\Sigma(L;A,m)$ on $(-\rm{i}\textit{A},\rm{i}\textit{A})$ separated by $2A-1$ symmetric
open gaps. 
The central gap (i.e., the gap surrounding the origin) contains an eigenvalue at $z=0$. 
This eigenvalue is periodic when $A$ is even and antiperiodic when $A$ is odd.
\item
Each of the open $2A-1$ gaps on  $(-\rm{i}\textit{A},\rm{i}\textit{A})$ contains exactly one movable Dirichlet eigenvalue.
Thus, all of the $2A-1$ movable Dirichlet eigenvalues of the finite-band solution are located in the gaps of the interval $(-\rm{i}\textit{A},\rm{i}\textit{A})$.
\end{enumerate}
\end{theorem}

\textit{Proof.}
The idea of the proof is based on continuous deformation of the elliptic parameter $m$, starting from $m=0$ and going into $m\in(0,1)$. 
The proof is based on the following three main steps, each of which will be discussed more fully below:

1. \textit{Analysis of the spectrum for $m=0$.} 
When $m=0$, $\dn(x,0)\equiv 1$, and the ZS system \eqref{e:Diraceigenvalueproblem} has a simple solution $\Phi(x;z,m)$. 
The monodromy matrix $M(z;m)$
based on $\Phi(x;z,m)$
was given explicitly in~\eqref{e:M_m=0} for $m=0$, 
and the Lax spectrum is $\Sigma(L)= \Real\cup[-\i A,\i A]$ in this case.
In particular, the vertical segment $[-\i A,\i A]$ is a single band that contains $2A-1$ double periodic/antiperiodic eigenvalues,
which we consider as being closed gaps.
Each of these closed gaps contains a Dirichlet eigenvalue (a zero of $s(z;m)$,
see~\eqref{D,sm=0}), 
which for $m=0$ is immovable according to Lemma~\ref{lem-Dcs} (see also~\eqref{Mx0a}).

2. \textit{Analysis of the spectrum for small nonzero values of $m$.} 
Corollary \ref{cor-z=0} states that
for all $m\in(0,1)$  the double eigenvalue at $z=0$ is embedded in the central gap $\g_0\subset(-\i A,\i A)$. 
Moreover, 
Corollary~\ref{cor-z=0} and Lemma \ref{lem-Dir} 
show that there is at least one movable Dirichlet eigenvalue on $\g_0$. Next, we show that
under a small deformation $m>0$ all the remaining closed gaps on $(-\i A,\i A)$ must open, creating $2A$ bands and $2A-1$ gaps on $(-\i A,\i A)$, with each gap containing exactly one movable Dirichlet eigenvalue.  Our proof of this statement is based on the fact that any periodic/antiperiodic eigenvalue on $(-\i A,\i A)\setminus\{0\}$ has geometric multiplicity one 
(see Theorem \ref{t:z-mult-1}), 
whereas  a double eigenvalue at a closed gap would have  geometric multiplicity two. 
By a continuity argument, each gap on $(-\i A,\i A)\setminus\{0\}$ has exactly one movable Dirichlet eigenvalue.
Thus, items 2,6 are proved for small $m>0$.

3. \textit{Control of the spectrum for arbitrary values of $m\in(0,1)$.}
We finally prove that the number $2A$ of separate bands on $(-\i A,\i A)$, as well as the fact that each gap on $(-\i A,\i A)$ contains exactly one movable Dirichlet eigenvalue,  cannot change when $m$ varies on $(0,1)$.

In what follows, we prove all the statements in items 1--3 above.

1. From \eqref{D,sm=0} it follows that, for $m=0$, we have
$\Sigma(L)=\R\cup[-\i A,\i A] $ with 
\be
z_n^2=n^2-A^2, ~~n=0,\pm 1,\dots,\pm A
\ee
being (interlaced) periodic or antiperiodic eigenvalues on $[-\i A,\i A]$.  
Note that $\D_{z}(z;0)$ has a simple zero at each $z_n$ for $n\not=0,\pm A$.
Therefore, each $z_n\ne 0,\pm \i A$ identifies a closed gap. 
On the other hand, since $\D_{z}(\pm \i A;0)\not=0$, the endpoints $\pm \i A$ of the spectrum are simple periodic eigenvalues. 
[By Lemma~\ref{l:combinedlemma}, when $m>0$ it follows $\pm \i A\not\in\Sigma(L)$ .]
Finally, note that 
$\D_{z}(z;0)$ has a double zero at $z=0$, which will be relevant in the discussion of item 3 below.

Recall that the Dirichlet eigenvalues are the zeros of $s(z;0)$.
By~\eqref{D,sm=0},
each closed gap on  $[-\i A,\i A]\setminus \{0\}$  contains exactly one Dirichlet eigenvalue $\mu_n$, $n=\pm 1,\dots, \pm (A-1)$, which is a simple zero of $s(z;0)$.  
Note also that there are exactly $2A-1$ closed gaps on $(-\i A,\i A)$,
and at each such gap with the exception of $z=0$ there is exactly one zero level curve $\Gamma_n$ of $\Im \D$ 
orthogonally crossing $\i\R$.
Moreover, by Lemma~\ref{l:del_zz=0}, there are eight zero level curves of $\Im \D$ passing through $z=0$, 
including the real and imaginary axes
(e.g., see Fig.~\ref{f:spectrum_A}, upper right panel).
Note that, by Lemma \ref{lem-Dir}, a Dirichlet eigenvalue $\mu_n$ becomes movable if the closed gap opens up as $m$ 
is deformed away from $m=0$.

2. Recall that the monodromy matrix $M(z;m)$ is entire in $z$ and $A$ and analytic in $m\in[0,1)$ 
(cf.\ Lemma~\ref{l:combinedlemma}). 
Also recall that $s(z;m)$ is real-valued on $\i\R$.
Finally, recall that by Lemma~\ref{lem-Dcs} zeros of $s(z;m)$ cannot lie in the interior of a band. 
Let $\mu_n$ be the zeros of $s(z;0)$ for $z\in\i\Real$, i.e., the Dirichlet eigenvalues along the imaginary axis when $m=0$. 
Since the zeros of $s(z;m)$ are isolated, for sufficiently small $m>0$, 
each Dirichlet eigenvalue $\mu_n$ must remain on $(-\i A,\i A)$ by continuity.
Thus, for sufficiently small values of~$m$, 
all the gaps on $(-\i A,\i A)\setminus\{0\}$ (independently of whether they are open or closed) must survive the small~$m$ deformation, with exactly one Dirichlet eigenvalue in each gap. 

Importantly, the above arguments imply that, for small $m\in(0,1)$, all the gaps on 
$(-\i A,\i A)\setminus\{0\}$ must be open.
Indeed,
the assumption that for small $m\in(0,1)$ there exists a closed gap at $z_*\in (-\i A,\i A)\setminus\{0\}$ leads to a contradiction, because by Corollary \ref{cor-FlDir},
$s(z;m)\ne0$ at the endpoints of each band.
There are $2A-2$ such open gaps. 
By continuity, each of them contains a zero of $s(z;m)$ and therefore a movable Dirichlet eigenvalue by Lemma \ref{lem-Dir}. 
Moreover, by continuity, $s(z;m)$ must have opposite signs at the endpoints of any gap in $(-\i A,\i A)\setminus\{0\}$.

Next, recall that by Theorem~\ref{t:mainresult_p1}, $\Sigma(L)\subset \R\cup (-\i A,\i A)$ for all $m\in(0,1)$.
It follows from Corollary~\ref{cor-z=0} that, for all $m\in(0,1)$, 
the (double) Floquet eigenvalue $z=0$ is immersed in a gap $\g_o\subset (-\i A,\i A)$. 
Then, by Lemma~\ref{lem-Dir}, for small $m>0$ there are exactly $2A-1$ open gaps on 
$(-\i A,\i A)$,  
with  each gap containing a movable Dirichlet eigenvalue.   
Therefore there are $2A$ (disjoint) bands on $(-\i A,\i A)$.

Finally, differentiation of $s(z;m)$ in~\eqref{D,sm=0} yields
\vspace*{-1ex}
\be
s_z(z;0)=\frac{zA}{z^2+A^2}\left[  \frac{\sin\left(\sqrt{z^2+A^2}{\pi}\right)}{\sqrt{z^2+A^2}}-\cos\left(\sqrt{z^2+A^2}{\pi}\right) \right],
\ee
which shows that $s_z(z;0)$ has a simple zero at the origin.
Therefore, $s_1(0)=s_{zz}(0;0)\neq 0$ and so, by Remark \ref{rem-Dir}, there is a unique movable Dirichlet eigenvalue in a vicinity of $z=0$ and it is situated on $\g_o$.
Hence there is exactly one movable Dirichlet eigenvalue in each gap implying that the genus of the corresponding Riemann surface in $2A-1$.
Thus, items 2 and 6 are proved for small $m>0$.

3. 
It remains to  prove that the number of bands on $(-\i A,\i A)$ and the number of movable Dirichlet eigenvalues
(which were established for small $m\in(0,1)$ in item 2 above) do not change as $m$ varies in $(0,1)$. 
Let us consider the deformation of the collection of bands (with genus $2A-1$) established for small $m\in(0,1)$.  
A possible  change of the genus can be caused  only by one of the following four possibilities:
(a) a collapse of a band into a point; 
(b) a splitting of a band into two or more separate bands; 
(c) a splitting of a gap into two or more separate gaps; 
(d) a collapse of an open gap into a closed one.

We next prove that none of these possibilities can occur. 
Indeed, regarding~(a), the collapse of a band into a point would contradict the analyticity of $\D(z;m)$,
since it would imply that the same value of $z$ is simultaneously a periodic and antiperiodic eigenvalue.
(Note that each band along $(-\i A,\i A)$ must necessarily start at a periodic eigenvalue and end at an antiperiodic one or vice versa,
since otherwise there would necessarily be a critical point $z_o$ inside the band. 
But a critical point $z_o$ inside the band would imply the existence of a second band emanating transversally
from the imaginary axis, contradicting Theorem~\ref{t:mainresult_p1}.)  
Similarly, regarding~(b), 
the splitting of a band would require a critical point of $\D(z;m)$ at some $z_o$ inside the band. 
But, again, a critical point at $z_o$ would mean that there is a zero-level curve of $\Im \D$ crossing $\i\R$ at $z_0$,
which in turn would contradict Theorem~\ref{t:mainresult_p1}.  

For the same reasons we have $\D_z(z_*;m)\neq 0$ at any non-periodic and non-antiperiodic Floquet eigenvalue $z_*$, 
separating a band and a gap on $(-\i A,\i A)$.  
Indeed, the contrary would lead to $\D_{z}$ having an even-order zero at $z_*$.
But that would imply at least two pairs of zero-level curves emanating from $\i\R$ at $z_*$ and, thus, again would 
contradict Theorem~\ref{t:mainresult_p1}.
 
We now turn our attention to the gaps, and specifically to the possibility~(c) listed above.
The splitting of a gap into two or more separate gaps 
would imply that $\D(z;m) $ has a local minimum $z_0$ on the gap at some $m\in(0,1)$ and simultaneously $|\D(z_0;m)|\leq  1$. 
That, again, would contradict Theorem~\ref{t:mainresult_p1}.  

Thus, it remains to exclude possibility~(d), namely the collapse of a gap.   
By Lemma \ref{lem-Dir}, the central gap $\g_o$ containing $z=0$ stays open for any $m\in(0,1)$. 
Also, for small $m\in(0,1)$, it was shown in the proof of item~2 above that the signs of $s(z;m)$ at the endpoints of any gap $\g\subset (-\i A,\i A)\setminus\{0\}$ are opposite. 
These signs cannot change in the course of a deformation with respect to $m\in(0,1)$, 
by Corollary~\ref{cor-FlDir}.
Thus, each gap that was open for small  $m\in(0,1)$ must  contain a zero of $s(z;m)$ and, therefore, as it was proven in item~2 above, 
must stay open for all $m\in(0,1)$. 
So, the genus $2A-1$ is preserved for all $m\in(0,1)$. 
Moreover, each gap contains a zero of $s(z;m)$ and thus, according to Lemma~\ref{lem-Dir}, a movable Dirichlet eigenvalue.

So, we proved that  each  gap on $(-\i A,\i A)$ contains exactly one movable Dirichlet eigenvalue.
It is well known~\cite{ForestLee,GW1998,mclaughlinoverman} that the number of movable Dirichlet eigenvalues is equal to the genus $2A-1$.
This, completes the proof of Theorem~\ref{t:Ainteger2} for all $m\in(0,1)$.
\qed
      
\begin{remark}
It follows from Remark \ref{rem-Dir}  that $s_1(m)\neq 0$ for  all $m\in(0,1)$. 
\end{remark}

\section{Dynamics of the spectrum as a function of $A$ and $m$}
\label{s:plots}

\begin{figure}[t!]
\centerline{\qquad\begin{minipage}[t!]{0.5\textwidth}
\kern-1.1\bigskipamount
\includegraphics[width=45mm]{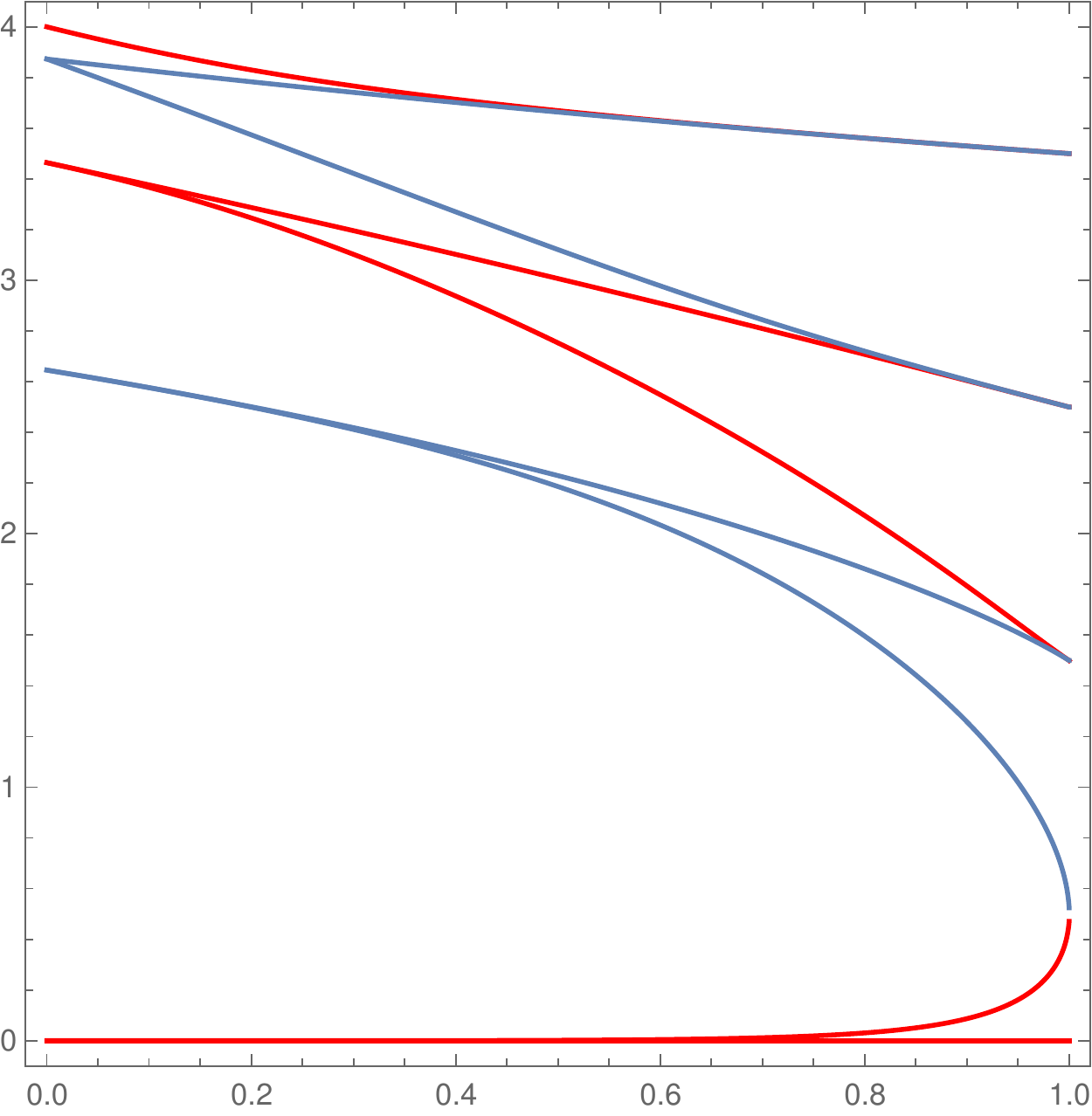}\\[0.2ex]
\hglue-0.25em\includegraphics[width=45.75mm]{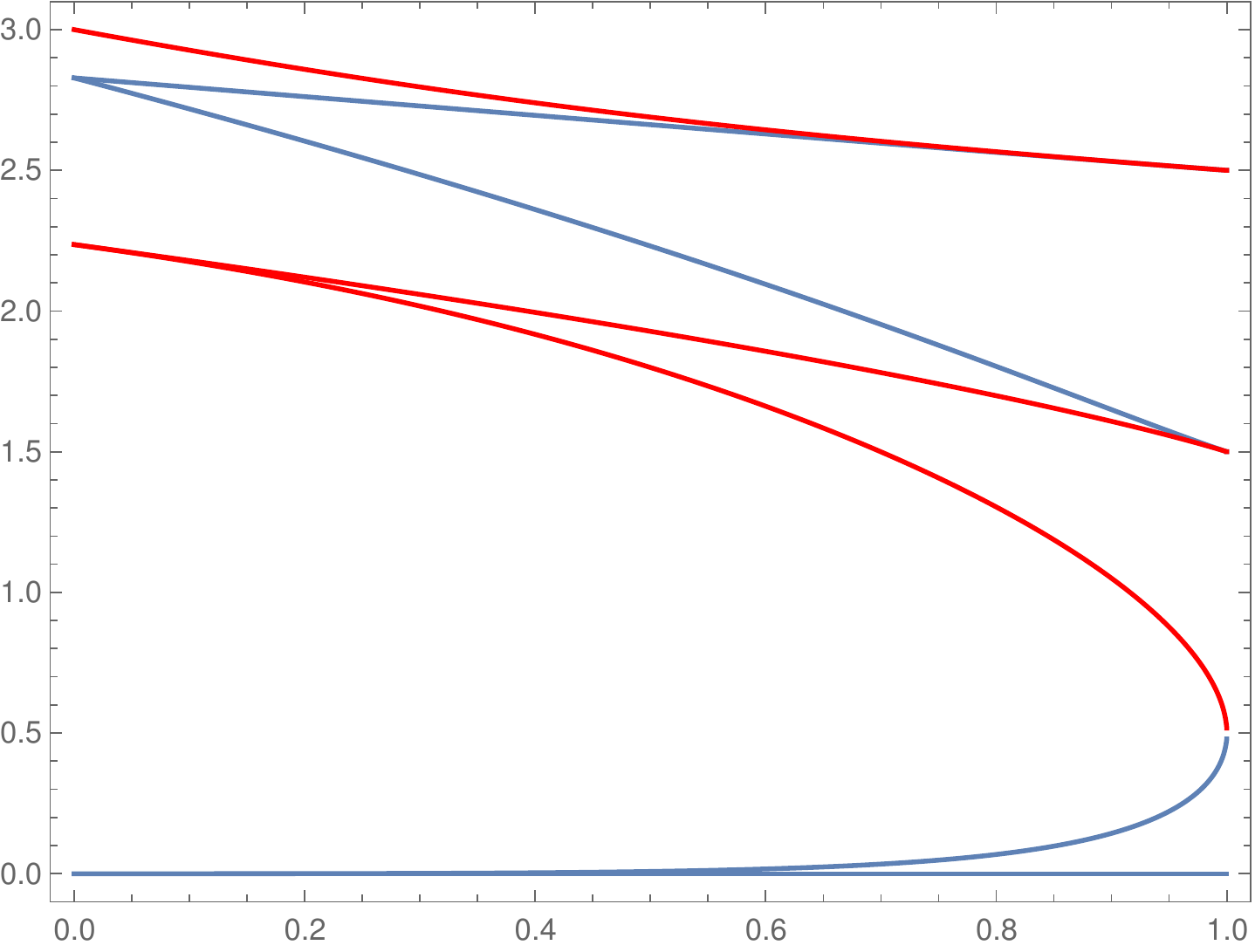}
\end{minipage}\kern-8em%
\begin{minipage}[t!]{0.5\textwidth}
\includegraphics[width=46.25mm]{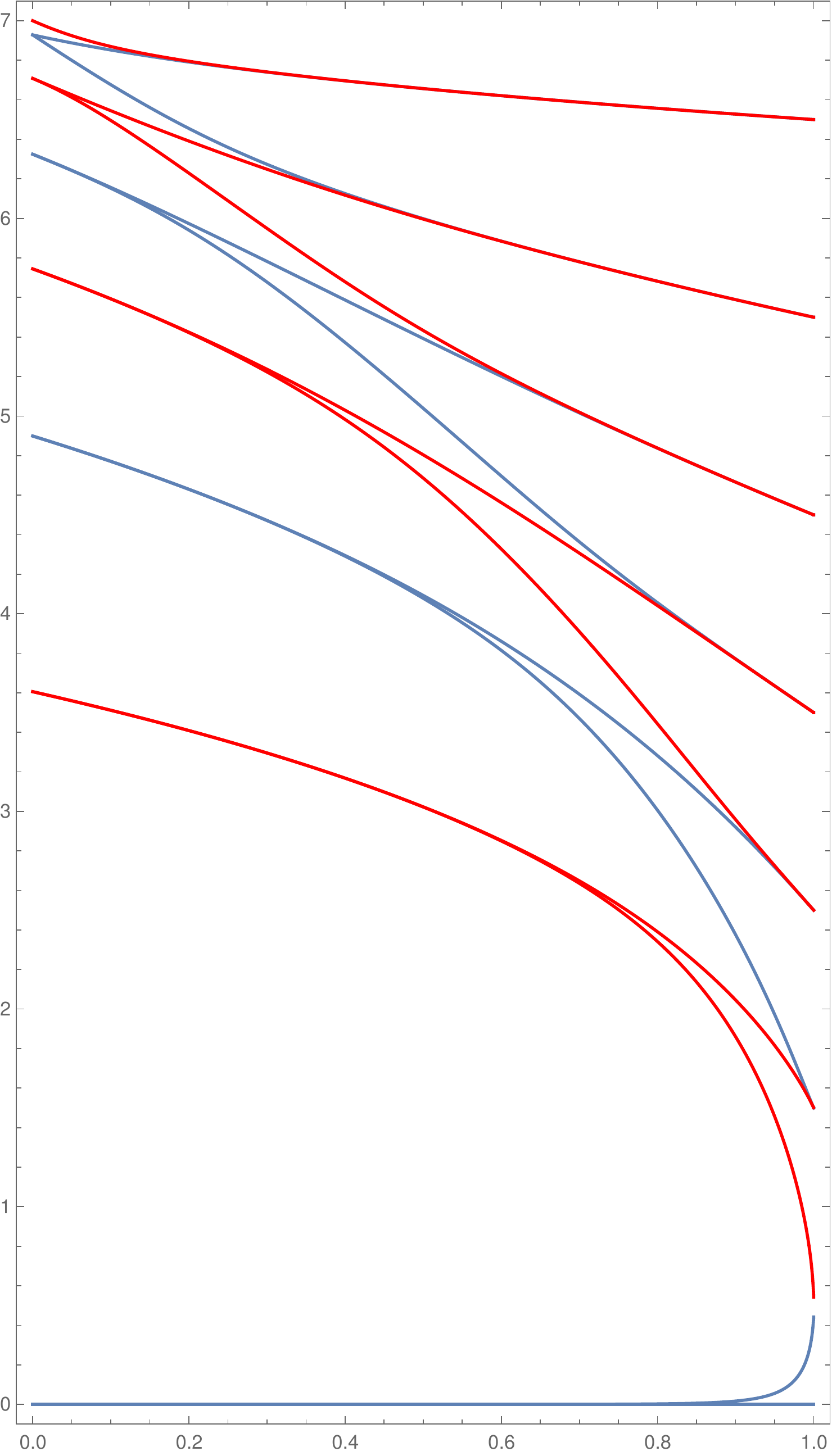}\\
\end{minipage}\kern-4em}
\caption{%
Periodic (red) and antiperiodic (blue) eigenvalues (vertical axis) of the spectrum as a function of the elliptic parameter~$m$ (horizontal axis) for a few integer values of $A$.
Bottom left: $A=3$. 
Top left: $A=4$.
Right: $A=7$.}
\label{f:spectrum_m}
\end{figure}

We further illustrate the results of this work by presenting some concrete plots of the spectrum.
We begin with the case of $A\in\Natural$.
Figure~\ref{f:spectrum_m} shows the periodic (red) and antiperiodic (blue) eigenvalues along the imaginary $z$-axis 
(vertical axis in the plot) 
as a function of the elliptic parameter~$m$ (horizontal axis) for a few integer values of $A$, namely:
$A=3$ (bottom left), 
$A=4$ (top left) and  
$A=7$ (right).
Note how all gaps are closed when $m=0$ and how they open immediately as soon as $m>0$ and remain open 
for all $m\in(0,1)$.
In the singular limit $m\to1^-$, the band widths tend to zero, and the periodic and antiperiodic eigenvalues 
``collide'' into the point spectrum of the operator $L$ on the line.

\begin{figure}[t!]
\newdimen\newfigwidth
\newfigwidth 55.5mm
\vglue-2.4ex
\centerline{\includegraphics[width=\newfigwidth]{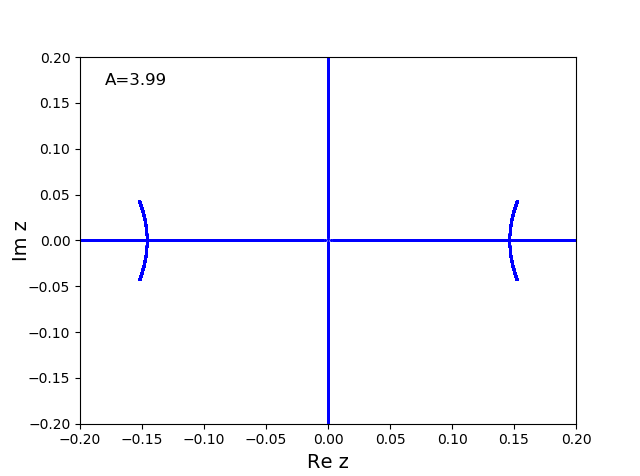}
\hspace{-6mm}
\includegraphics[width=\newfigwidth]{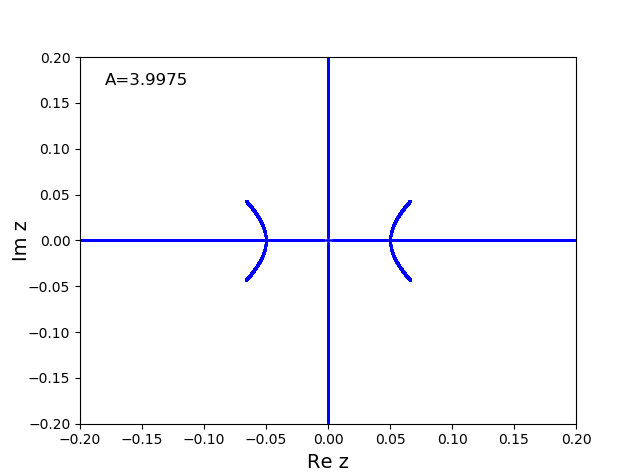}
\hspace{-6mm}
\includegraphics[width=\newfigwidth]{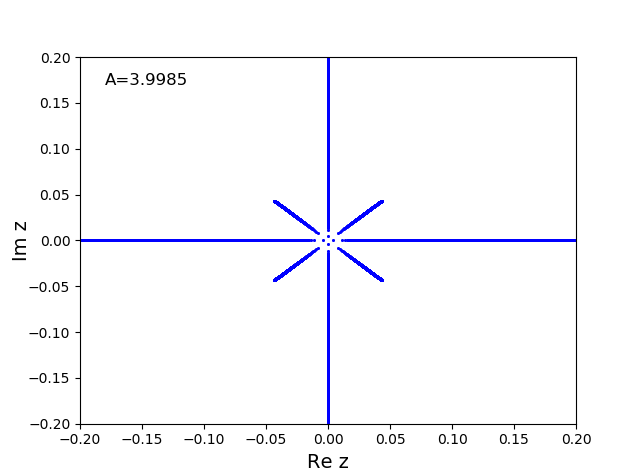}}
\vspace{-1mm}
\centerline{\includegraphics[width=\newfigwidth]{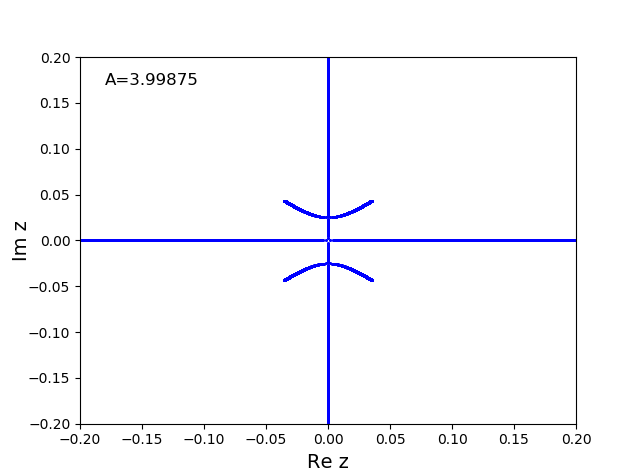}
\hspace{-6mm}
\includegraphics[width=\newfigwidth]{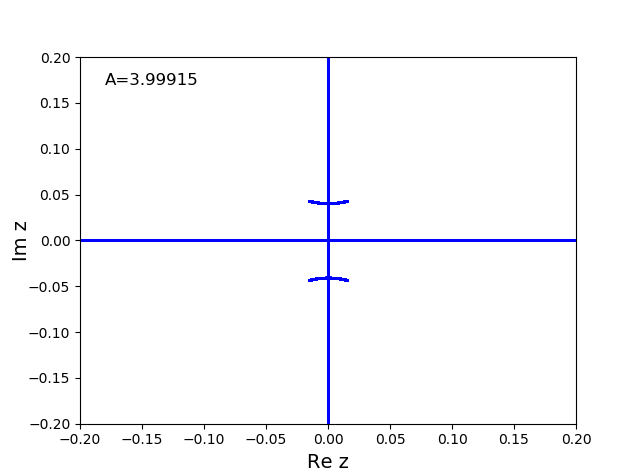}
\hspace{-6mm}
\includegraphics[width=\newfigwidth]{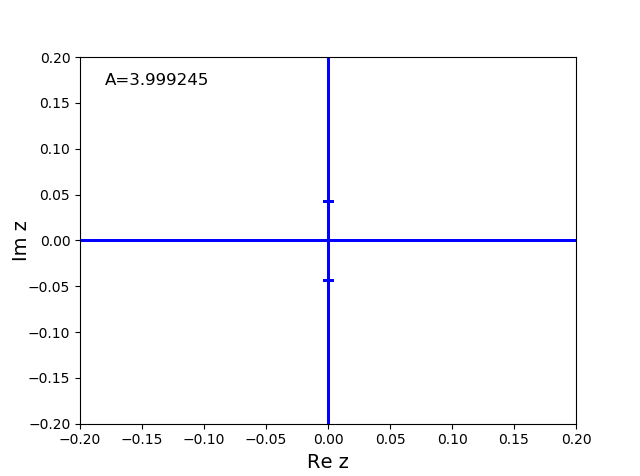}}
\vspace{-1mm}
\centerline{\includegraphics[width=\newfigwidth]{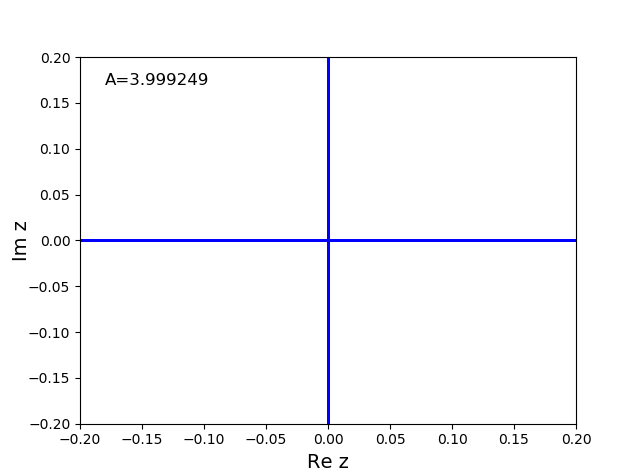}
\hspace{-6mm}
\includegraphics[width=\newfigwidth]{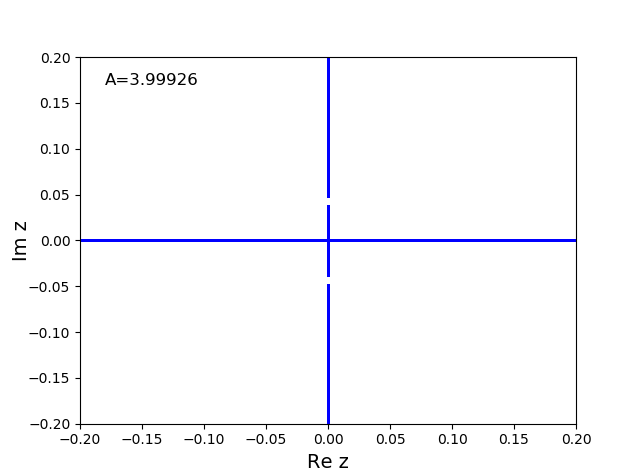}
\hspace{-6mm}
\includegraphics[width=\newfigwidth]{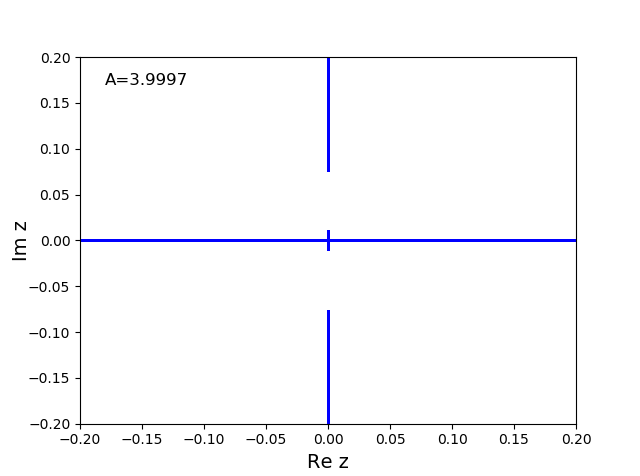}}
\vspace{-1mm}
\centerline{\includegraphics[width=\newfigwidth]{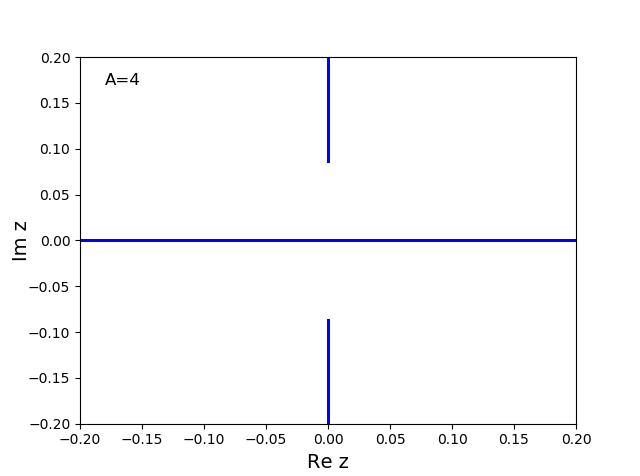}
\hspace{-6mm}
\includegraphics[width=\newfigwidth]{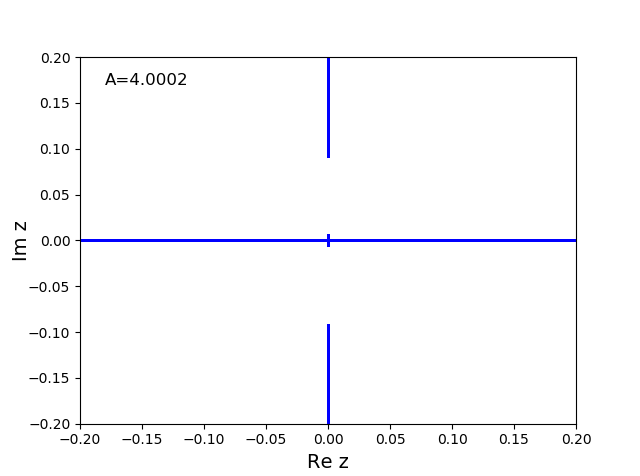}
\hspace{-6mm}
\includegraphics[width=\newfigwidth]{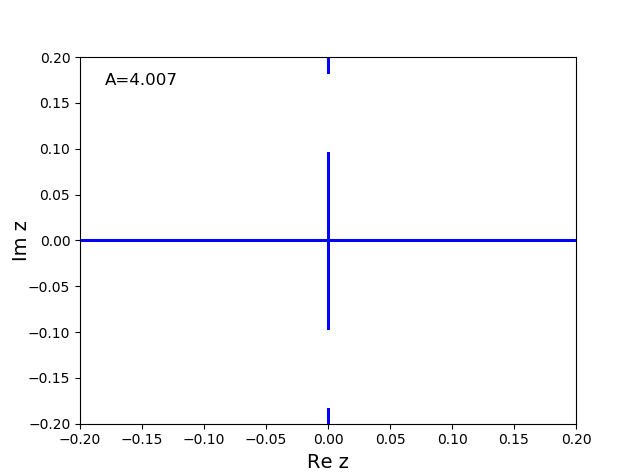}}
\vspace{-1mm}
\centerline{\includegraphics[width=\newfigwidth]{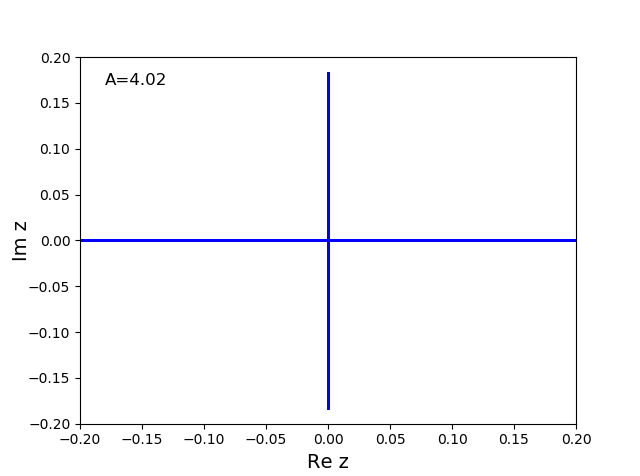}
\hspace{-6mm}
\includegraphics[width=\newfigwidth]{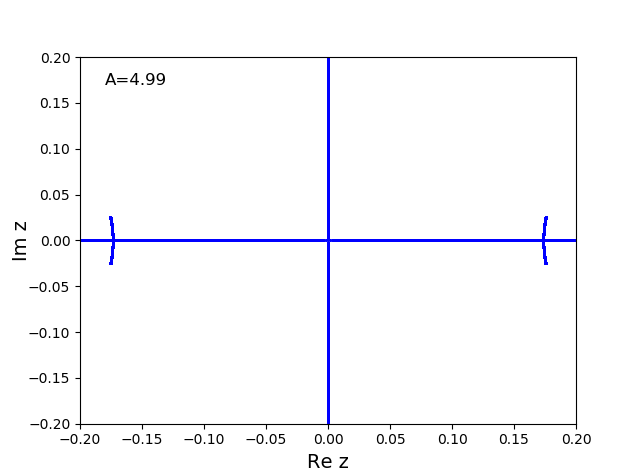}
\hspace{-6mm}
\includegraphics[width=\newfigwidth]{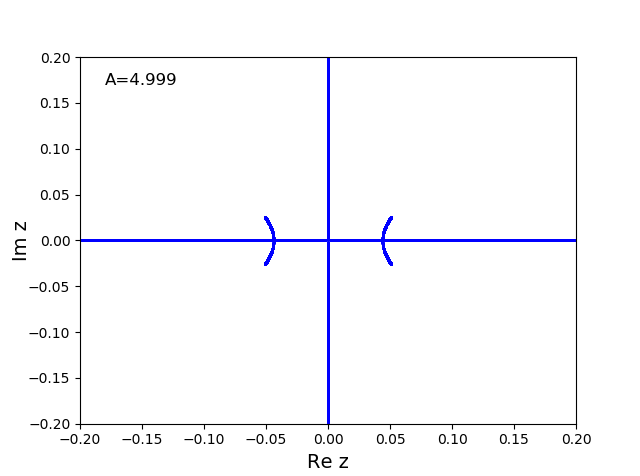}}
\kern\medskipamount
\caption{The Lax spectrum $\Sigma(L)$ [computed numerically using Hill's method (see~\cite{DK2006})] with potential $q\equiv A\dn(x;m)$, $m=0.9$, and increasing values of $A$, illustrating the formation of new bands and gaps for non-integer values of $A$.
}
\kern-\bigskipamount
\label{f:spectrum_A}
\end{figure}

Next, Figure~\ref{f:spectrum_A} shows the Lax spectrum (blue curves) in the complex $z$-plane for several
non-integer values of~$A$,
illustrating the formation of new bands and gaps as function of~$A$.
Note that the range of values for the real and imaginary parts of $z$ allows one to see only a small portion of 
the Lax spectrum.
For example, not visible outside the plot window are various bands and gaps along the imaginary axis
(cf.~Fig.~\ref{f:spectrum_m}) as well as the infinite number of spines growing off the real axis when $A\notin\N$.
However, selecting a larger portion of the complex $z$-plane would have made it more difficult 
to see the dynamics of the bands and gaps near the origin.
Starting from the smallest value of $A$ in the set ($A=3.99$, top left panel), one can see how,
as $A$ increases, a spine is pulled towards the origin, which it reaches at approximately $A=3.9985$ (top right panel).
As $A$ increases further, the spine moves along the imaginary axis, simultaneously shrinking to zero at
approximately $A=3.999249$ (left plot in the third row).
(Note that, even though the band is effectively gone at this value of $A$, the corresponding potential 
is still not finite-band due to the infinitely many spines that are still present outside the plot window.) 
As $A$ increases further, the band edges of the previous spine bifurcate along the imaginary axis, 
giving rise to a new gap.
Finally, at $A=4$ (left plot in the fourth row), the lower edge of this new gap reaches the origin.  
This is also exactly the value of $A$ at which the infinitely many spines shrink to zero.
As $A$ increases further, the band centered at $z=0$ reappears, 
a new spine gets sucked towards the origin, and the cycle repeats.

In summary,
every time $A$ increases by one unit, two more spines from the real axis gets pulled into the imaginary axis,
and a two new Scwarz symmetric gaps open on the imaginary axis.
When $A$ hits the next integer value, 
the lower band edges in the upper-half plane and the corresponding one in the lower-half plane reach the point $z=0$.
Simultaneously, all remaining spines emanating from the real $z$-axis shrink to zero,
giving rise to a finite-band potential.
As $A$ keeps increasing, the spines grow back, and the process repeats.

It is also interesting to briefly describe the dynamics of the zero-level curves of $\Im \D(z;m)$ near $z=0$.
Thanks to \eqref{Del''m=0} and~\eqref{Del''zm=0},
we know that $\D_{zz}(0;0)=0$ if and only if $A\in\Z$. 
One can also see that $\D_{zzzz}(0;0)\not=0$ when $A\in\Z$.
So, when $A\in\Z$, there are exactly eight zero-level curves of  $\Im \D$ emanating from $z=0$. 
As it follows from Corollary \ref{cor-z=0}, 
under a small $m>0$ deformation from $m=0$, 
a pair of these level curves will move up along the imaginary $z$-axis, 
while the other pair will symmetrically move down along $\i\R^-$.
It also follows from \eqref{D,sm=0} that $s(z;0)$ has a second-order zero at $z=0$;
that is, $s_1(0)\neq 0$,
which will remain in place under a small $m>0$ deformation according to \eqref{det-eq}, \eqref{Dsc} 
and Corollary \ref{cor-z=0}.  
This is another way to show that $s_1(m)\neq 0$ for small $m>0$.

\section{Discussion and concluding remarks}
\label{s:conclusions}
The results of this work provide an extension to 
the non-self-adjoint operator \eqref{e:Diracoperator}
of the classical works of Ince \cite{Ince,Ince2,Ince3}.
The results of this work also provide: 
(i) an example of Hill's equation with a complex, PT-symmetric potential 
(and a corresponding complex deformation of Ince's equation)
whose spectrum is purely real, which is especially relevant, 
since the study of quantum mechanics with non-Hermitian, PT-symmetric potential continues to attract considerable interest 
(e.g. see~\cite{Bagarello2015,Bender2005,NaturePhys2018} and references therein),
(ii) an example of an exactly solvable connection problem for Heun's ODE,
and
(iii) for the fist time a perturbation approach to study the determination of the genus, and the movable Dirichlet eigenvalues was presented.

We point out that the fact that the elliptic potential~\eqref{e:ellipticpotential} is finite-band for any $A\in\Natural$ 
can also be obtained as a consequence of the results of~\cite{GW1998},
where the potential $q(x) = n\,(\zeta(x) - \zeta(x-\omega_2) - \zeta(\omega_2))$ was studied
(where $\zeta(x)$ is Weierstrass' zeta function and $\omega_2$ one of the lattice generators~\cite{NIST})
and was shown to be finite-band when $n\in\Natural$ using the criteria introduced there
(see Appendix~\ref{a:gesztesy} for details).
In Appendix~\ref{a:gesztesy} we also discuss other elliptic potentials satisfying the criteria laid out in \cite{GW1998}. 
On the other hand, no discussion of the spectrum (i.e., location of the periodic/antiperiodic eigenvalues and 
of the spectral bands) was present in \cite{GW1998}. 

It is also the case that the elliptic potential~\eqref{e:ellipticpotential} 
is associated with the so-called Trebich-Verdier potentials \cite{TreibichVerdier} for Hill's equation 
(which are known to be algebro-geometric finite-band, see \cite{Takemura})
if and only if $A\in\Natural$,
as we show in Appendix~\ref{a:treibichverdier}.
To the best of our knowledge, this connection had not been previously made in the literature.

The family of elliptic potentials~\eqref{e:ellipticpotential} is especially important from an applicative point of view,
since (as was discussed in Section~\ref{s:introduction}) 
it interpolates between the plane  wave potential $q(x) \equiv A$ when $m=0$ and the Satsuma-Yajima (i.e., sech) potential 
$q(x)\equiv A\,\sech x$ when $m=1$, 
which, when $A\in\Natural$, gives rise to the celebrated $A$-soliton bound-state solution of the focusing NLS equation [\eqref{e:nls} with $s=1$].

The potential $q(x)\equiv A\,\sech x$ has also been used in relation to the semiclassical limit of the focusing NLS equation.  
This is because, by letting $A = 1/\epsilon$ and performing a simple rescaling $t\mapsto \epsilon t$ of the temporal variable, 
\eqref{e:nls} is mapped into the semiclassical focusing NLS equation 
\be
\i\epsilon q_t + \epsilon^2 q_{xx} + 2|q|^2q = 0\,,
\label{e:semiclassicalnls}
\ee
with the rescaled initial data $q(x,0)\equiv \dn(x;m)$.
The dynamics of solutions of \eqref{e:semiclassicalnls} has been studied extensively in the literature
(e.g., see \cite{bertolatovbis,bronskikutz1999,ClarkeMiller2002,Conti09,EKT2016,JM2013,KMM2003,lyng2012,TVZ2004}).
In particular, it is known that, for a rather broad class of single-lobe initial conditions
(including $q(x)\equiv \sech x$),
the dynamics gives rise to a focusing singularity (gradient catastrophe) that is regularized by the formation
of high-intensity peaks regularly arranged in the pattern of genus-2 solutions of the NLS equation.
The caustic (i.e., breaking) curve along which the genus-2 region breaks off from the genus-0 region 
(characterized by a slowly modulated plane wave, in which the solution does not exhibit short-scale oscillations)
has also been characterized,
and it is conjectured that additional breaking curves exist, giving rise to regions of higher genus.

All of the above-cited works studied localized potentials on the line.
However, similar behavior was observed for \eqref{e:semiclassicalnls} with periodic potentials in \cite{BO2020},
where a formal asymptotic characterization of the spectrum of the Zakharov-Shabat system~\eqref{e:ZS}
in the semiclassical scaling was obtained using WKB methods,
and in particular it was shown that the genus is $O(1/\epsilon)$ as $\epsilon\to0^+$. 
Some of the numerical results of \cite{BO2020} about the localization of the spectrum were rigorously proved in \cite{BOT2020}.
The results of the present work provide some rigorous evidence, for the dn potential~\eqref{e:ellipticpotential}, 
in support of the formal results of~\cite{BO2020} about the genus of the potential as a function of~$\epsilon$.

We emphasize that, even though we limited our attention to the focusing NLS equation for simplicity,
all the equations of the infinite NLS hierarchy 
(including the modified KdV equation, higher-order NLS equation, sine-Gordon equation, etc.) 
share the same Zakharov-Shabat scattering problem~\eqref{e:ZS}.
Therefore, the results of this work provide a two-parameter family of finite-band potentials for all the equations in the focusing NLS hierarchy. 

The results of this work open up a number of interesting avenues for further study.  
In particular, an obvious question is whether these potentials are stable under pertubations.
The stability of genus-1 solutions of the focusing NLS equation was recently studied in \cite{DeconinckSegal}
by taking advantage of the machinery associated with the Lax representation. 
A natural question is therefore whether similar results can also be used for the higher-genus potentials when $A>1$ or whether different methods are necessary.

Another interesting question is whether more general elliptic finite-band potentials related to~\eqref{e:ellipticpotential} exist.
Recall that, for the focusing NLS equation on the line, 
the potential $q(x) = A\,\sech x\,\e^{-\i a\log(\cosh x)}$ 
(which reduces to $q(x) = A\,\sech x$ when $a=0$) 
was shown in \cite{TovbisVenakides,TVZ2004} to be amenable to exact analytical treatment.
It is then natural to ask whether exactly solvable periodic potentials also exist related to~$q(x) = A\,\dn(x;m)$
but with an extra non-trivial periodic phase.

Yet another question is related to the time evolution of 
the potential~\eqref{e:ellipticpotential} according to the focusing NLS equation.
When $A=1$, time evolution is trivial, and the corresponding solution of the NLS equation is simply 
$q(x,t) = \dn(x;m)\,\e^{\i(2-m) t}$.
That is not the case when $A>1$, however.
For the Dirac operator~\eqref{e:Diracoperator} on the line with reflectionless potentials,
sufficient conditions were obtained in \cite{SLEPJP2018} guaranteeing that,
if the discrete spectrum is purely imaginary, 
the corresponding solution of the focusing NLS equation is periodic in time.
The natural question is then whether a similar result is also true for the elliptic potential~\eqref{e:ellipticpotential},
namely, whether such potentials generates a time-periodic solution of the focusing NLS equation when $A\in\Natural$.

The semiclassical limit of certain clases of periodic potentials (including the potential $\dn(x;m)$)
generates a so-called breather gas for the focusing NLS equation (e.g., see \cite{BO2020,TW2022}), 
which is to be understood as the thermodynamic limit of a finite-band solution of the focusing NLS equation where the genus $G\ra\infty$ and simultaneously all bands but one shrink in size exponentially fast in $G$ (see \cite{ET2020} for details).
It was proposed that such gases be called periodic breather gases.
Periodic gases have the important feature that, together with their spectral data (i.e., independent of the phase variables) such as the density of states, one also can obtain some information on a ``realization" of the gas, namely, on the  semiclassical evolution of the given periodic potential.   
Thus, progress in studying the $(2A-1)$-band solutions of the focusing NLS equation (with $A\in\N$) generated by 
the potential $q(x)=A\dn(x;m)$, and especially its large $A$ limit, is of definite interest.
In fact, finite-band solutions to integrable systems (such as the KdV and NLS equations) generated by elliptic potentials, 
based on the work of Krichever \cite{Krichever1980},
were studied in the literature. 
We will not go into the details of those results here, 
but it should be clear that any details about the family of finite-band solutions of the focusing NLS equation that homotopically ``connect" the known behavior of the plane wave and the multi-soliton solutions 
will be very interesting to obtain and analyze.

\subsection*{Acknowledgments}

This work was partially supported by the U.S.\ National Science Foundation under grant numbers\break DMS-2009487 (GB), and DMS-2009647 (AT),
and the National Science Foundation of China under grant number 12101590 (XDL).
Part of this work was supported by the U.S.\ National Science Foundation under grant number 1440140 while JO was in residence at the Mathematical Sciences Research Institute (MSRI) in Berkeley, California, during the Fall 2021 semester program 
``Universality and Integrability in Random Matrix Theory and Interacting Particle Systems''.

\setcounter{section}1
\setcounter{subsection}0
\addcontentsline{toc}{section}{Appendix}
\def\thesubsection{\Alph{section}.\arabic{subsection}}
\def\theequation{\Alph{section}.\arabic{equation}}
\def\thetheorem{\Alph{section}.\arabic{theorem}}
\def\thefigure{\Alph{section}.\arabic{figure}}

\section*{Appendix}
\label{a:appendix}

\subsection{Notation and function spaces}
\label{a:notation}
The Pauli spin matrices, used throughout this work, are defined as
\vspace*{-1ex}
\be
\sigma_1 := \begin{pmatrix*}[r] 0 & 1 \\ 1 & 0 \end{pmatrix*}\,,
\qquad
\sigma_2 := \begin{pmatrix*}[r] 0 & -\i \\ \i & 0 \end{pmatrix*}\,,
\qquad
\sigma_3 := \begin{pmatrix*}[r] 1 & 0 \\ 0 & -1 \end{pmatrix*}\,.
\ee
Moreover,
$L^\infty(\R;\C^2)$ is the space of essentially bounded Lebesgue measurable two-component vector functions with the essential supremum norm.
Given the interval $I_{x_o} = [x_o,x_o+l]$ we define the inner product between two-component Lebesgue measurable functions $\phi$ and $\psi$ as
\vspace*{-1ex}
\be
\langle \phi,\psi \rangle := \int_{x_o}^{x_o+l} \big(\phi_{1}\overline{\psi}_1 + \phi_2\overline{\psi}_2\big) \,\d x\,.
\label{e:innerproduct}
\ee
Then $L^2(I_{x_o};\Complex^2)$ denotes the set of two-component Lebesgue measurable vector functions that are square integrable, i.e.,
$\|\phi\|_{L^2(I_{x_o};\C^2)} := \langle \phi,\phi \rangle^{1/2} < \infty$.
Similarly, we define the inner product
of two scalar Lebesgue measurable functions $f$ and $g$ as
\be
\langle f,g \rangle := \int_{x_o}^{x_o+l} f\hspace*{0.1ex}\overline{g} \,\d x\,.
\label{e:innerproductscalar}
\ee
Then $L^2(I_{x_o};\C)$ denotes the set of scalar Lebesgue measurable functions that are square integrable, i.e.,
$\|f\|_{L^2(I_{x_o};\C)} := \langle f,f \rangle^{1/2} < \infty\,$.
Finally,
we define the inner product 
of two infinite sequences $c=\{c_n\}_{n\in\Z}$ and $d=\{d_n\}_{n\in\Z}$ as
\be
\label{e:innerproductsum}
\langle c,d \rangle := \sum_{n\in\Z} c_n\overline{d}_n\,.
\ee
Then $\ell^2(\Z)$ denotes the set of square-summable sequences, i.e.,
$\|c\|_{\ell^2(\Z)} := \langle c,c \rangle^{1/2} < \infty\,$.
Finally,
the space of continuous functions on the real axis is denoted $C(\R)$,
and $\N_o:=\N\cup\{0\}$. 
%

\subsection{Proof of two lemmas}
\label{a:newboundproof}

\textit{Proof of Lemma~\ref{l:combinedlemma}.}
To prove part (i) we begin by writing \eqref{e:Diraceigenvalueproblem} as the coupled system of linear differential equations~\eqref{e:ZS}.
By Floquet theory $z \in \Sigma(L)$ if and only if $\phi = \e^{\i \nu x}\psi$,
where $\psi=(\psi_1,\psi_2)^{T}$ with $\psi(x+l;z) = \psi(x;z)$,
and $\nu \in [0,2\pi/l)$.
Plugging this expression into~\eqref{e:ZS} yields the modified system 
\vspace*{-1ex}
\be
\label{e:zsmod}
\i \psi_{1,x} - \i q\psi_{2} = (z + \nu)\psi_{1}\,,
\qquad
\i \psi_{2,x} + \i q\psi_{1} = (-z + \nu)\psi_{2}\,,
\ee
Multiplying the first of these equations by $\overline{\psi}_{1}$ and taking the complex conjugate yields two equations, which we integrate over a full period,
arriving at the expressions 
\be
\i \langle q\psi_2, \psi_1 \rangle = -\i\langle \psi_1, \psi_{1,x} \rangle - (z+\nu) \| \psi_1\|_{L^2([0,l])}^2 \,,
\qquad
\i \langle \psi_1, q\psi_2 \rangle = \i\langle \psi_1, \psi_{1,x} \rangle + (\overline{z}+\nu) \|\psi_1\|_{L^2([0,l])}^2 \,,
\nonumber
\ee
where boundary terms vanish since $\psi_{1}(x+l;z)=\psi_{1}(x;z)$. 
Adding these two expressions then one gets 
\be
-\Im z \|\psi_{1}\|_{L^2([0,l])}^{2} = 
\Re \langle q \psi_{2}, \psi_{1} \rangle\,.
\label{e:eq1}
\ee
Similarly, the second equation of \eqref{e:zsmod} 
yields 
\be
\i\langle \psi_1, q\psi_2 \rangle = \i\langle \psi_2, \psi_{2,x} \rangle + (-z+\nu) \|\psi_2\|_{L^2([0,l])}^2\,,
\qquad
\i\langle q\psi_2, \psi_1\rangle = -\i\langle \psi_2, \psi_{2,x} \rangle + (\overline{z}-\nu) \|\psi_2\|_{L^2([0,l])}^2 \,,
\nonumber
\ee
as well as
\be
-\Im z \|\psi_{2}\|_{L^2([0,l])}^{2} = 
\Re \langle q\psi_2, \psi_1 \rangle\,.
\label{e:eq2}
\ee
Equating \eqref{e:eq1} and \eqref{e:eq2} we conclude that if $|\Im z|>0$, 
then
\be
\label{e:intermediate}
\|\psi_{1}\|_{L^2([0,l])} = \|\psi_{2}\|_{L^2([0,l])}\,.
\ee
Next, note that
$|\langle q\psi_2, \psi_1 \rangle| \leq \langle |q\psi_2|, |\psi_1| \rangle$.
Also, since $q$ is not constant, there exists $(a,b)\subset(0,l)$ such that
$|q(x)| < \|q\|_{\infty}$ for $x\in(a,b)$.
Thus,
for $|\Im z| > 0$ it follows from~\eqref{e:intermediate} and the H\"{o}lder inequality that
\begin{align*}
0 < |\Im z|\|\psi_1\|_{L^{2}([0,l])}^{2} = |\Re \langle q\psi_2, \psi_1 \rangle| 
&< \|q\|_{\infty}\langle |\psi_2|, |\psi_1| \rangle \\
&\leq \|q\|_{\infty}\|\psi_1\|^2_{L^2([0,l])}\,.
\end{align*} 
Hence $|\Im z| < \|q\|_{\infty}$ for $z\in\Sigma(L)$.
The proof of (ii) which can be found in \cite{BOT2020} follows 
from Lemma~\ref{l:quartets} and since $\D(z)$ is an analytic function of $z$.  
\qed

\begin{lemma}
\label{l:M_analytic_m}
Consider the Dirac operator~\eqref{e:Diracoperator}.
If the potential $q\equiv A\dn(x;m)$,
then the monodromy matrix $M(z;m)$ is an analytic function of $m$ for any $m\in[0,1)$. 
\end{lemma}

\textit{Proof.}
The result follows from two key facts. 
One is that $\dn(x;m)$ is an analytic function of $m$ for all $|m|\le 1$ \cite{Walker}.
The second is the fact that solutions of ODEs with analytic dependence on variables and parameters are analytic (see~\cite{CoddingtonLevinson} pp.~23--32 and~\cite{Ince1956} p.~72). 
\qed

\subsection{Solution of the ZS system at $z=0$}
\label{a:exactsolution}

We have seen that, 
for $q$ real, 
the eigenvalue problem~\eqref{e:Diraceigenvalueproblem} can be reduced to second-order scalar ODEs~\eqref{e:tindschrod5}.
Consider~\eqref{e:tindschrod5} with $\lambda=0$, namely 
$v_{xx}^{\pm} + (\pm \i q_x + q^2)v^{\pm} = 0$.
Using the ansatz $v^{\pm} = \e^{\pm f}$, 
one gets
$\pm f_{xx} + (f_{x})^2 \pm \i q_{x} + q^2 = 0$.
Then, letting $g = f_{x}$ 
yields
$\pm g_{x} + g^2 = \mp \i q_x - q^2$
with a solution given by $g = \mp \i q$.
In particular,
it follows $g^2 = -q^2$.
Hence,
we have derived the following solution to the ODEs~\eqref{e:tindschrod5} for $\lambda=0$,
namely,
\be
v^{\pm}(x;0) = \e^{\mp \i \int^{x}_{0}q(s)\,\d s}\,.
\label{e:y1}
\ee
Next, 
using the invertible change of variables~\eqref{e:changevar_5}, 
one gets the following solution to the eigenvalue problem~\eqref{e:Diraceigenvalueproblem} when $z = 0$:
\bse
\be
\phi(x;0) = \bigg( \cos \Big( \int_{0}^{x}q(s)\,\d s\Big), -\sin \Big( \int_{0}^{x}q(s)\,\d s\Big) \bigg)^\T\,.
\label{e:zssoln1_z0}
\ee
Moreover,
using Rofe-Beketov's formula~\cite{rofebek84}, 
one obtains a second linearly independent solution as:
\be
\widetilde{\phi}(x;0) = \bigg( \sin \Big( \int_{0}^{x}q(s)\,\d s\Big), \cos \Big( \int_{0}^{x}q(s)\,\d s\Big) \bigg)^\T\,.
\label{e:zssoln2_z0}
\ee
\ese
Thus, the Floquet discriminant~\eqref{e:DeltaM} for eigenvalue problem~\eqref{e:Diraceigenvalueproblem} at $z=0$ is
\be
\D(0) = \cos \Big( \int_{0}^{l}\, q(s)\, \d s \Big)\,,
\label{e:Dsoln_z0}
\ee
where $l$ is the period of the potential.
We can now prove Lemma~\ref{l:zequalzero}. 

\textit{Proof of Lemma~\ref{l:zequalzero}.}
Using well-known properties of the Jacobi elliptic functions (see \cite{Gradshteyn, NIST}), 
when $q(x)=A\,\dn(x;m)$, 
\eqref{e:zssoln1_z0}, \eqref{e:zssoln2_z0} and~\eqref{e:Dsoln_z0} yield, respectively
\vspace*{-1ex}
\bse
\begin{gather}
\phi(x;0,A,m) = \Big( \cos( A\,\am(x;m)), -\sin( A\,\am(x;m)) \Big)^\T\,,
\label{e:v_amx}
\\
\widetilde{\phi}(x;0,A,m) = \Big(\sin(A\,\am(x;m)), \cos(A\,\am(x;m))\Big)^\T\,,
\label{e:v_amx2}
\\
\D(0;A,m) = \cos(A\pi)\,.
\label{e:cosnpi}
\end{gather}
\ese
In particular,
$\phi(0;0,A,m)=(1,0)^\T$,
and $\phi(2K;0,A,m)=(\cos(A\pi),\sin(A\pi))^\T$.
Therefore, we have that
$\phi(x+2K;0,A,m) = \phi(x;0,A,m)$ if and only if $A \in 2\Z$,
and 
$\phi(x+2K;0,A,m) = -\phi(x;0,A,m)$ if and only if $A \in 2\Z+1$,
with similar relations for $\widetilde\phi$.
Thus,
when $A$ is an even integer, $z=0\in\Sigma_+(L)$,
whereas when $A$ is an odd integer, $z=0\in\Sigma_-(L)$.
Finally, the above calculations also show that,  
for $q\equiv A\dn(x;m)$ with $A\in \Z$ the eigenvalue $z=0$ has geometric multiplicity two. 
\qed

\subsection{Transformation of the ZS system into a Heun system}
\label{a:heun}

If the potential $q$ is real, 
then the change of dependent variable
$\phi = \half\diag(1,-\i)(\s_3+\s_1)\,v$,
maps the ZS system \eqref{e:ZS} into the equivalent system
\be
v_x = -\i(z\s_1 + q\s_3)v\,.
\label{e:zs2}
\ee
Then the transformation $t=2\am(x;m)$ maps~\eqref{e:zs2} to the trigonometric first-order system
\be
\label{e:zs3}
v_t=-\halfi \Bigg(A\s_3+\frac{z\s_1}{\sqrt{1-m\sin^2\frac t2}}\Bigg)\,v\,, 
\ee
which is equivalent to~\eqref{e:trigonometricODE}. 
Finally, 
the transformation $\z=\e^{\i t}$ maps~\eqref{e:zs3} to 
\be
\label{e:zs4a}
\z v_\z=-\half \Bigg(A\s_3+\frac{z\s_1}{\sqrt{1-\frac m2\big(1-\frac 12(\z+\z^{-1})\big)}}\Bigg)v \,,
\ee
and the transformation
\be
v = \Xi w\,,\quad
\Xi = \diag\left(1, \txtfrac{1}{z}{\sqrt{1-\txtfrac{m}{2}\big(1-\half (\z+\z^{-1})\big)}}\right)\,,
\ee
maps~\eqref{e:zs4a} to the Heun system~\eqref{zs4} where $\l=z^2$. 
The Heun system \eqref{zs4} has four (regular) singular points, located at $\zeta = 0,\z_{1,2},\infty$, 
where $\z_{1,2}$ are zeros of the denominator in~\eqref{zs4}. 
The Frobenius exponents at the singularities can be derived directly from~\eqref{zs4}. 

\subsection{Augmented convergence and Perron's rule}
\label{s:perron}

In general the Frobenius series~\eqref{e:frobenius1} with base point $\z=0$ only converges for $|\z| < |\z_1|$ and 
the series~\eqref{e:frobenius2} with base point $\z=\infty$ only converges for $|\z| > |\z_2|$.
Therefore, neither expansion is convergent on $|\zeta|=1$ in general.
However, there exist certain values of $\lambda$ for which one or both of the Frobenius series have a larger 
(i.e., augmented) radius of convergence.
These are precisely the periodic/antiperiodic eigenvalues of the problem,
and Perron's rule provides a constructive way to identify them
(see also~\cite{Arscott,erdelyi,Ince2,Ince1956,Ronveaux} for further details).

We begin by noting that, by dividing all coefficients by $n^2$,
all four three-term recurrence relations~\eqref{e:threetermrecurrence@z=0}, \eqref{e:heunrecur2}, \eqref{e:heunrecur3} and \eqref{e:heunrecur4} can be rewritten as 
\vspace*{-1ex}
\bse
\label{e:perronrecurrence}
\begin{gather}
e_0c_0 + f_0c_1 = 0\,,\qquad n=0\,,
\label{e:recurex1}
\\
d_nc_{n-1}+e_nc_n+f_nc_{n+1} = 0\,,\qquad n=1,2,\dots
\label{e:recurex2}
\end{gather}
\ese 
with $f_n\ne 0$, 
and 
$d_n \to d$, 
$e_n\to e$, 
and $f_n \to f$ as $n\to \infty$.
Perron's rule~\cite{erdelyi, Perron} states that, 
if $\xi_\pm$ are the roots of the quadratic equation
\be
f\xi^2 + e\xi + d = 0\,,
\label{e:perronquadratic}
\ee
with $|\xi_-|<|\xi_+|$,
then 
$\lim_{n\to\infty} c_{n+1}/c_{n}=\xi_+$, unless 
the coefficients $d_n,e_n,f_n$ satisfy the infinite continued fraction equation 
\vspace*{0ex}
\be
\frac{e_o}{f_o} = \cfrac{d_1}{e_1 - \cfrac{d_2f_1}{e_2 - \cfrac{d_3f_2}{e_3 - \cdots}}}\,,
\label{e:continuedfractionequation}
\ee
in which case $\lim_{n\to\infty} c_{n+1}/c_{n} = \xi_-$.
That is, Perron's rule implies that, generically,
the radius of convergence of the series $\Sigma_{n=0}^{\infty}c_{n}\z^{n}$ is $1/|\xi_+|$.
However, if and only if \eqref{e:continuedfractionequation} holds,
the radius of convergence is $1/|\xi_-|$, and therefore larger.
In our case, the roots $\xi_\pm$ of~\eqref{e:perronquadratic} are exactly the singular points $\z_{1,2}$ of 
Heun's ODE~\eqref{e:HeunODE}.
Then, since $e_n$ depends on $\lambda$, 
\eqref{e:continuedfractionequation} is a condition that determines the exceptional values of $\lambda$
that guarantee augmented convergence.
Indeed, 
\eqref{e:continuedfractionequation} is equivalent to requiring that~$\lambda$ is an eigenvalue of $T_{o}^\pm$ (resp. $T_{\infty}^\pm$).

\subsection{Generalized convergence of closed operators}
\label{a:genconv}
Here we briefly discuss the generalized convergence of closed operators,
(see \cite{kato} p. 197 for a detailed discussion).
Consider $\mathfrak{C}(\mathcal{X},\mathcal{Y})$ the space of closed linear operators
between Banach spaces.
If $T$, $S \in \mathfrak{C}(\mathcal{X},\mathcal{Y})$,
their graphs $G(T)$, $G(S)$ are closed linear manifolds on the product space $\mathcal{X}\times\mathcal{Y}$.
Set
$\hat{\delta}(T,S) = \hat{\delta}(G(T),G(S))$,
i.e.,
the \textit{gap} between $T$ and $S$. 
(See \cite{kato} p. 197 for the definition of $\hat{\delta}(T,S)$.)
Similarly,
we can define the \textit{distance} $\hat{d}(T,S)$ between $T$ and $S$ as equal to $\hat{d}(G(T),G(S))$.
(See \cite{kato} p. 198 for the definition of $\hat{d}(T,S)$.)
With this distance function $\mathfrak{C}(\mathcal{X},\mathcal{Y})$ becomes a metric space.

In this space the convergence of a sequence $T_n\to T\in \mathfrak{C}(\mathcal{X},\mathcal{Y})$ is defined by $\hat(T_n,T)\to 0$.
Since $\hat{\delta}(T,S)\leq\hat{d}(T,S)\leq 2\hat{\delta}(T,S)$ (\cite{kato} p. 198)
this is true if and only if $\hat{\delta}(T_n,T)\to 0$.
In this case we say $T_n\to T$ \textit{in the generalized sense}.
This notion of generalized convergence for closed operators is a generalization of convergence \textit{in norm}. 
Importantly,
the convergence of closed operators in the generalized sense implies the continuity of a finite system of eigenvalues (\cite{kato} p. 213).  

\subsection{Gesztesy-Weikard criterion for finite-band potentials}
\label{a:gesztesy}

According to Theorem~1.2 from~\cite{GW1998}, an elliptic potential $Q(x)$ of the Dirac operator~\eqref{e:Diracoperator} is finite-band if and only if its 
fundamental matrix solution $\Phi(x;z)$ is meromorphic in $x$ for all $z\in\C$.

\begin{theorem}
\label{t:geszt}
Consider \eqref{e:Diracoperator}, then $q\equiv A\dn(x;m)$ with $m\in(0,1)$ is finite-band if and only if $A\in\Z$.
\end{theorem}

\textit{Proof.} 
The (simple) poles of $\dn(x;m)$ within the fundamental period $2jK + 4n\i K'$ are at $x=\i K'$ and $x=3\i K'$ where $K':=K(1-m)$.
By the Schwarz symmetry, 
it is sufficient to consider only the pole at $\i K'$.
The residue at $\i K'$ is $-\i$ (\cite{Gradshteyn}, 8.151) and the local Laurent expansion is odd.
Let $\Phi(u):=\Phi(x-\i K';z)$ and note $\dn(x+\i K';m)=-\i\cn(u;m)/\sn(u;m)$ (\cite{BirdFriedman}, p. 20).
Substitution into \eqref{e:zs2} gives 
\be
\label{zs1-pole-1}
u\Phi_{u}(u)=[-\i zu\s_3 + (A+u^2p(u))\s_2]\Phi(u)=:B(u)\Phi(u),
\ee
where $p(u)$ is analytic near $u=0$ and even,  
and is meromorphic near $u=0$ for all $z\in\C$. 
The leading order term of $B(u) $ is $A\s_2$ with eigenvalues
$\l=\pm A$. Thus, meromorphic $\Phi(u)$ requires $A\in\Z$. 

To show that $A\in\Z$ is also a sufficient condition we need to show that $\Phi(u)$
does not contain logarithms, i.e., regular singular point $u=0$ is non-resonant. 
To do so we need to shift the spectrum of the leading term of $B(u)$ to a single point, for example, $-A$. Without loss of generality, we can assume $A>0$. Since
\be
\half(\1-\i\s_1)\s_2(\1+\i\s_1)=\s_3, \qquad \half(\1-\i\s_1)\s_3(\1+\i\s_1)=-\s_2,
\ee
we first diagonalize the leading term $A\s_2$ by the transformation $\Phi=(\1+\i\s_1)\tilde \Phi$. Then \eqref{zs1-pole-1} becomes
\be
\label{zs1-pole-2}
u\tilde\Phi_{u}=[\i zu\s_2+(A+u^2p(u))\s_3]\tilde\Phi.
\ee
Then the shearing transformation $\tilde\Phi={\rm diag}(u,1)\Psi$ reduces \eqref{zs1-pole-2} to
\be
\label{zs1-pole-3}
u\Psi_{u}=\left[
\begin{pmatrix*}
A-1  & z  \\
0  & -A 
\end{pmatrix*}+
u^2\begin{pmatrix*}
p(u)  & 0  \\
-z  & -p(u)
\end{pmatrix*}
\right]\Psi.
\ee
After diagonalizing the leading term, we obtain 
\be
\label{zs1-pole-4}
u\tilde\Psi_{u}=\left[
\begin{pmatrix*}
A-1  & 0  \\
0  & -A 
\end{pmatrix*}+
u^2\begin{pmatrix*}
\tilde p(u)  & r(u)  \\
-z  & -\tilde p(u)
\end{pmatrix*}
\right]\tilde\Psi,
\ee
where $\tilde p(u), r(u)$ are even and analytic at $u=0$ functions. Thus, the coefficient of \eqref{zs1-pole-4} is an analytic and even  matrix function. 

If $A=1$, one more shearing transformation would produce leading order term $-\1$, which is non-resonant
(no non trivial Jordan block) and, thus, the result would follow.
If $A>1$, we apply shearing transformations with the matrix  ${\rm diag}(u^2,1)$ with consequent diagonalizations that will shift the $(1,1)$ entry of the leading term by $-2$ and preserve the analyticity and evenness of the coefficient. When the difference of the eigenvalues of the leading term becomes one, we repeat the last step of the case $A=1$.
\qed

\begin{corollary}
For the Dirac operator~\eqref{e:Diracoperator},
$q\equiv A\,\cn(x;m)$ with $m\in(0,1)$ and $A>0$ is finite-band if and only if $A=\sqrt{m}n$ with $n\in\Z$,
while $q = A\,\sn(x;m)$ is not finite-band for any~$A>0$.
\end{corollary}

\textit{Proof.} 
The function $\cn(x;m)$ has the same locations of simple poles as $\dn(x;m)$. 
Given that the residues of the poles $2jK+\i K'$ for $j\in\Z$ of $\cn(x;m)$ are $(-1)^{j-1}\i/\sqrt m$, 
it is clear that the choice of $A$ given above leads to integer Frobenius exponents. 
To prove the non-resonance conditions, 
we notice that in Theorem \ref{t:geszt} we used only the fact that $\dn(x;m)$ has an odd Laurent expansion at the pole. Since this is also true for (\cite{Gradshteyn}, 8.151) the proof is complete.
Similar arguments show that $A\,\sn(x;m)$ is never finite-band.
\qed

\subsection{Connection between Heun's equation and Treibich-Verdier potentials}
\label{a:treibichverdier}

It was shown in \cite{Takemura} that the Heun ODE in standard form:
\begin{equation}
\label{e:generic Heun}
\frac{\d^{2}y}{\d \z^{2}} + \left(\frac{\gamma}{\z}+\frac{\delta}{\z-1}+\frac{\epsilon}{\z-a}\right)\frac{\d y}{\d\z}+\frac{\alpha \beta \z - \xi}{\z(\z-1)(\z-a)}y=0
\end{equation}
is associated with the so-called Treibich-Verdier potentials (defined below) for Hill's equation \cite{TreibichVerdier}, 
where $\alpha, \beta, \gamma, \delta, \epsilon, \xi, a$ (with each of them $\ne 0, 1$) are complex parameters linked by the relation $\gamma+\delta+\epsilon=\alpha+\beta+1$. 
Specifically, the Heun equation~\ref{e:generic Heun} can be transformed into 
\begin{equation}
\label{Heun-ell}
\left(-\frac{\d^{2}}{\d x^{2}}+\sum_{i=0}^{3}l_{i}(l_{i}+1)\wp(x+\omega_{i})-E\right)f(x)=0
\end{equation}
via the transformation $f(x)=y\z^{-l_{1}/2}(\z-1)^{-l_{2}/2}(\z-a)^{-l_{3}/2}$, where 
$\wp(x)$ is the Weierstrass $\wp$-function
with periods $\{2\omega_{1}, 2\omega_{3}\}$,  
where $\omega_{1}/\omega_{3} \not \in \R$ 
and
\begin{gather*}
\omega_{0}=0\,,\quad
\omega_{2}=-\omega_{1}-\omega_{3}\,,\quad
e_{i}=\wp(\omega_{i})\,,\quad
z=\frac{\wp(x)-e_{1}}{e_{2}-e_{1}}\,,\quad
a=\frac{e_{3}-e_{1}}{e_{2}-e_{1}}\,,
\end{gather*}
$E=(e_{2}-e_{1})[-4\xi+(-(\alpha-\beta)^{2}+2\gamma^{2}+6\gamma \epsilon+2\epsilon^{2}-4\gamma-4\epsilon-\delta^{2}+2\delta+1)/3
+ (-(\alpha-\beta)^{2}+2\gamma^{2}+6\gamma\delta+2\delta^{2}-4\gamma-4\delta-\epsilon^{2}+2\epsilon+1)a/3]$,
and the coefficients $l_{i}$ in \eqref{Heun-ell} are connected with the parameters in~\eqref{e:generic Heun} as follows:
\be
l_{0}=\beta-\alpha-\txtfrac{1}{2}\,,\quad 
l_{1}=-\gamma+\txtfrac{1}{2}\,,\quad 
l_{2}=-\delta+\txtfrac{1}{2}\,,\quad 
l_{3}=-\epsilon+\txtfrac{1}{2}\,.
\ee
It is known that the potential $\sum_{i=0}^{3}l_{i}(l_{i}+1)\wp(x+\omega_{i})$ is a (finite-band) Treibich-Verdier potential if and only if $l_{i}\in \Z$, $i=0, 1, 2, 3$~\cite{TreibichVerdier}.  
Note that the periods $\{2\omega_{1}, 2\omega_{3}\}$ of $\wp(x)$ are not uniquely determined.

In this subsection we show the special case of the Heun equation~(\ref{e:HeunODE}) considered in this work 
corresponds to a Treibich-Verdier potential if and only if $A\in \Z$.
To show this, we employ the conformal mapping $\tilde{\z}:={\z}/{\z_{1}}$ 
Under this transformation,
and recalling the relation $\z_{2}={1}/{\z_{1}}$,
the Heun equation~(\ref{e:HeunODE}) is mapped into
\begin{equation}
\label{e:Heun'}
\frac{\d^{2}y}{\d\tilde{\z}^{2}}+\frac{\frac{3}{2}\tilde{\z}^{2}-\big(\frac{2m-4}{m\z_{1}}\big)\tilde{\z}+\frac{1}{2\z_{1}^{2}}}{\tilde{\z}(\tilde{\z}-1)(\tilde{\z}-1/\z_{1}^{2})} \frac{\d y}{\d\tilde{\z}}-\frac{\frac{1}{4}A(A+1)\tilde{\z}^{2}+\big(\frac{\lambda+A^{2}(1-m/2)}{m\z_{1}}\big)\tilde{\z}+\frac{1}{4\z_{1}^{2}}A(A-1)}{\tilde{\z}^{2}(\tilde{\z}-1)(\tilde{\z}-1/\z_{1}^{2})}y=0\,.
\end{equation}
The four regular singularities $\{0, \z_{1}, \z_{2}, \infty\}$ of \eqref{e:HeunODE} are mapped into $\{0, 1, 1/\z_{1}^{2}, \infty\}$.
Moreover, applying the change of dependent variable $y(\z)=\z^{\rho}\tilde{y}(\z)$ to \eqref{e:generic Heun}) yields
\begin{equation}
\label{e:generic Heun'}
\tilde{y}_{\z\z} + P(\z)\tilde{y}_{\z} + Q(\z)\tilde{y} = 0\,,
\end{equation}
where 
\[
P(\z)=\frac{\gamma+2\rho}{\z}+\frac{\delta}{\z-1}+\frac{\epsilon}{\z-a}\,,\qquad 
Q(\z)=\frac{\rho(\rho-1+\gamma)}{\z^{2}}+\frac{\delta\rho}{\z(\z-1)}+\frac{\epsilon \rho}{\z(\z-a)}+\frac{\alpha\beta \z-\xi}{\z(\z-1)(\z-a)}\,.
\]
Note that \eqref{e:Heun'}) and \eqref{e:generic Heun'} are of the same form with 
$a=1/\z_{1}^{2}$ and 
$\z_{1} = [m-2 + 2\sqrt{1-m}]/m$. 
Reducing to a common denominator for $P(\z)$ and comparing the corresponding coefficients with \eqref{e:Heun'} leads to
$\gamma+2\rho=1/2$, $\delta=1/2$ and $\epsilon=1/2$, 
which implies that $l_{2}=l_{3}=0$.
Repeating the same procedure for $Q(\z)$, we find that: 
(i)~$\rho=(A-1)/2$ or $\rho=-A/2$, and; 
(ii)~$-A(A+1)/4=-\rho(\rho-1/2)+\alpha\beta$. 

Now we discuss the two possible cases for $\rho$:
If $\rho=(A-1)/2$, then $\gamma=3/2-A$ and $\alpha\beta=1/2-A$. 
Combining $\alpha+\beta=\gamma+\delta+\epsilon-1=\gamma$, one obtains $l_{1}=A-1$ and $l_{0}=-1-A ~\text{or}~A$. 
Alternatively, if $\rho=-A/2$, then $\gamma=1/2+A$, $\alpha\beta=0$, and $\alpha+\beta=1/2+A$. It turns out that $l_{1}=-A$ and $l_{0}=-1-A~\text{or}~ A$.
Either way, we therefore have that $l_{0}$, $l_{1} \in \Z$ if and only if $A \in \Z$.

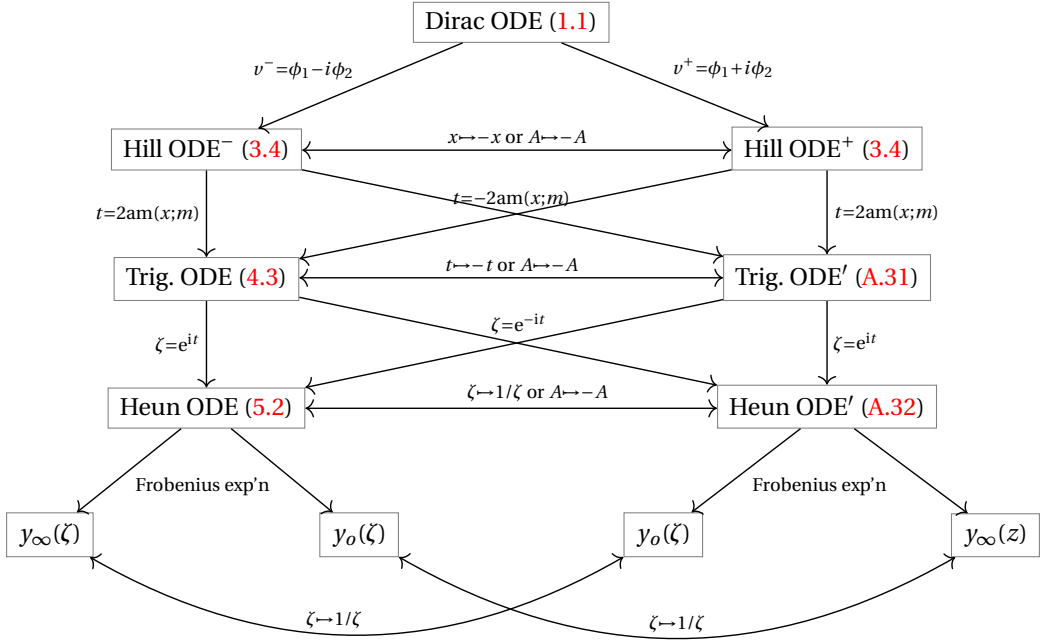
\begin{figure}[t!]
\begingroup
\small
\vspace*{-1ex}
\be
\begin{tikzcd}[row sep=huge,column sep=tiny,cells={nodes={draw=gray}}]
& & & \text{Dirac ODE}~\eqref{e:Diraceigenvalueproblem} \arrow[dll,"v^- = \phi_1 - i\phi_2"'] \arrow[drr,"v^+ = \phi_1 + i\phi_2"] & & &
\\
& \text{Hill ODE${}^-$}~\eqref{e:hillevalue} \arrow[d,"t=2\am(x;m)"'] \arrow[drrrr, "\hspace{-9mm} t=-2\am(x;m)"''] \arrow[rrrr, "x\mapsto-x \,\,\text{or}\,\,A\mapsto-A",leftrightarrow] & & & & \text{Hill ODE${}^+$}~\eqref{e:hillevalue} \arrow[d,"t=2\am(x;m)"] \arrow[dllll] & 
\\
& \text{Trig.\ ODE}~\eqref{e:trigvalue} \arrow[d, "\zeta=\e^{\i t}"'] \arrow[drrrr, "\kern-1em\zeta = \e^{-\i t}"] & & & & \text{Trig.\  ODE}'~\eqref{e:trigonometricODE'} \arrow[d, "\zeta = \e^{\i t}"] \arrow[dllll] \arrow[llll, "t\mapsto-t \,\,\text{or}\,\,A\mapsto-A"', leftrightarrow] & 
\\
& \text{Heun ODE}~\eqref{e:HeunODE} \arrow[dl, "\text{Frobenius~exp'n}"] \arrow[dr] & & & & \text{Heun ODE}'~\eqref{e:HeunODE'} \arrow[llll, "\hspace*{7mm}\zeta\mapsto 1/\zeta \,\,\text{or}\,\, A\mapsto-A"', leftrightarrow] \arrow[dl, "\text{Frobenius~exp'n}"] \arrow[dr] & 
\\
y_\infty(\z) \arrow[rrrr, bend right, leftrightarrow, "\hspace{-6mm}\zeta\mapsto 1/\zeta" ] & & y_o(\z) \arrow[rrrr, bend right, leftrightarrow,  "\zeta\mapsto 1/\zeta" ] & & y_o(\z) & & y_\infty(z)
\end{tikzcd}
\nonumber
\ee
\endgroup
\vglue-\bigskipamount
\caption{Relations between the various ODEs and solutions discussed throughout this work.
}
\label{f:commutative}
\end{figure}

\subsection{Transformations $A\mapsto-A$ and $\zeta\mapsto1/\zeta$}
\label{a:FourierHeun_connection_Anoninteger}

The maps $A\mapsto -A$ and $\zeta\mapsto1/\zeta$ allow one to establish a connection between several related objects.
Specifically, using the change of independent variable~(\ref{e:amp}), 
Hill's equation $H^{+}v^{+} = \lambda v^{+}$ is mapped into the following second-order ODE with trigonometric coefficients
\be
4(1-m\sin^2 \halft)y_{tt} - (m\sin t)y_{t} + 
(\lambda + A^2(1-m\sin^2 \halft) - \halfi Am\sin t)y = 0\,.
\label{e:trigonometricODE'}
\ee
Next, applying the transformation $\z=\e^{\i t}$ to~\eqref{e:trigonometricODE'} yields another Heun ODE,
namely,
\be
\z^2F(\z;m)y_{\z\z} + \z G(\z;m)y_{\z} + \tilde{H}(\z;\lambda,A,m)y = 0\,,
\label{e:HeunODE'}
\ee
where $F(\z;m)$ and $G(\z;m)$ are still given by~\eqref{e:p1} and~\eqref{e:p2}, respectively, and with
$\tilde{H} := H(\z;\l,-A,m)$.
Note that the four regular singular points of~\eqref{e:HeunODE}) and~\eqref{e:HeunODE'} are the same.
The full chain of transformations and correspondences is summarized in the commutative diagram in Fig.~\ref{f:commutative}.

\catcode`\@ 11

\end{document}